%% file: paper_knaf_selder_spindler.tex
\hsize = 6.5 in 
\vsize = 9.1 in
\overfullrule 0pt
\magnification = \magstep1 
\font\titlefont = cmssdc10 at 16 pt 
\font\titlefontsmall = cmssdc10 at 12 pt 
\font\tenmsb = msbm10
\font\sevenmsb = msbm7
\font\fivemsb = msbm5

\newfam\msbfam
\textfont\msbfam=\tenmsb
\scriptfont\msbfam=\sevenmsb
\scriptscriptfont\msbfam=\fivemsb
\def\N{\fam\msbfam N}
\def\Z{\fam\msbfam Z} 
\def\Q{\fam\msbfam Q}
\def\R{\fam\msbfam R}

\def\P{\fam\msbfam P}
\input epsf.tex
\input colorbox.tex

\input bfall.tex

\def\doppelbildbox#1#2#3#4{\smallskip\noindent \vbox{\hrule\hbox{\noindent \vrule
   \vbox{\noindent\medskip\vbox{$$\epsfxsize #1\epsffile{#2}\qquad 
           \epsfxsize #1\epsffile{#3}$$}\medskip}\vrule }
   \hrule\smallskip\hbox{\vbox{\noindent #4}}}\smallskip}
\hoffset -0.3 true cm
%
%
\leftline{\titlefontsmall Hagen Knaf, Erich Selder, Karlheinz Spindler:} 
\bigskip\noindent 
\rightline{\titlefont Explicit transformation of an intersection of two}
\rightline{\titlefont quadrics to an elliptic curve in Weierstra\ss{} form}
\bigskip 

\leftline{\titlefont 1. Introduction} \medskip 
This paper is motivated by problems in Diophantine analysis which can be formulated as 
problems of finding rational points on the intersection of two quadrics. Let us look at 
two typical examples. \medskip 

$\bullet$ {\bf Example 1: Euler's problem of concordant forms.} Let $M,N\in{\Z}$ be two different 
nonzero integers. Euler called $M$ and $N$ concordant or discordant depending on whether or not the 
system of the equations $X^2+MY^2=Z^2$ and $X^2+NY^2=W^2$ possesses a nontrivial solution $(X,Y,Z,W)
\in{\Z}^4$ with $Y\not= 0$ (see [3], [7]). This amounts to the question whether or not the intersection 
$Q_{M,N}$ of the two quadrics $X^2+MY^2=Z^2$ and $X^2+NY^2=W^2$ possesses a rational point which is 
nontrivial in the sense that $Y\not= 0$. \medskip 

$\bullet$ {\bf Example 2: Rational squares in arithmetic progression.} We can ask when four rational 
numbers $\alpha,\beta,\gamma,\delta\in{\Q}$ form a progression $\alpha^2<\beta^2<\gamma^2<
\delta^2$ which is part of an arithmetic progression, which is the case if and only if differences 
between subsequent terms are integral multiples of a (necessarily rational) step size $s$, say 
$$\beta^2-\alpha^2=ks,\qquad \gamma^2-\beta^2=\ell s,\qquad 
  \delta^2-\gamma^2=ms\leqno{(1)}$$ 
where $k,\ell,m\in{\N}$. If (1) holds for a fixed triplet $(k,\ell,m)$, then $\ell (\beta^2-\alpha^2) 
= \ell ks = k\ell s = k(\gamma^2-\beta^2)$ and $m(\gamma^2-\beta^2) = m\ell s = \ell ms = \ell 
(\delta^2-\gamma^2)$ so that 
$$(k+\ell)\beta^2\, -k\gamma^2\, -\, \ell\alpha^2\ =\ 0\qquad\hbox{and}\qquad 
  -m\beta^2\, +\, (m+\ell)\gamma^2\, -\,\ell\delta^2\ =\ 0.\leqno{(2)}$$ 
Conversely, if $\alpha,\beta,\gamma,\delta$ are such that (2) holds for a given 
triplet $(k,\ell, m)\in{\N}^3$, we can define
$$s\ :=\ {{\beta^2-\alpha^2}\over k}\ =\ {{\gamma^2-\beta^2}\over{\ell}}\ =
  \ {{\delta^2-\gamma^2}\over{m}}\leqno{(3)}$$ 
to get a solution of the original problem (1). Thus given a triplet $(k,\ell, m)$, we ask 
whether or not the system (2) admits a nontrivial solution $(\alpha,\beta,\gamma,\delta)\not= 
\lambda\cdot (\pm 1,\pm 1,\pm 1,\pm 1)$ in rational numbers. Now since the 
equations in (2) are homogeneous, we can interpret $\alpha,\beta,\gamma,\delta$ 
as projective coordinates of a point in ${\P}^3({\Q})$, and condition (2) can be 
reformulated by stating that the point $(X_0,X_1,X_2,X_3):=(\beta,\gamma,\alpha,
\delta)$ is a rational point in the intersection of the two quadrics 
$$\eqalign{
  Q_1(k,\ell,m)\ &:=\ \{ (X_0,X_1,X_2,X_3)\in{\P}^3({\R})\mid (k+\ell) X_0^2-kX_1^2-\ell X_2^2\ =\ 0\},\cr 
  Q_2(k,\ell,m)\ &:=\ \{ (X_0,X_1,X_2,X_3)\in{\P}^3({\R})\mid -mX_0^2+(m+\ell)X_1^2-\ell X_3^2\ =\ 0\}.\cr}
  \leqno{(4)}$$ 
\eject 
 
The purpose of this paper is to present, in detail, the construction of a rationally defined 
isomorphism (biregular mapping) between a rationally defined smooth intersection of two quadrics 
in projective three-space ${\P}^3(K)$ where $K$ is a (not yet specified) field and an elliptic 
curve in Weierstra\ss{} form which maps a distinguished rational point to the point at infinity. 
The existence of such an isomorphism makes available the elaborate theory of elliptic curves 
to study number-theoretical problems of the type described before. For example, the 
construction will yield a rationally defined isomorphism between the intersection 
$Q_{k,\ell ,m}:=Q_1(k,\ell,m)\cap Q_2(k,\ell,m)$ and the plane Weierstra\ss{} cubic 
$E_{k,\ell ,m}$ given by the affine equation $y^2=x(x+km)\bigl( x+(k\! +\!\ell)(\ell\! +\! 
m)\bigr)$. It turns out that, up to the very last step, the construction involved does 
not depend on any special characteristics of the equations (2), but works quite generally, 
and the purpose of this paper is to describe this construction in full generality and in 
full detail. (The conclusions to be drawn from this construction for the original 
number-theoretical problem will be discussed elsewhere.) \smallskip 

Most of the ideas and calculations can be found scattered in the literature, albeit 
in some cases only in a sketchy way, only illustrated by way of example or without 
considering all special cases which can occur. (References will be given as we 
progress.) While being only interested in the base-field ${\Q}$, we note that most of 
the calculations are valid for an arbitrary base-field $K$, sometimes with the 
restriction that the characteristic be different from $2$ or $3$. The calculations 
consist of several steps, which can be summarized as follows: \smallskip 
\item{$\bullet$} transformation of the intersection of two quadrics to a smooth plane 
  cubic curve; 
\item{$\bullet$} transformation of a rationally defined smooth plane cubic curve to
  an elliptic curve given by a Weierstra\ss{} equation; 
\item{$\bullet$} in the special situation of $Q_{k,l,m}$ transformation of the 
  Weierstra\ss{} equation to the equation $y^2 = x(x+km)\bigl( x+(k\! +\!\ell)
  (\ell\! +\! m)\bigr)$.
\smallskip\noindent 
The calculations will be illustrated by a series of diagrams showing the curves that 
occur during the procedure, both in an affine as well as in a projective coordinate 
system. The relation between projective coordinates $(X,Y,Z)$ and affine coordinates 
$(x,y)$ is given by $x=X/Z$ and $y=Y/Z$. The projective coordinates are represented 
by a triangle whose vertices are the projective points $(0,0,1)$ (leftmost point), 
$(1,0,0)$ (rightmost point) and $(0,1,0)$ (upper middle point). In each case the 
drawings for the curves are based on the true equations; however, suitable scaling 
factors were used in order to make the relevant geometric properties of the curves 
visible. In addition to the curves themselves, we show those points which are relevant 
for the actual transformation. \smallskip\noindent 
Several times throughout the paper the following result will be invoked:
\bigskip\noindent 
{\bf Regularity Theorem:} \it Let $C$ be a nonsingular curve, $V\subseteq{\P}^n(K)$ a projective 
variety and $\varphi:C\rightarrow V$ a rational mapping. Then $\varphi$ is a morphism, i. e., 
a mapping which is regular at any point of $C$.
\bigskip\noindent\rm 
This is a standard result in algebraic geometry which can be found, for example, in 
[4] (see Chap. I, \S 6, Prop. 6.8), in [10] (see Chap. II, Sect. 3, Thm. 3, Cor. 1), or 
in [11] (see Chap. II, \S 2, Prop. 2.1).
\vfill\eject 

\leftline{\titlefont 2. Transformation of a quadric intersection to a plane cubic} \medskip 

The transformation of a smooth intersection of two quadrics in projective three-space 
to a smooth plane cubic curve is briefly sketched in [1], Chap. 8, pp. 36; 
some special cases are discussed in [2], Sect. 1.4.3, pp. 123-125. The most 
interesting part of the following calculations is the inverse mapping, which is not 
to be obtained in a completely trivial way. \smallskip 

The starting point is as follows: We are given two quadrics $Q_1$ and $Q_2$ in 
projective three-space whose intersection is a smooth irreducible curve over an algebraic 
closure of the base-field $K$, and we are given a $K$-rational point $x=(x_0,x_1,x_2,x_3)\in 
Q_1\cap Q_2$. (The condition on $Q_1$ and $Q_2$ is satisfied, for example, if $Q_1$ and 
$Q_2$ generate a separable pencil in the sense of [8], Sect. 2.4.3, p. 74; also see [9], 
chap. XIII.) The quadrics can be assumed to be given as  
$$Q_1=\{ X\in{\P}^3(K)\mid X^TAX=0\}\quad\hbox{and}\quad 
  Q_2=\{ X\in{\P}^3(K)\mid X^TBX=0\}\leqno{(5)}$$ 
with two symmetric $(4\times 4)$-matrices $A$ and $B$. (Each quadric can be written in this 
way over a base-field whose characteristic is different from $2$.) Moreover, after re-indexing 
variables if necessary, we may assume that $x_3\not= 0$ and hence even $x_3=1$. This will be 
assumed from now on. \smallskip 

To rewrite (5), we consider the coordinate transformation 
$$\left[\matrix{Y_0\cr Y_1\cr Y_2\cr Y_3\cr}\right] = 
  \left[\matrix{1&0&0&-x_0\cr 0&1&0&-x_1\cr 0&0&1&-x_2\cr 0&0&0&\phantom{-}1\cr}
  \right]\left[\matrix{X_0\cr X_1\cr X_2\cr X_3\cr}\right]\! ,
  \quad\hbox{i.e.},\quad  
  \left[\matrix{X_0\cr X_1\cr X_2\cr X_3\cr}\right] = 
  \left[\matrix{1&0&0&x_0\cr 0&1&0&x_1\cr 0&0&1&x_2\cr 0&0&0&1\cr}\right]
  \left[\matrix{Y_0\cr Y_1\cr Y_2\cr Y_3\cr}\right],\leqno{(6)}$$ 
which we write for short as $Y=PX$ and $X=QY$ where $Q=P^{-1}$. $\bigl($Geometrically, 
the transition from $X=(X_0,X_1,X_2,X_3)$ to $Y=(Y_0,Y_1,Y_2,Y_3)$ is the translation 
which maps the point $x$ to the point $(0,0,0,1)$.$\bigr)$ The equations of the 
quadrics $Q_1$ and $Q_2$ are then transformed to
$$\eqalign{
  &0\ =\ X^TAX\ =\ Y^TQ^TAQY\ =\ Y^T\!\left[\matrix{a_{00}&a_{01}&a_{02}&u_0\cr  
     a_{10}&a_{11}&a_{12}&u_1\cr a_{20}&a_{21}&a_{22}&u_2\cr u_0&u_1&u_2&0\cr}
     \right]\! Y,\cr 
  &0\ =\ X^TBX\ =\ Y^TQ^TBQY\ =\ Y^T\!\left[\matrix{b_{00}&b_{01}&b_{02}&v_0\cr  
     b_{10}&b_{11}&b_{12}&v_1\cr b_{20}&b_{21}&b_{22}&v_2\cr v_0&v_1&v_2&0\cr}
     \right]\! Y\cr} 
\leqno{(7)}$$ 
where 
$$u_i:=\sum_{j=0}^3 a_{ij}x_j\quad\hbox{and}\quad v_i:=\sum_{j=0}^3 b_{ij}x_j
  \quad\hbox{for}\ 0\leq i\leq 2. \leqno{(8)}$$ 
The main point of this step is that in the coefficient matrices occurring in $(7)$ 
the entry at the bottom right vanishes, which means that in the quadrics defined by 
$(7)$ the variable $Y_3$ occurs only linearly. Thus after the transformation we 
have obtained the intersection of two quadrics of the form 
$$\eqalign{
  \widehat{Q}_1\ &=\ \{ (Y_0,Y_1,Y_2,Y_3)\in{\P}^3(K)\mid q_1(Y_0,Y_1,Y_2) 
     + \ell_1(Y_0,Y_1,Y_2)Y_3= 0\},\cr 
  \widehat{Q}_2\ &=\ \{ (Y_0,Y_1,Y_2,Y_3)\in{\P}^3(K)\mid q_2(Y_0,Y_1,Y_2) 
     + \ell_2(Y_0,Y_1,Y_2)Y_3= 0\}\cr}\leqno{(9)}$$ 
where $q_1,q_2$ are quadratic and $\ell_1,\ell_2$ are linear. Specifically, we 
have 
$$\eqalign{
  q_1(Y_0,Y_1,Y_2)\ &=\ \hbox{$\sum_{i,j=0}^2 a_{ij}Y_iY_j$},\cr 
  q_2(Y_0,Y_1,Y_2)\ &=\ \hbox{$\sum_{i,j=0}^2 b_{ij}Y_iY_j$},\cr 
  \ell_1(Y_0,Y_1,Y_2)\ &=\ 2(u_0Y_0+u_1Y_1+u_2Y_2),\cr
  \ell_2(Y_0,Y_1,Y_2)\ &=\ 2(v_0Y_0+v_1Y_1+v_2Y_2).\cr}\leqno{(10)}$$ 
Moreover, the rational point $x$ is mapped to the rational point $(0,0,0,1)$. 
Consequently, on $\widehat{Q}_1\cap\widehat{Q}_2$ we have $q_1\ell_2 = -\ell_1
\ell_2Y_3 = -\ell_2\ell_1Y_3 = q_2\ell_1$. Hence if we define 
$$(11)\quad C:=\{ (Y_0,Y_1,Y_2)\!\in\!{\P}^2(K)\mid q_1(Y_0,Y_1,Y_2)\ell_2(Y_0,Y_1,Y_2) 
   = q_2(Y_0,Y_1,Y_2)\ell_1(Y_0,Y_1,Y_2)\}$$ 
(which is a smooth cubic in the projective plane) then the assignment $(Y_0,Y_1,Y_2,Y_3)
\mapsto (Y_0,Y_1,Y_2)$ (which is everywhere defined except at $(0,0,0,1)$) maps 
$\widehat{Q}_1\cap\widehat{Q}_2$ into $C$. As we shall see in a moment, this 
mapping can be redefined around the point $(0,0,0,1)$ to yield an everywhere defined 
regular mapping 
$$\varphi:\widehat{Q}_1\cap\widehat{Q}_2\rightarrow C\leqno{(12)}$$ 
which is then automatically an isomorphism (due to the regularity theorem quoted in the 
introduction), a fact which in our situation can also be established by explicitly 
revealing the inverse mapping 
$$\psi:C\rightarrow\widehat{Q}_1\cap\widehat{Q}_2\leqno{(13)}$$ 
rather than by invoking a general principle. Let $(Y_0,Y_1,Y_2)\in C$. If 
$\ell_1(Y_0,Y_1,Y_2)\not= 0$ then 
$$\psi (Y_0,Y_1,Y_2)\ :=\ \left[\matrix{\ell_1(Y_0,\! Y_1,\! Y_2)Y_0\cr 
  \ell_1(Y_0,\! Y_1,\! Y_2)Y_1\cr \ell_1(Y_0,\! Y_1,\! Y_2)Y_2\cr -q_1(Y_0,\! Y_1,\! Y_2)\cr}
  \right] \leqno{(14)}$$ 
lies in $\widehat{Q}_1\cap\widehat{Q}_2$, and $(\varphi\circ\psi)(Y_0,Y_1,Y_2) = 
(Y_0,Y_1,Y_2)$. Analogously, if $\ell_2(Y_0,Y_1,Y_2)\not= 0$ then   
$$\psi (Y_0,Y_1,Y_2)\ :=\ \left[\matrix{\ell_2(Y_0,\! Y_1,\! Y_2)Y_0\cr 
  \ell_2(Y_0,\! Y_1,\! Y_2)Y_1\cr \ell_2(Y_0,\! Y_1,\! Y_2)Y_2\cr -q_2(Y_0,\! Y_1,\! Y_2)\cr}
  \right] \leqno{(15)}$$ 
lies in $\widehat{Q}_1\cap\widehat{Q}_2$, and $(\varphi\circ\psi)(Y_0,Y_1,Y_2) = 
(Y_0,Y_1,Y_2)$. What if $\ell_1(Y_0,Y_1,Y_2) = \ell_2(Y_0,Y_1,Y_2) = 0$ (in which case 
$(Y_0,Y_1,Y_2)$ clearly also lies in $C$)? In this case $(Y_0,Y_1,Y_2)$ 
is the unique point of intersection of the lines $u_0Y_0+u_1Y_1+u_2Y_2=0$ and 
$v_0Y_0+v_1Y_1+v_2Y_2=0$, which is the (rational) point 
$$(u_1v_2-u_2v_1,\, u_2v_0-u_0v_2,\, u_0v_1-u_1v_0)\ =:\ z.\leqno{(16)}$$ 
We could again invoke the regularity theorem quoted in the introduction to conclude that $\varphi$ 
and $\psi$ can be extended to become mutually inverse regular mappings in such a way that $\varphi 
(0,0,0,1) = z$ and $\psi (z)=(0,0,0,1)$, but it is also possible to see this in an elementary way 
by explicitly redefining $\varphi$ around the point $(0,0,0,1)$. By another elementary argument it 
can be seen that $\psi$ does not even need to be redefined around the point $z$, but that one of the 
two representations (14) and (15) is necessarily defined at this point. Let us now redefine $\varphi$ 
around the point $(0,0,0,1)$. To begin with, we split each of the polynomials defining 
$\widehat{Q}_1$ and $\widehat{Q}_2$ into the part containing $Y_2$ and the part 
not containing $Y_2$, thereby writing 
$$\eqalign{q_1+\ell_1Y_3\ &=\ \alpha_1Y_2 + \alpha_2,\cr 
           q_2+\ell_2Y_3\ &=\ \beta_1Y_2 + \beta_2\cr}\leqno{(17)}$$ 
where $\alpha_1,\beta_1$ are linear polynomials in $(Y_0,Y_1,Y_2,Y_3)$ and where 
$\alpha_2,\beta_2$ are quadratic polynomials in $(Y_0,Y_1,Y_3)$. Specifically, 
we have 
$$\eqalign{
  \alpha_1\ &=\ 2a_{02}Y_0 + 2a_{12}Y_1 + a_{22}Y_2 + 2u_2Y_3,\cr 
  \alpha_2\ &=\ a_{00}Y_0^2 + a_{11}Y_1^2 + 2a_{01}Y_0Y_1 + 2u_0Y_0Y_3 + 2u_1Y_1Y_3,\cr 
  \beta_1\ &=\ 2b_{02}Y_0 + 2b_{12}Y_1 + b_{22}Y_2 + 2v_2Y_3,\cr 
  \beta_2\ &=\ b_{00}Y_0^2 + b_{11}Y_1^2 + 2b_{01}Y_0Y_1 + 2v_0Y_0Y_3 + 2v_1Y_1Y_3.\cr}
  \leqno{(18)}$$ 
Next, we  split each of the polynomials $\alpha_2$ and $\beta_2$ into a part containing 
$Y_0$ and a part containing $Y_1$ (where the terms with $Y_0Y_1$ can be arbitrarily 
assigned to either term), thereby writing 
$$\eqalign{\alpha_2\ &=\ \gamma_0Y_0\, +\,\gamma_1Y_1,\cr 
           \beta_2\ &=\ \delta_0Y_0\, +\,\delta_1Y_1\cr}\leqno{(19)}$$ 
where $\gamma_0,\gamma_1,\delta_0,\delta_1$ are all linear in $(Y_0,Y_1,Y_2,Y_3)$. 
A specific choice is given by 
$$\eqalign{
  \gamma_0\ &=\ a_{00}Y_0 + a_{01}Y_1 + 2u_0Y_3, \cr 
  \gamma_1\ &=\ a_{11}Y_1 + a_{01}Y_0 + 2u_1Y_3, \cr 
  \delta_0\ &=\ b_{00}Y_0 + b_{01}Y_1 + 2v_0Y_3, \cr 
  \delta_1\ &=\ b_{11}Y_1 + b_{01}Y_0 + 2v_1Y_3.\cr}
\leqno{(20)}$$ 
Now on $\widehat{Q}_1\cap\widehat{Q}_2$ we have $q_1+\ell_1Y_3 = q_2+\ell_2Y_3 = 0$, 
hence $\alpha_1Y_2=-\alpha_2$ and $\beta_1Y_2 = -\beta_2$, hence $\alpha_1\beta_2 = 
-\alpha_1\beta_1Y_2 = -\beta_1\alpha_1Y_2 = \beta_1\alpha_2$ and therefore 
$\alpha_1\delta_0Y_0 + \alpha_1\delta_1Y_1 = \beta_1\gamma_0Y_0 + \beta_1\gamma_1Y_1$, 
i.e., 
$$(\alpha_1\delta_1-\beta_1\gamma_1)Y_1\ =\ (\beta_1\gamma_0-\alpha_1\delta_0)Y_0.
  \leqno{(21)}$$ 
Thus the polynomial 
$$\lambda\ :=\ \alpha_1\beta_1(\alpha_1\delta_1-\beta_1\gamma_1)\leqno{(22)}$$ 
satisfies on $\widehat{Q}_1\cap\widehat{Q}_2$ the equations  
$$\lambda Y_1\ =\ \alpha_1\beta_1(\beta_1\gamma_0-\alpha_1\delta_0)Y_0\leqno{(23)}$$ 
and 
$$\eqalign{
  \lambda Y_2\ &=\ (\alpha_1\delta_1-\beta_1\gamma_1)\alpha_1\beta_1Y_2\cr 
       &=\ -(\alpha_1\delta_1-\beta_1\gamma_1)\alpha_1\beta_2\cr 
       &=\ -(\alpha_1\delta_1-\beta_1\gamma_1)\alpha_1(\delta_0Y_0 + \delta_1Y_1)\cr 
       &=\ -(\alpha_1\delta_1-\beta_1\gamma_1)\alpha_1\delta_0Y_0 
                -\alpha_1\delta_1(\alpha_1\delta_1-\beta_1\gamma_1)Y_1\cr 
       &=\ -(\alpha_1\delta_1-\beta_1\gamma_1)\alpha_1\delta_0Y_0 
                -\alpha_1\delta_1(\beta_1\gamma_0-\alpha_1\delta_0)Y_0\cr  
       &=\ \alpha_1\beta_1(\gamma_1\delta_0-\delta_1\gamma_0)Y_0\cr}\leqno{(24)}$$ 
so that  
$$\left[\matrix{Y_0\cr Y_1\cr Y_2\cr}\right]\ =\ \left[\matrix{\lambda Y_0\cr 
  \lambda Y_1\cr \lambda Y_2\cr}\right]\ =\ \left[\matrix{
     \alpha_1\beta_1(\alpha_1\delta_1-\beta_1\gamma_1)Y_0\cr 
     \alpha_1\beta_1(\beta_1\gamma_0-\alpha_1\delta_0)Y_0\cr 
     \alpha_1\beta_1(\gamma_1\delta_0-\delta_1\gamma_0)Y_0\cr}\right]
\, .\leqno{(25)}$$ 
Discarding common factors, the mapping $\varphi$ can thus be rewritten as 
$$\varphi(Y_0,Y_1,Y_2,Y_3)\ =\ \left[\matrix{
     \alpha_1\delta_1-\beta_1\gamma_1\cr 
     \beta_1\gamma_0-\alpha_1\delta_0\cr 
     \gamma_1\delta_0-\delta_1\gamma_0\cr}\right]
\, .\leqno{(26)}$$ 
Written in this way, the mapping $\varphi$ is defined at the point $(0,0,0,1)$; namely, 
plugging in $(0,0,0,1)$ for $(Y_0,Y_1,Y_2,Y_3)$ in (18) and (20), we find that 
$$\varphi (0,0,0,1) = \left[\matrix{2u_2\cdot 2v_1-2v_2\cdot 2u_1\cr 
  2v_2\cdot 2u_0-2u_2\cdot 2v_0\cr 2u_1\cdot 2v_0-2v_1\cdot 2u_0\cr}\right]
  = \left[\matrix{-4(u_1v_2\! -\! u_2v_1)\cr -4(u_2v_0\! -\! u_0v_2)\cr -4(u_0v_1\! 
  -\! u_1v_0)\cr}\right] = \left[\matrix{u_1v_2\! -\! u_2v_1\cr u_2v_0\! -\! u_0v_2\cr 
  u_0v_1\! -\! u_1v_0\cr}\right]\leqno{(27)}$$ 
which shows that $\varphi$ maps $(0,0,0,1)$ to $z$. Thus $\varphi$ provides an isomorphism 
from the quadric intersection $\widehat{Q}_1\cap\widehat{Q}_2$ onto a smooth cubic curve. 
We want to verify the regularity of the inverse mapping $\psi$ at the point $z$ in an 
elementary way (without invoking the regularity theorem quoted in the introduction). We have the 
two representations (14) and (15) for $\psi$, and it is immediately clear that at least one of 
these is well-defined at any point of $C$ other than $z$. However, we claim that at 
least one of these must be also defined at the point $z$ (so that no redefinition of $\psi$ 
is required about this point). The point $z$ is uniquely defined by the conditions $\ell_1(z)=0$ 
and $\ell_2(z)=0$; we need to rule out that, in addition, the conditions $q_1(z)=0$ and 
$q_2(z)=0$ can simultaneously hold. Now if this were the case then $(z_1,z_2,z_3,Y_3)$ would 
be contained in $\widehat{Q}_1\cap\widehat{Q}_2$ for all $Y_3$, as is clear from equation 
$(9)$, which means that $\widehat{Q}_1\cap\widehat{Q}_2$ would contain a whole line over an 
algebraically closed field containing $K$, contradicting our assumption that $Q_1\cap Q_2$ and 
hence $\widehat{Q}_1\cap\widehat{Q}_2$ is irreducible.\hfill\break\smallskip 
 
{\bf Example.} We consider the quadrics $Q_1$ and $Q_2$ which are defined as the solution 
set of the equations 
$$X_0^2 + 2X_0X_1 + 2X_1^2 - 6X_1X_2 - 2X_2X_3 + 3X_3^2\ =\ 0$$ 
and 
$$-2X_0^2 + X_1^2 + 2X_2^2 - X_3^2\ =\ 0,$$ 
respectively. Note that $Q_1\cap Q_2$ contains the rational point $(1,1,1,1)$. The isomorphism 
described before maps $Q_1\cap Q_2$ to the smooth cubic $C$ given by the equation 
$$-2Y_0^3 + 3Y_0^2Y_1 + 6Y_0^2Y_2 + 4Y_0Y_1^2 - 16Y_0Y_1Y_2 + 4Y_0Y_2^2 - 2Y_1^3 - 2Y_1^2Y_2 + 
  12Y_1Y_2^2 - 8Y_2^3\ =\ 0$$ 
so that the coefficients of $C$ are given by 
$$\eqalign{
  &C_{300}=-2,\ C_{210}=3,\ C_{201}=6,\ C_{120}=4,\ C_{111}=-16,\cr 
  &C_{102}=4,\ C_{030}=-2,\ C_{021}=-2,\ C_{012}=12,\ C_{003}=-8.\cr}$$ 
The above isomorphism maps $(1,1,1,1)$ to the rational point $(2,2,1)$. 
\smallskip 
\doppelbildbox{7 true cm}{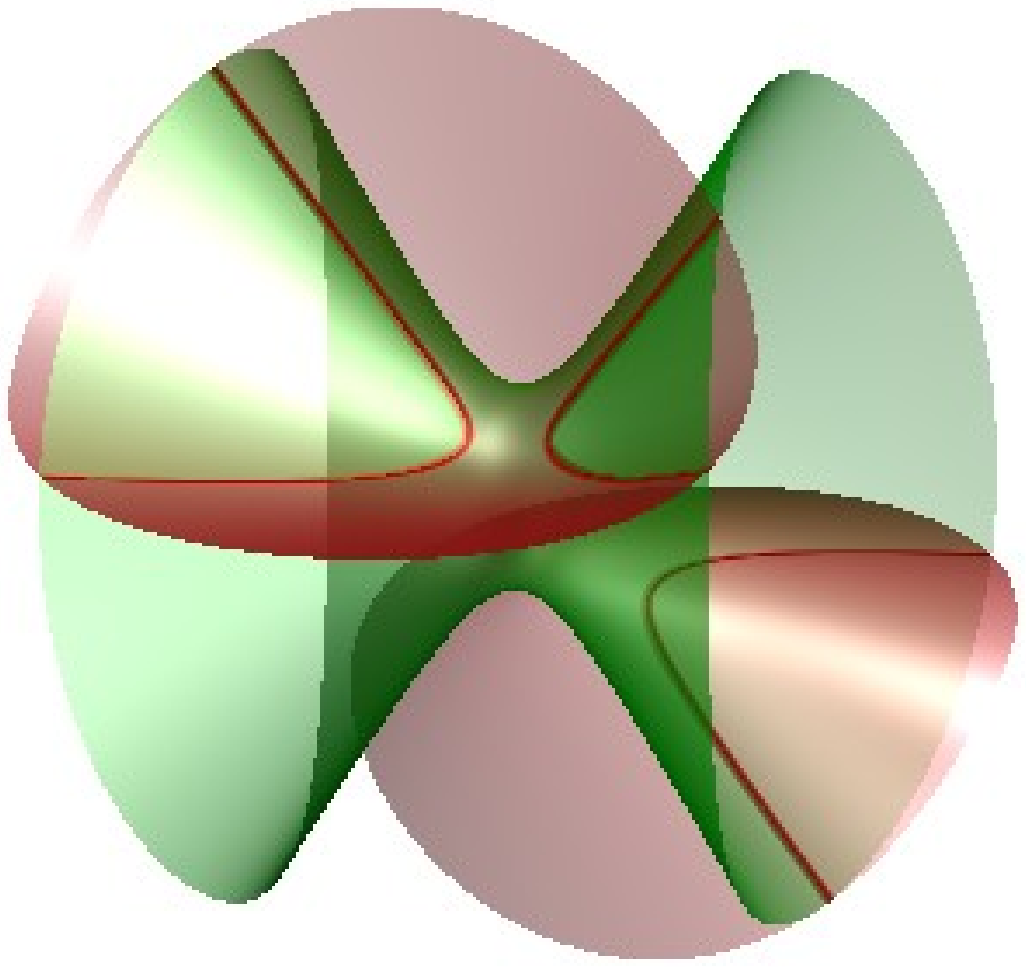}{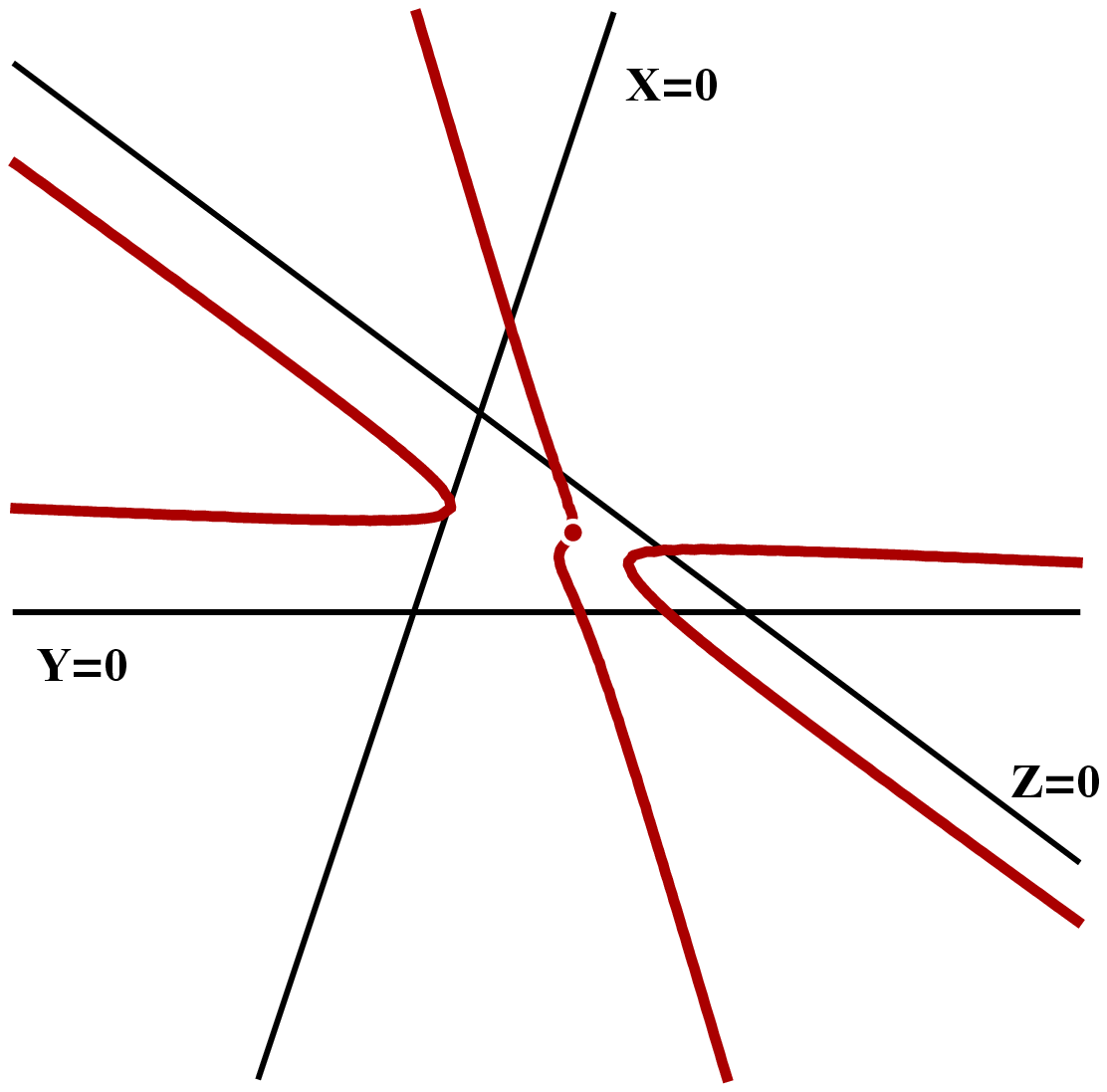} 
{{\bf Left:} Affine view $X_3=1$ of $Q_1\cap Q_2$. {\bf Right:} Projective view $Z=1$ of 
the cubic $C$ (with the distinguished rational point $(2,2,1)$ marked by a red dot).}
\vfill\eject 

\leftline{\titlefont 3. Transformation of a smooth plane cubic to Weierstra\ss{} form} 
\medskip 

It is a well-known fact that any smooth cubic curve in the (projective) plane can be transformed 
into one given by a Weierstra\ss{} equation. In nearly any textbook on algebraic geometry or 
algebraic curves, this is done by choosing an inflection point and transforming this inflection 
point to the point at infinity (with respect to certain coordinates) in such a way that its 
tangent is given by the line at infinity. However, if the curve is defined over some field 
$K$ which is not algebraically closed (in our case $K={\Q}$), there may be no inflection 
point defined over $K$. By a long known but rarely used construction due to Nagell (see 
[6]) the transformation can be done in such a way that a given $K$-rational point, 
which is not necessarily an inflection point, is mapped to the point at infinity. \smallskip 

Nagell's construction seems to be half-forgotten, to the detriment of various 
approaches found in the literature to translate number-theoretical problems as 
problems of finding rational points on elliptic curves. (For example, in  [5], pp. 
6 ff., and in [7], the authors use a mapping of degree $4$ instead of an isomorphism, 
thus losing some information on effects corresponding to torsion points on the 
Weierstra\ss{} curve in question.) We suspect that one reason for Nagell's construction 
to sink into oblivion lies in the rather unwieldy calculations to which this 
construction leads. (This can be seen in the following examples by the fast-growing 
size of the coefficients showing up in the equations obtained; see also the explicit 
calculations in Appendix B of [12].) \smallskip 

The following calculations represent a complete and explicit version of Nagell's construction. 
However, our calculations follow a geometrically more convenient concept which is briefly 
sketched in [12], pp. 17 ff (also, see again the explicit example in Appendix B, pp. 
311 ff). Our calculations are complete in the sense that non-generic exceptional situations 
are covered, and they are explicit in the sense that all formulas may directly be realized
in any programming language. Let us quickly explain our notations. \smallskip 

Generally, we consider our cubics to be given by homogeneous equations in three variables $X,Y,Z$, 
i.e., by equations of the form $\sum_{i+j+k=3}\Gamma_{ijk}X^iY^jZ^k=0$ where $\Gamma_{ijk}\in K$ 
denotes the coefficient of the monomial $X^iY^jZ^k$. Since our cubics are considered to be 
($K$-)rationally defined, they contain a fixed ($K$-)rational point, which we denote by $p=
(p_x, p_y ,p_z)\in K^3$. The transformation from a general cubic $C$ in the plane to a cubic 
in Weierstra\ss{} form is accomplished by a sequence of transformations from one cubic to another,  
starting with the original cubic $C=C_{(0)}$ given in terms of the variables $(X,Y,Z) = (X_0,Y_0,
Z_0)$. The $r$-th step in our algorithm consists of a coordinate transformation (mostly linear, 
in one instance quadratic) transforming the coordinates $(X_{r-1},Y_{r-1},Z_{r-1})$ used for the 
old curve $C_{(r-1)}$ to the coordinates $(X_r,Y_r,Z_r)$ used for the new curve $C_{(r)}$, thereby 
mapping the ($K$-)rational point $p^{(r-1)}=(p_x^{(r-1)}, p_y^{(r-1)}, p_z^{(r-1)})$ to the 
($K$-)rational point $p^{(r)} = (p_x^{(r)},p_y^{(r)},p_z^{(r)})$. At any step in the calculations, 
all the data which are used for the calculations can be projectively simplified by cancelling 
common factors, which effects coefficients of cubic equations, coordinates of points in projective 
space and coordinates defining the slopes of lines. (For example, over the base-field ${\Q}$ we can 
always assume the rational point $p^{(r)}$ to have coprime integer coefficients.) We denote by 
$\Gamma^{(r)}_{ijk}$ the coefficients of the monomial $X_r^iY_r^jZ_r^k$ in the equation describing 
the curve $C_{(r)}$; i.e., the equation for $C_{(r)}$ is written as  
$$0\ =\ \sum_{i+j+k=3} \Gamma_{ijk}^{(r)}X^iY^jZ^k \leqno{(29)}$$ 
where $(X,Y,Z) = (X_r,Y_r,Z_r)$ are the coordinates used after the $r$-th transformation.
Without loss of generality we may assume that initially $p_x\neq 0$ and hence even that $p_x=1$. 
This assumption will be made now. Each of the various coordinate transformations will be 
explained in geometric terms before being written down explicitly, and subsequently the arithmetical 
effect of the transformation on the equations of the various cubics involved will be explained. 
To ease readability, each step is presented on a double-page with dia\-grams visualizing the 
transformations used on one page and the associated arithmetical explanations on the other page. 
\vfill\eject 

{\bf Step 1:} We start with an arbitrary cubic $C_{(0)}$ with a distinguished 
  rational point $p^{(0)} = (p_x,p_y,p_z) = (1,p_y,p_z)$. The first transformation 
  is the translation which maps $p^{(0)}$ to $p^{(1)} = (1,0,0)$ and hence yields 
  a cubic $C_{(1)}$ with the special distinguished point $p^{(1)} = (1,0,0)$.
  \doppelbildbox{7 true cm}{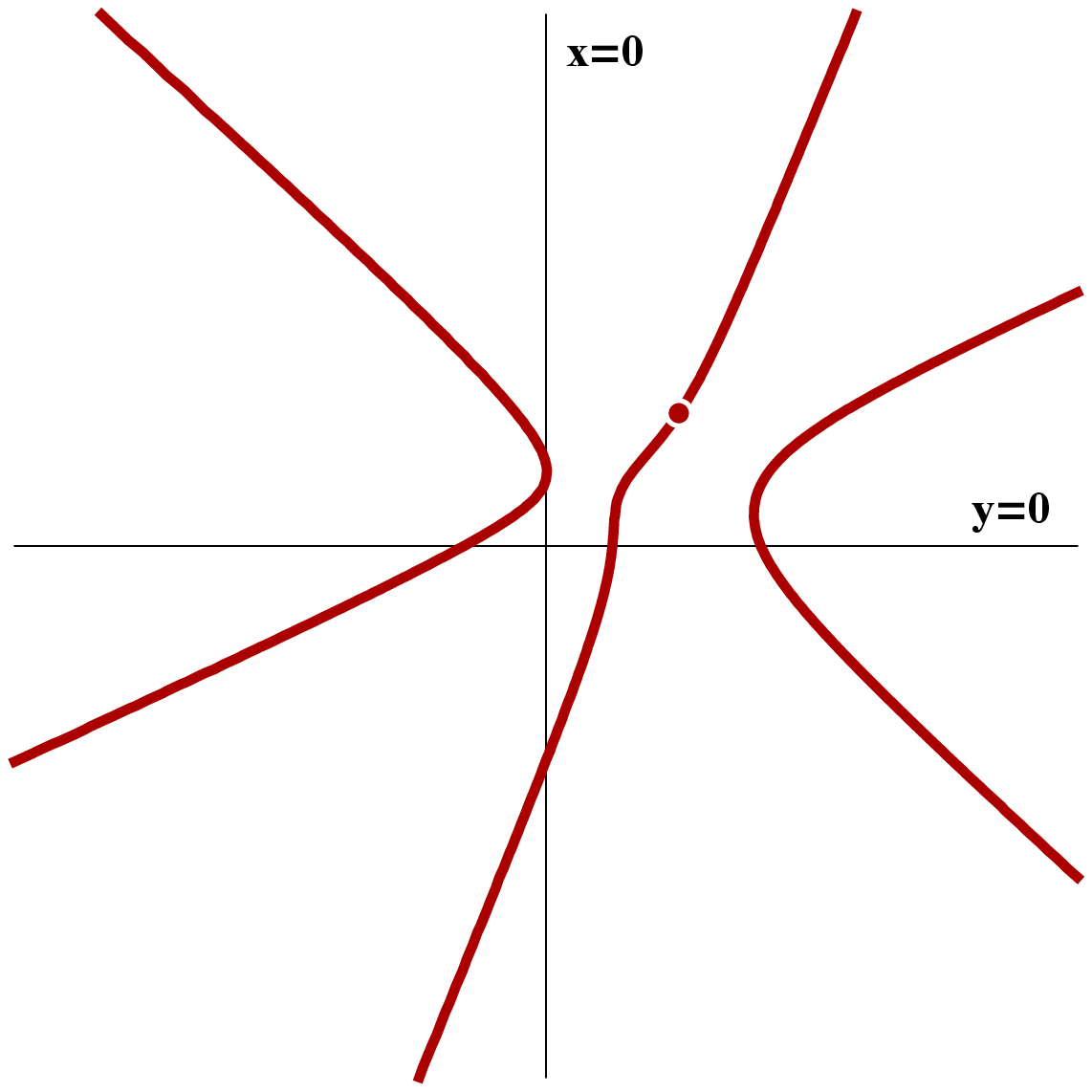}{C0gen_proj.eps}{{\bf Before:} 
  Affine view (left) and projective view (right) of the curve $C_{(0)}$.}
  \bigskip 
  \doppelbildbox{7 true cm}{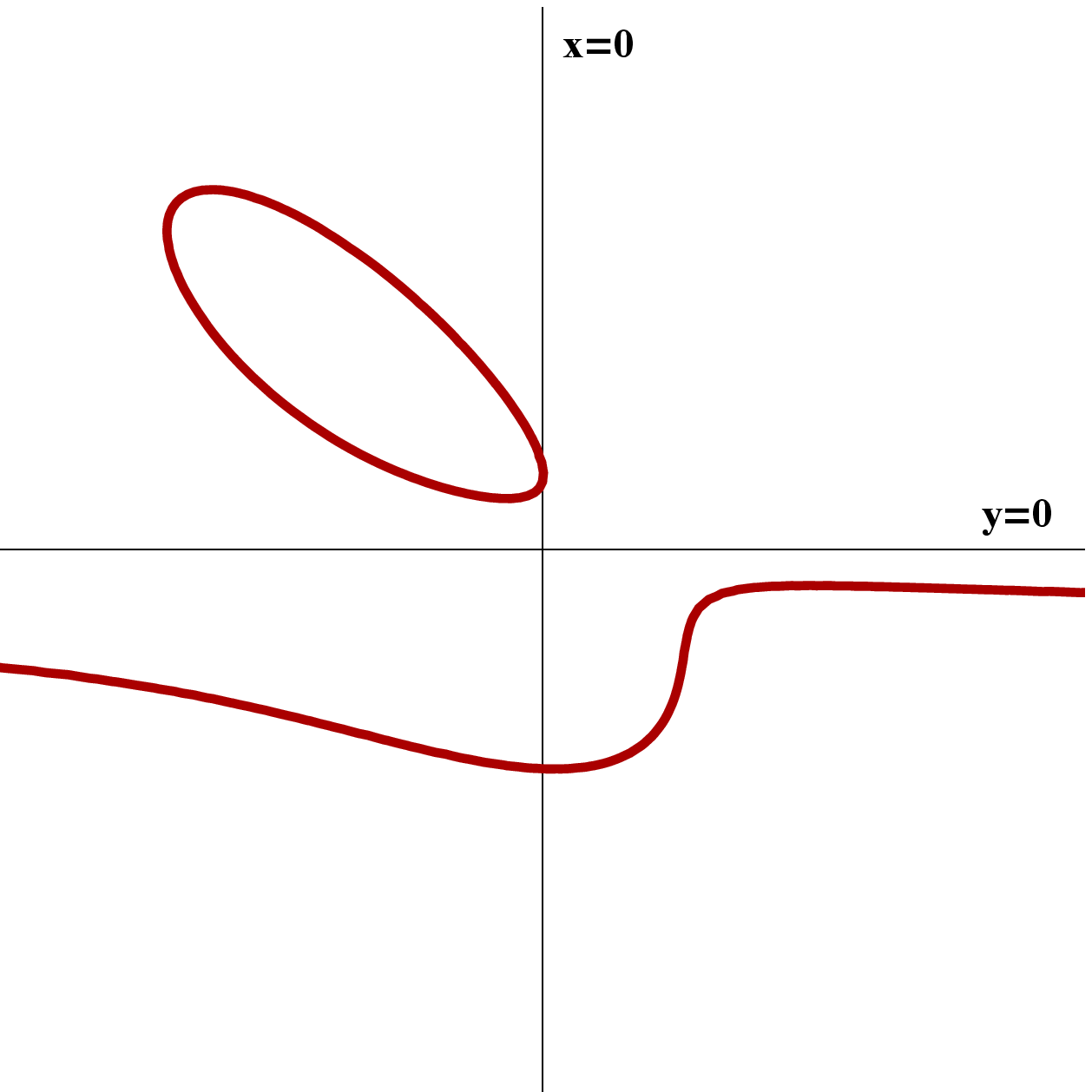}{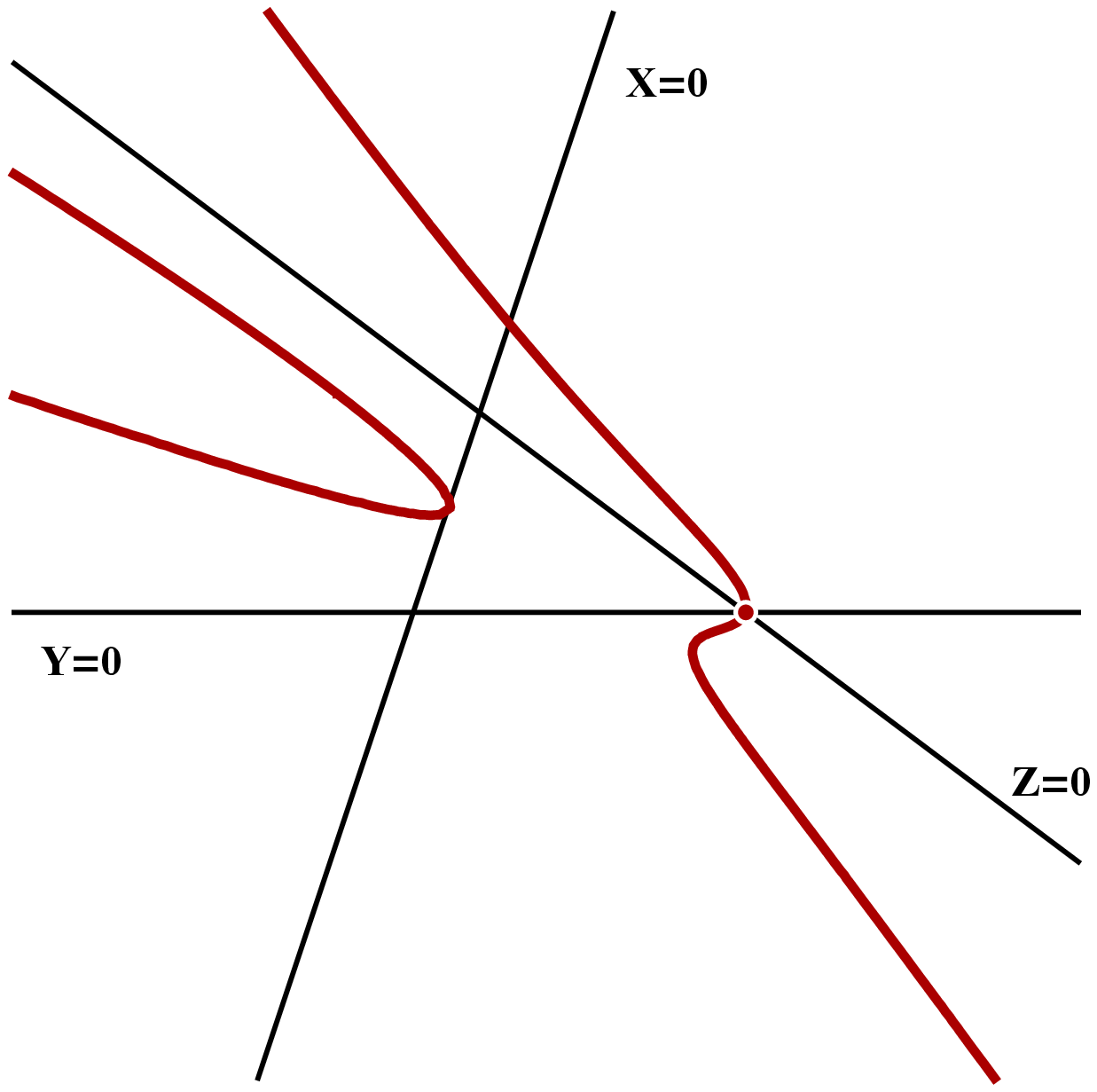}{}{{\bf After:} 
  Affine view (left) and projective view (right) of the curve $C_{(1)}$.}
  \vfill\eject 

{\bf Execution of step 1:} The translation which transforms the distinguished point $p^{(0)} = (p_x,p_y,p_z) = 
(1,p_y,p_z)$ to the special distinguished point $p^{(1)} = (1,0,0)$ is given by the 
linear coordinate transformation 
$$\left[\matrix{X_1\cr Y_1\cr Z_1\cr}\right]\ =\ \left[\matrix{\phantom{-}1&0&0\cr 
  -p_y&1&0\cr -p_z&0&1\cr}\right]\left[\matrix{X_0\cr Y_0\cr Z_0\cr}\right]\quad 
  \hbox{with inverse}\quad \left[\matrix{X_0\cr Y_0\cr Z_0\cr}\right]\ =\ \left[
  \matrix{1&0&0\cr p_y&1&0\cr p_z&0&1\cr}\right]\left[\matrix{X_1\cr Y_1\cr Z_1\cr}\right]\, .
\leqno{(30)}$$ 
In arithmetical terms, the purpose of this transformation is to make the term $X^3$ disappear 
from the original equation, i.e., to render the coefficient $\Gamma^{(1)}_{300}$ zero. In 
fact, plugging $(30)$ into the original equation $(29)$ (with $r=0$) results in the transformed 
cubic $C_{(1)}$ with the equation $\sum_{i+j+k=3}\Gamma^{(1)}_{ijk}X_1^iY_1^jZ_1^k=0$ where 
$$\eqalign{
  \Gamma_{210}^{(1)}\ &=\ \Gamma_{210}^{(0)}p_x^2 + 2\Gamma_{120}^{(0)}p_xp_y + \Gamma_{111}^{(0)}p_z
       + 3\Gamma_{030}^{(0)}p_y^2 + 2\Gamma_{021}^{(0)}p_yp_z + \Gamma_{012}^{(0)}p_z^2,\cr 
  \Gamma_{201}^{(1)}\ &=\ \Gamma_{201}^{(0)}p_x^2 + \Gamma_{111}^{(0)}p_xp_y + 2\Gamma_{102}^{(0)}p_xp_z 
       + \Gamma_{021}^{(0)}p_y^2 + 2\Gamma_{012}^{(0)}p_yp_z + 3\Gamma_{003}^{(0)}p_z^2,\cr 
  \Gamma_{120}^{(1)}\ &=\ \Gamma_{120}^{(0)}p_x^2 + 3\Gamma_{030}^{(0)}p_xp_y + \Gamma_{021}^{(0)}p_xp_z,\cr 
  \Gamma_{111}^{(1)}\ &=\ \Gamma_{111}^{(0)}p_x^2 + 2\Gamma_{021}^{(0)}p_xp_y + 2\Gamma_{012}^{(0)}p_xp_z,\cr 
  \Gamma_{102}^{(1)}\ &=\ \Gamma_{102}^{(0)}p_x^2 + \Gamma_{012}^{(0)}p_xp_y + 3\Gamma_{003}^{(0)}p_xp_z,\cr 
  \Gamma_{030}^{(1)}\ &=\ \Gamma_{030}^{(0)}p_x^2,\cr 
  \Gamma_{021}^{(1)}\ &=\ \Gamma_{021}p_x^2,\cr 
  \Gamma_{012}^{(1)}\ &=\ \Gamma_{012}^{(0)}p_x^2,\cr 
  \Gamma_{003}^{(1)}\ &=\ \Gamma_{003}^{(0)}p_x^2.\cr}\leqno{(31)}$$
(Note that we did not replace $p_x$ by $1$ in order to exhibit the homogeneity of the equations.) In 
the above example, we have 
$$\eqalign{
  &\Gamma_{300}^{(1)}=0,\quad \Gamma_{210}^{(1)}=-2,\quad \Gamma_{201}^{(1)}=-2,\quad \Gamma_{120}^{(1)}=-3,
     \quad \Gamma_{111}^{(1)}=-8,\cr 
  &\Gamma_{102}^{(1)}=4,\quad \Gamma_{030}^{(1)}=-2,\quad \Gamma_{021}^{(1)}=-2,\quad \Gamma_{012}^{(1)}=12, 
     \quad \Gamma_{003}^{(1)}=-8.\cr}$$ 
\vfill\eject  

{\bf Step 2:} In geometrical terms, the purpose of the second transformation is to 
  transform $C_{(1)}$ into a cubic $C_{(2)}$ such that the tangent at the distinguished 
  rational point $(1,0,0)$ is given by the equation $Z=0$. 
  \doppelbildbox{7 true cm}{C1gen_aff.eps}{C1gen_proj.eps}{{\bf Before:} 
   Affine view (left) and projective view (right) of the curve $C_{(1)}$.}
  \bigskip 
  \doppelbildbox{7 true cm}{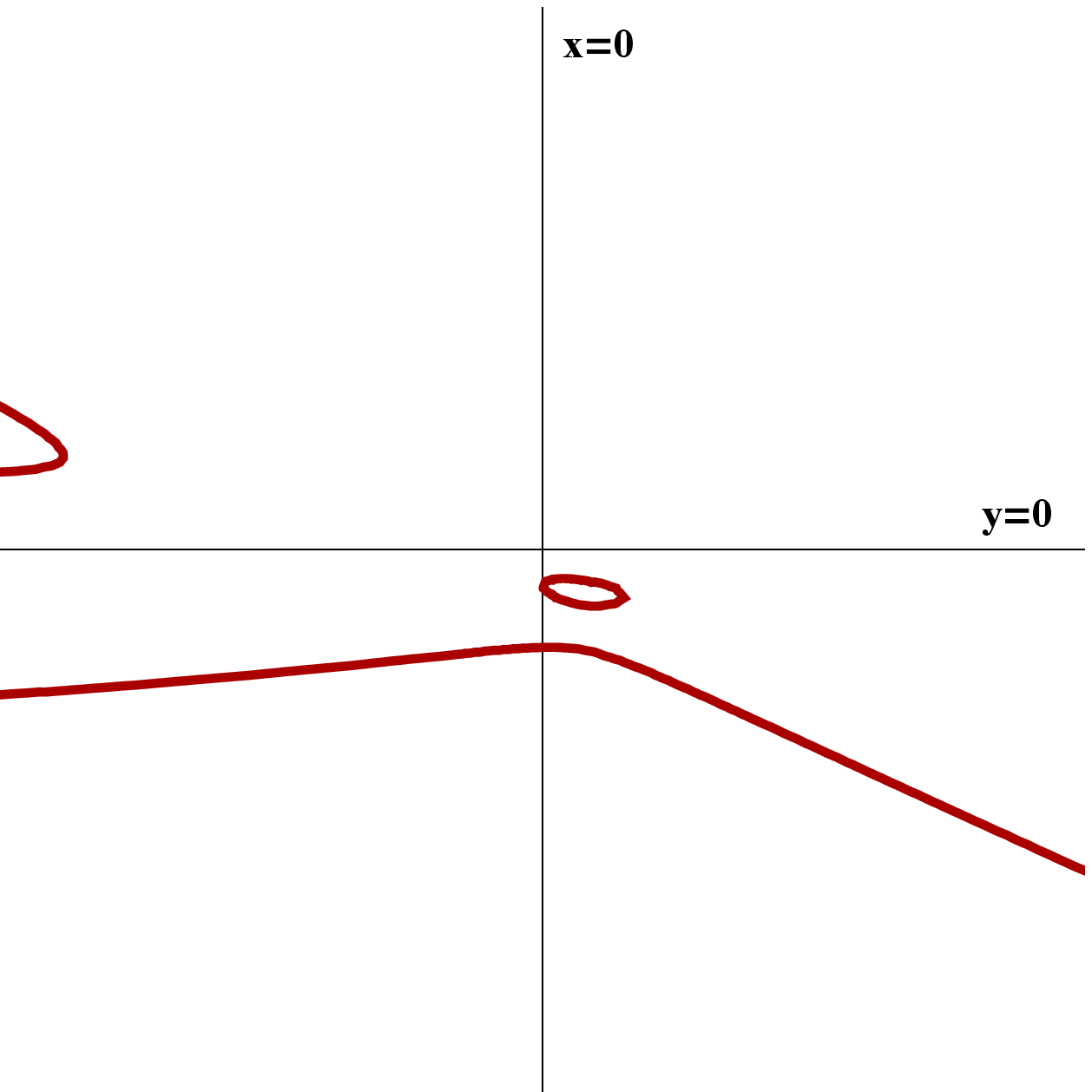}{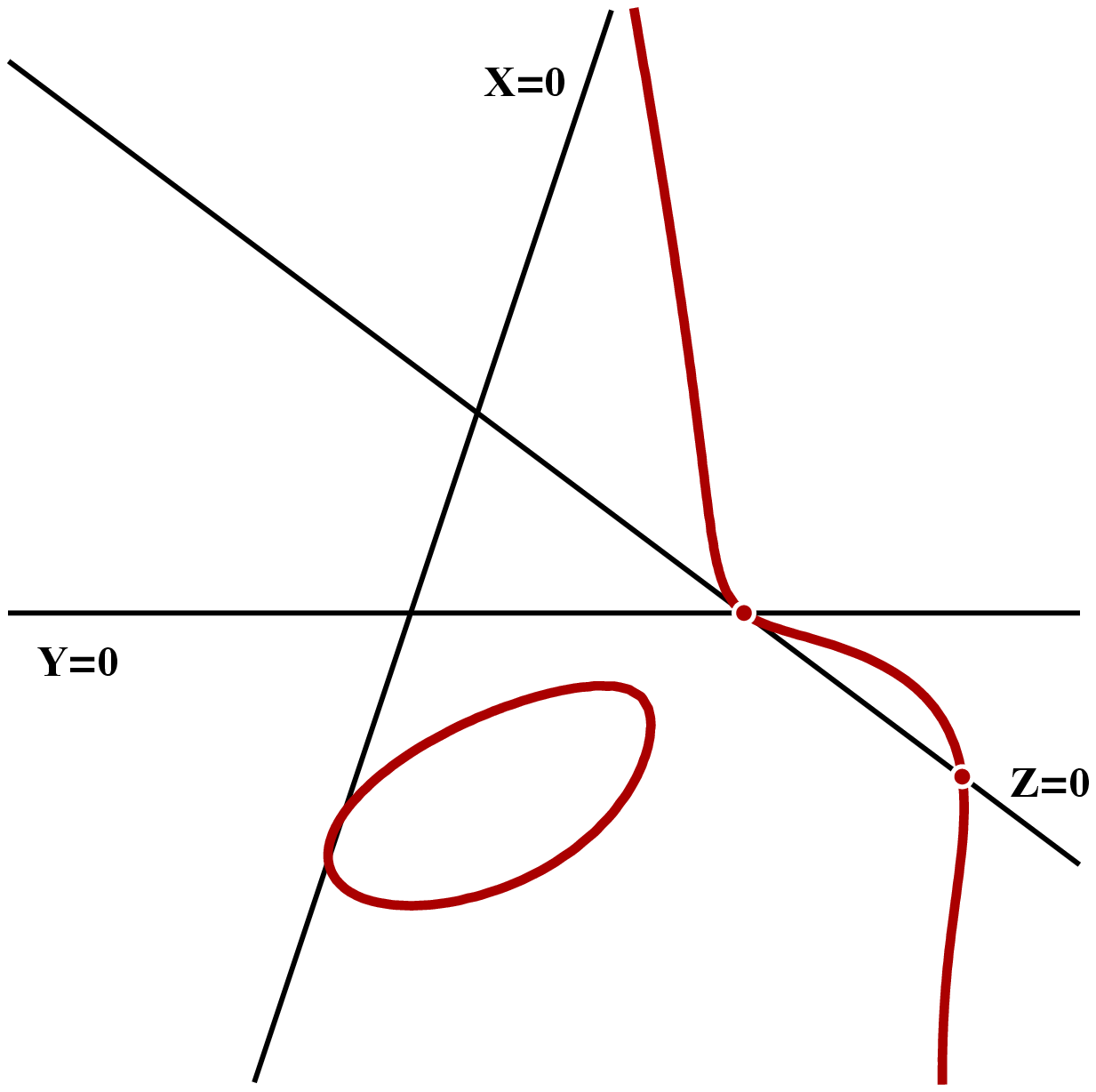}{}{{\bf After:} 
    Affine view (left) and projective view (right) of the curve $C_{(2)}$.}
  \vfill\eject 

{\bf Execution of step 2:} In arithmetical terms, the purpose of this second transformation is to make not only the coefficient 
of $X^3$, but also the coefficient of $X^2Y$ vanish. Note that the tangent of the cubic $C_{(1)}$ 
with coefficients $(31)$ at the point $p^{(1)}=(1,0,0)$ is given by the equation $0=\Gamma_{210}^{(1)}
Y_1 + \Gamma_{201}^{(1)} Z_1 =: g_yY_1 + g_zZ_1$. We may assume that $\Gamma_{201}^{(1)}\neq 0$ 
(since otherwise we could simply exchange the coordinates $Y$ and $Z$ to yield this condition) 
and hence even that $g_z=1$. The goal of the transformation is then accomplished by the linear 
transformation 
$$\left[\matrix{X_2\cr Y_2\cr Z_2\cr}\right] = \left[\matrix{1&0&0\cr 0&1&0\cr 0&g_y&1\cr}\right]\left[\matrix{
  X_1\cr Y_1\cr Z_1\cr}\right]\quad \hbox{with inverse}\quad\left[\matrix{X_1\cr Y_1\cr Z_1\cr}\right]\ =
  \ \left[\matrix{1&\phantom{-}0&0\cr 0&\phantom{-}1&0\cr 0&-g_y&1\cr}\right]\left[\matrix{X_2\cr Y_2\cr Z_2\cr}
  \right].\leqno{(32)}$$ 
Plugging (32) into the equation for $C_{(1)}$ results in the transformed cubic $C_{(2)}$ with the equation 
$\sum_{i+j+k=3}\Gamma^{(2)}_{ijk}X_2^iY_2^jZ_2^k=0$ where 
$$\eqalign{
  \Gamma_{201}^{(2)}\ &=\ \Gamma_{201}^{(1)}g_z^2,\cr 
  \Gamma_{120}^{(2)}\ &=\ \Gamma_{120}^{(1)}g_z^3 - \Gamma_{111}^{(1)}g_yg_z^2 
        + \Gamma_{102}^{(1)}g_zg_y^2,\cr 
  \Gamma_{111}^{(2)}\ &=\ \Gamma_{111}^{(1)}g_z^2 - 2\Gamma_{102}^{(1)}g_yg_z,\cr 
  \Gamma_{102}^{(2)}\ &=\ \Gamma_{102}^{(1)}g_z,\cr 
  \Gamma_{030}^{(2)}\ &=\ \Gamma_{030}^{(1)}g_z^3 - \Gamma_{021}^{(1)}g_z^2g_y 
       + \Gamma_{012}^{(1)}g_zg_y^2 - \Gamma_{003}^{(1)}g_y^3,\cr 
  \Gamma_{021}^{(2)}\ &=\ \Gamma_{021}^{(1)}g_z^2 - 2\Gamma_{012}^{(1)}g_yg_z 
       + 3\Gamma_{003}^{(1)}g_y^2,\cr 
  \Gamma_{012}^{(2)}\ &=\ \Gamma_{012}^{(1)}g_z - 3\Gamma_{003}^{(1)}g_y,\cr 
  \Gamma_{003}^{(2)}\ &=\ \Gamma_{003}^{(1)}.\cr}\leqno{(33)}$$ 
If $(1,0,0)$ happens to be an inflection point of $C_{(2)}$ (i.e., if $\Gamma^{(2)}_{120}=0$), 
we can proceed directly to Step 5 (i.e., we let $C_{(5)}:=C_{(2)}$ after exchanging the 
variables $X_2$ and $Y_2$, because after this change of variables we obtain an equation of 
the form 
$$0 = \Gamma_{300}^{(5)}X_5^3 + \Gamma_{201}^{(5)}X_5^2Z_5 + \Gamma_{111}^{(5)}X_5Y_5Z_5 
  + \Gamma_{102}^{(5)}X_5Z_5^2 + \Gamma_{021}^{(5)}Y_5^2Z_5 + \Gamma_{012}^{(5)}Y_5Z_5^2 
  + \Gamma_{003}^{(5)}Z_5^3$$ 
which is already in Weierstra\ss{} form (with one nonzero coefficient more than in the form 
obtained for $C_{(5)}$ in the other case). Generically, however, the cubic $C_{(2)}$ intersects 
the tangent $Z_2=0$ at the point $(1,0,0)$ in a (simple) second point, namely $p^{(2)} = 
(p_x^{(2)},p_y^{(2)},0) = (\Gamma_{030}^{(2)},-\Gamma_{120}^{(2)},0)$ where $\Gamma^{
(2)}_{120}\not= 0$, and this will be assumed to be the starting point for the next step. 
In our example, we have
$$\eqalign{
  &\Gamma_{300}^{(2)}=0,\quad \Gamma_{210}^{(2)}=0,\quad \Gamma_{201}^{(2)}=-2,\quad \Gamma_{120}^{(2)}=-9,
     \quad \Gamma_{111}^{(2)}=-16,\cr 
  &\Gamma_{102}^{(2)}=-4,\quad \Gamma_{030}^{(2)}=-20,\quad \Gamma_{021}^{(2)}=-50,\quad \Gamma_{012}^{(2)}=-36, 
     \quad \Gamma_{003}^{(2)}=-8;\cr}$$ 
the distinguished rational point is $p^{(2)}=(-20,9,0)$. 
\vfill\eject 

{\bf Step 3:} We transform the cubic $C_{(2)}$ into a general Weierstra\ss{} cubic 
  in the usual form which, in affine coordinates, is $y^2 + a_1xy + a_3y = 
  x^3 + a_2x^2 + a_4x + a_6$. In projective coordinates $(X,Y,Z)$ where $(x,y) = 
  (X/Z,Y/Z)$ this is equivalent to saying that the point $(0,1,0)$ lies on the 
  projective curve and that the line $Z=0$ (the line at infinity) is tangent to $(0,1,0)$.  
 ´Thus the cubic $C_{(3)}$ has the distinguished 
  rational point $p=(1,0,0)$, the tangent at $p$ being given by $Z=0$ and the second intersection 
  point of this tangent with the cubic being given by $q=(0,1,0)$. 
  \doppelbildbox{7 true cm}{C2gen_aff.eps}{C2gen_proj.eps}{{\bf Before:} 
   Affine view (left) and projective view (right) of the curve $C_{(2)}$.}
  \doppelbildbox{7 true cm}{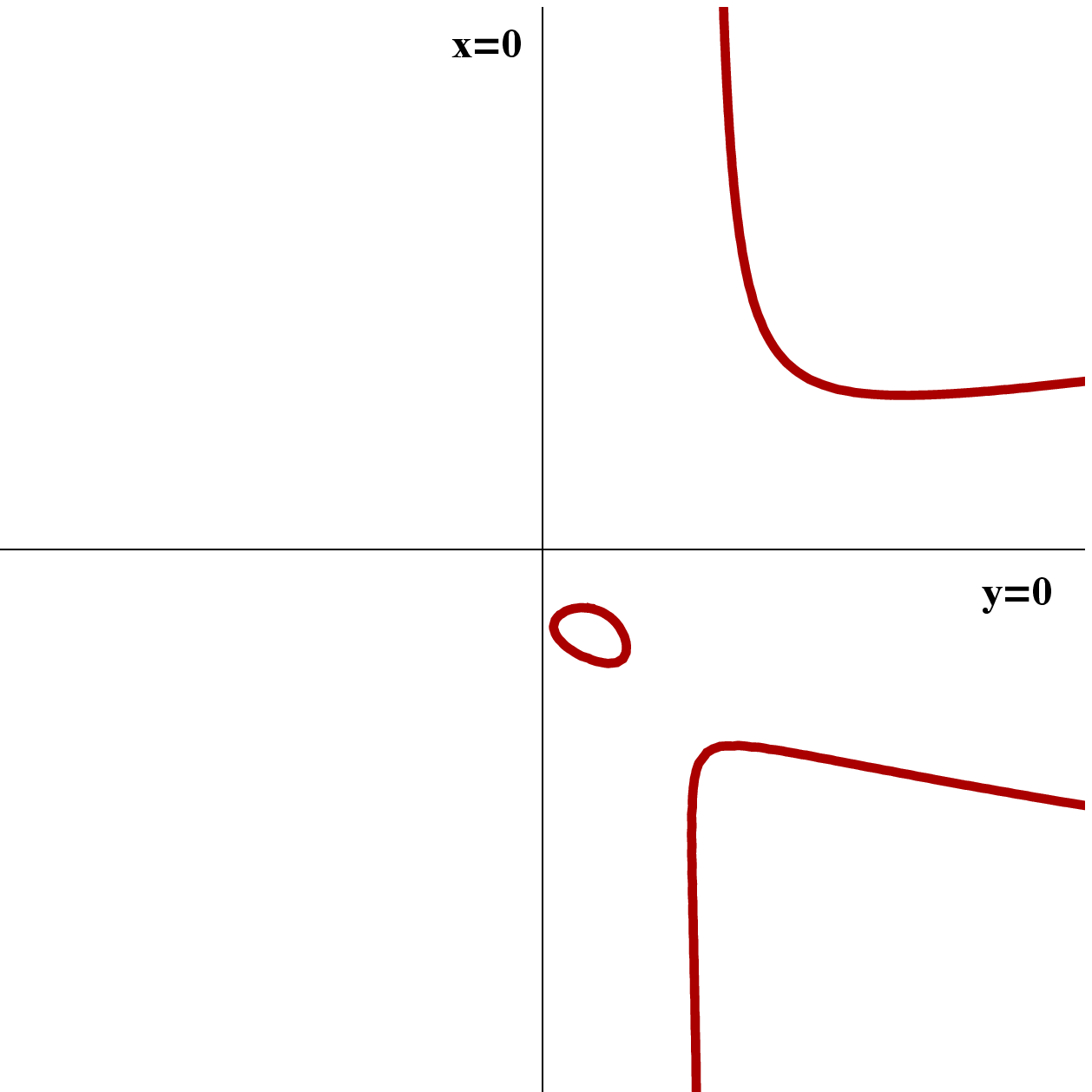}{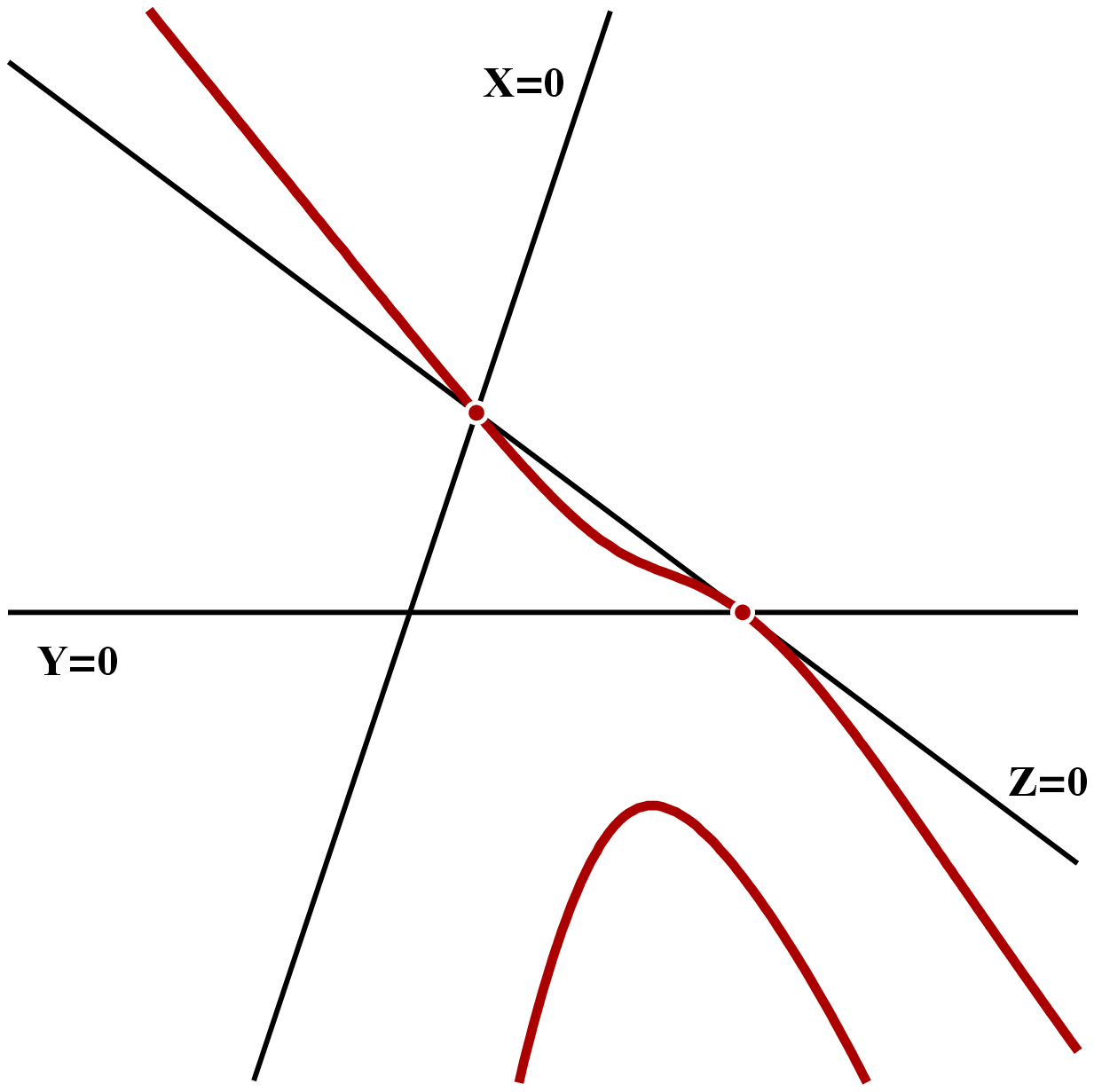}{}{{\bf After:} 
    Affine view (left) and projective view (right) of the curve $C_{(3)}$.}
  \eject 

{\bf Execution of step 3:} If the Weierstra\ss{} cubic $C_{(3)}$ to be constructed is given by 
the equation $0=\sum_{i+j+k=3}\Gamma_{ijk}X^iY^jZ^k$, the above conditions are equivalent to 
$\Gamma_{030} = \Gamma_{210} = \Gamma_{120} = 0$. Geometrically, we want to transform $C_{(2)}$ 
to a cubic $C_{(3)}$ such that the second point of intersection of $C_{(2)}$ with the tangent 
$Z=0$ at $p=(1,0,0)$, i.e., the point $(p_x^{(2)},p_y^{(2)},0)=:(q_x,q_y,0)$ in the previous 
step, is transformed into the point $(0,1,0)$. Since $q_y\not= 0$ by assumption, we may as well 
assume that $q_y=1$. Doing so, we use the coordinate transformation 
$$\left[\matrix{X_3\cr Y_3\cr Z_3\cr}\right]=\left[\matrix{-q_y&\! q_x&\! 0\cr \phantom{-}0&\! 1&\! 0\cr 
  \phantom{-}0&\! 0&\! 1\cr}\right]\left[\matrix{X_2\cr Y_2\cr Z_2\cr}\right]\ \hbox{with inverse}
  \ \left[\matrix{X_2\cr Y_2\cr Z_2\cr}\right] = \left[\matrix{1&\! -q_x&\!\phantom{-}0\cr 
  0&\! -q_y&\!\phantom{-}0\cr 0&\!\phantom{-}0&\! -q_y\cr}\right]\left[\matrix{X_3\cr Y_3\cr 
  Z_3\cr}\right].\leqno{(34)}$$ 
Plugging (34) into the equation for $C_{(2)}$ results in the transformed cubic $C_{(3)}$ with 
the equation $\sum_{i+j+k=3}\Gamma^{(3)}_{ijk}X_2^iY_2^jZ_2^k=0$ where 
$$\eqalign{
  \Gamma_{201}^{(3)}\ &=\ \Gamma_{201}^{(2)}, \cr 
  \Gamma_{120}^{(3)}\ &=\ -q_y\Gamma_{120}^{(2)}, \cr 
  \Gamma_{111}^{(3)}\ &=\ -q_y\Gamma_{111}^{(2)} - 2q_x\Gamma_{201}^{(2)}),\cr 
  \Gamma_{102}^{(3)}\ &=\ -q_y\Gamma_{102}^{(2)}, \cr 
  \Gamma_{021}^{(3)}\ &=\ q_y^2\Gamma_{021}^{(2)} + q_xq_y\Gamma_{111}^{(2)} + q_x^2\Gamma_{201}^{(2)}),\cr 
  \Gamma_{012}^{(3)}\ &=\ q_y^2\Gamma_{012}^{(2)} + q_xq_y\Gamma_{102}^{(2)}),\cr 
  \Gamma_{003}^{(3)}\ &=\ q_y^2\Gamma_{003}^{(2)}.\cr}\leqno{(35)}$$ 
The curve $C_{(3)}$ is characterized by the following properties:
\item{$\bullet$} the point $(1,0,0)$ lies on $C_{(3)}$, i.e., $\Gamma_{300}^{(3)}=0$; 
\item{$\bullet$} the point $(0,1,0)$ lies on $C_{(3)}$, i.e., $\Gamma_{030}^{(3)}=0$; 
\item{$\bullet$} the tangent to $C_{(3)}$ at $(1,0,0)$ is given by $Z_3=0$, i.e., 
   $\Gamma_{210}^{(3)}=0$; 
\item{$\bullet$} the point $(1,0,0)$ is not an inflection point of $C_{(3)}$, i.e., 
  $\Gamma_{120}^{(3)}\neq 0$.
\smallskip\noindent 
In this situation the tangent to $C_{(3)}$ at the point $(0,1,0)$ is given by the equation 
$h_xX_3 + h_zZ_3 = 0$ where $(h_x,h_z) = (\Gamma_{120}^{(3)},\Gamma_{021}^{(3)})$ projectively; 
furthermore, we have $\Gamma_{120}^{(3)}\not= 0$, i.e, $h_x\not= 0$. In our example, 
we have 
$$\eqalign{
  &\Gamma_{300}^{(3)}=0,\quad \Gamma_{210}^{(3)}=0,\quad \Gamma_{201}^{(3)}=-2,\quad \Gamma_{120}^{(3)}=81,
    \quad \Gamma_{111}^{(3)}=64,\cr 
  &\Gamma_{102}^{(3)}=36,\quad \Gamma_{030}^{(3)}=0,\quad \Gamma_{021}^{(3)}=-1970,\quad \Gamma_{012}^{(3)}=-2196, 
      \quad \Gamma_{003}^{(3)}=-648.\cr}$$ 
\vfill\eject 

{\bf Step 4:} We want to transform $C_{(3)}$ to a cubic $C_{(4)}$ such that the tangent to $C_{(4)}$ at the 
  point $(0,1,0)$ is given by the equation $X_4=0$.
  \doppelbildbox{7 true cm}{C3gen_aff.eps}{C3gen_proj.eps}{{\bf Before:} 
   Affine view (left) and projective view (right) of the curve $C_{(3)}$.}
  \bigskip 
  \doppelbildbox{7 true cm}{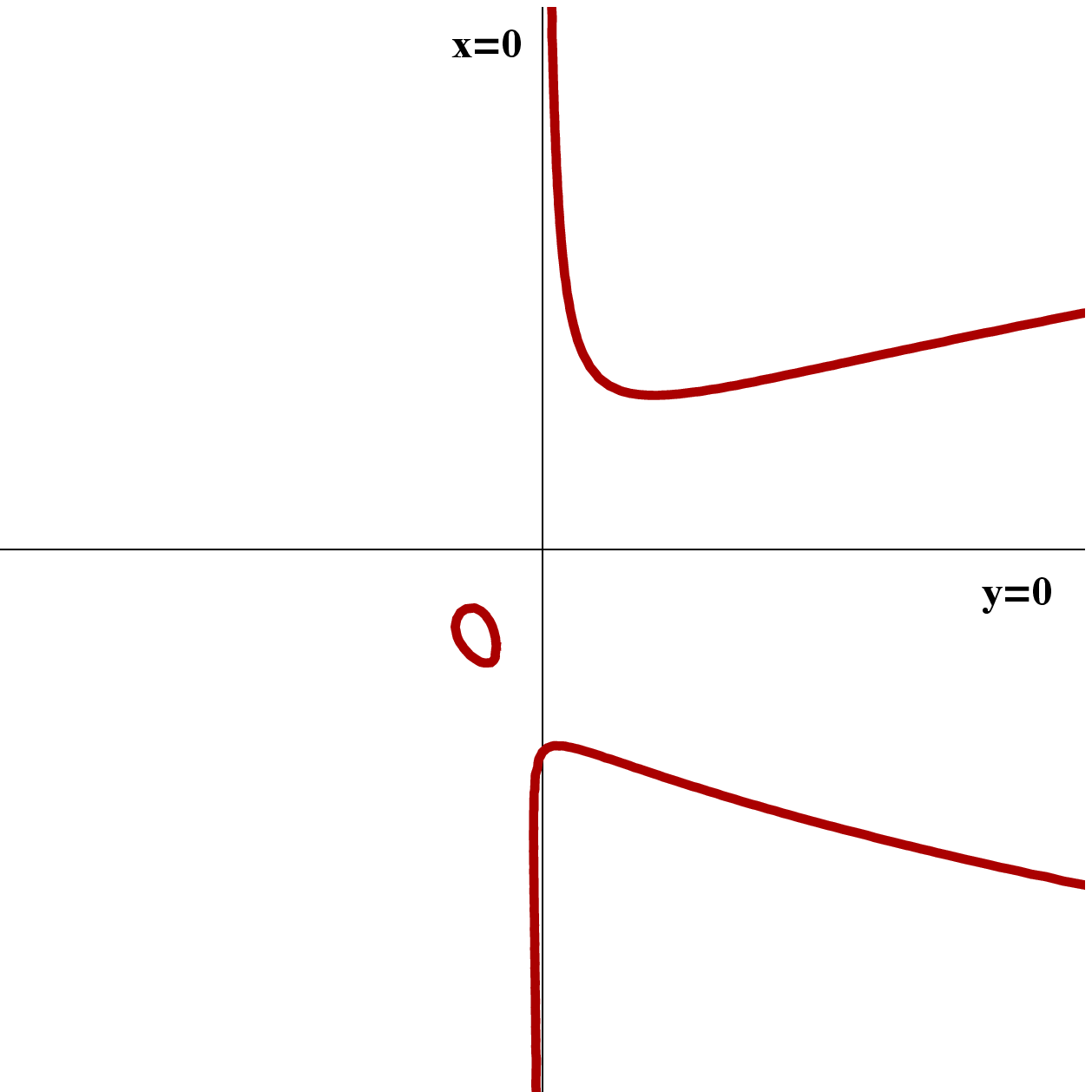}{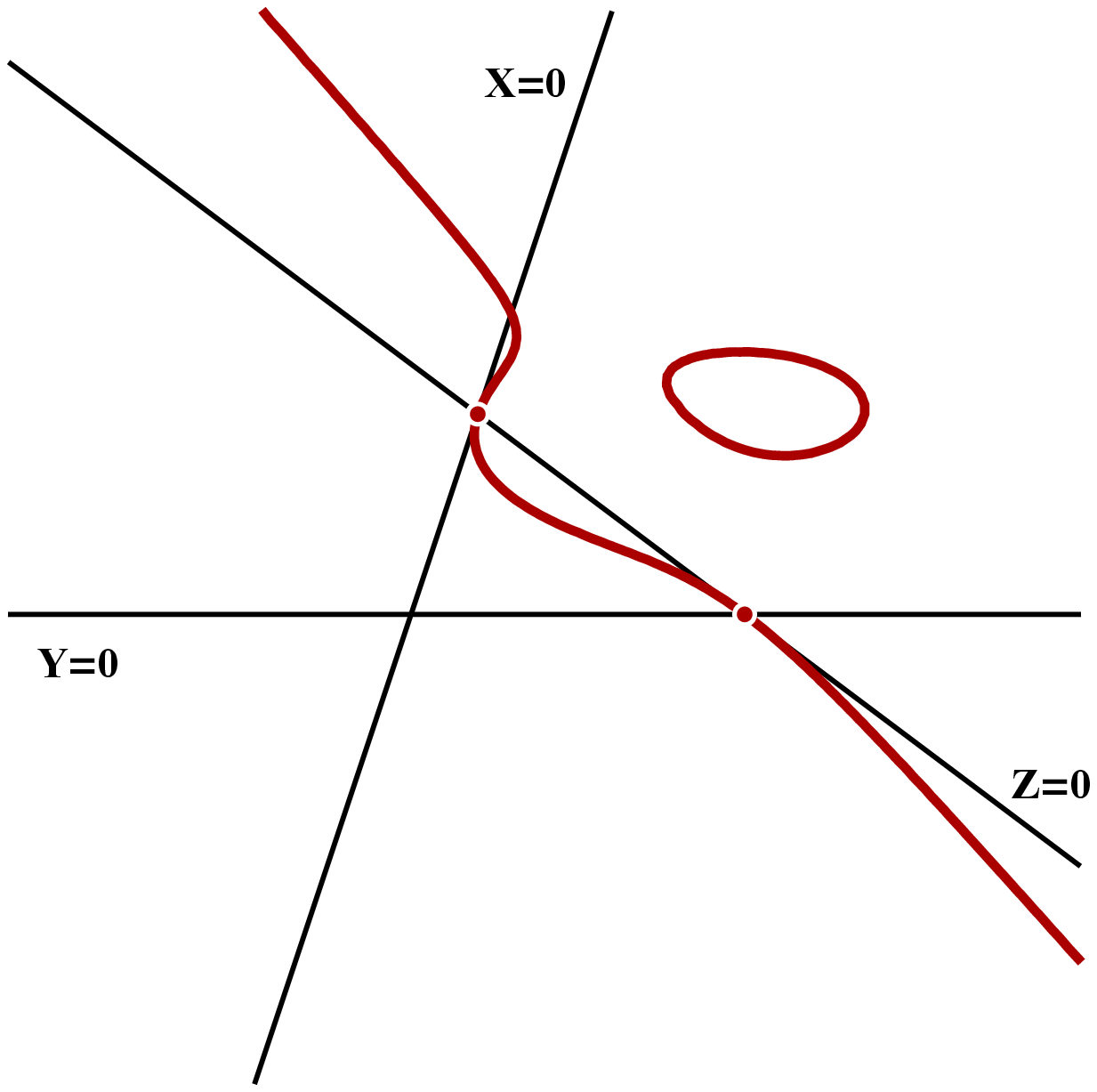}{}{{\bf After:} 
    Affine view (left) and projective view (right) of the curve $C_{(4)}$.}
  \vfill\eject 

{\bf Execution of step 4:} The goal of this transformation is accomplished by the linear transformation 
$$\left[\matrix{X_4\cr Y_4\cr Z_4\cr}\right] = \left[\matrix{h_x&0&h_z\cr 
  0&1&0\cr 0&0&1\cr}\right]\left[\matrix{X_3\cr Y_3\cr Z_3\cr}\right]\quad\hbox{with inverse}\quad 
  \left[\matrix{X_3\cr Y_3\cr Z_3\cr}\right] = \left[\matrix{1&0&-h_z\cr 0&h_x&\phantom{-}0\cr 
  0&0&\phantom{-}h_x\cr}\right]\left[\matrix{X_4\cr Y_4\cr Z_4\cr}\right].\leqno{(36)}$$ 
Plugging (36) into the equation for $C_{(3)}$ results in the transformed cubic $C_{(4)}$ with 
the equation $\sum_{i+j+k=3}\Gamma^{(4)}_{ijk}X_2^iY_2^jZ_2^k=0$ where 
$$\eqalign{
  \Gamma_{201}^{(4)}\ &=\ \Gamma_{201}^{(3)}, \cr 
  \Gamma_{120}^{(4)}\ &=\ h_x\Gamma_{120}^{(3)}, \cr 
  \Gamma_{111}^{(4)}\ &=\ h_x\Gamma_{111}^{(3)}, \cr 
  \Gamma_{102}^{(4)}\ &=\ h_x\Gamma_{102}^{(3)} - 2h_z\Gamma_{201}^{(3)}, \cr 
  \Gamma_{012}^{(4)}\ &=\ h_x^2\Gamma_{012}^{(3)} - h_xh_z\Gamma_{111}^{(3)}, \cr 
  \Gamma_{003}^{(4)}\ &=\ h_x^2\Gamma_{003}^{(3)} - h_xh_z\Gamma_{102}^{(3)} + h_z^2\Gamma_{201}^{(3)}.
  \cr}\leqno{(37)}$$
The curve $C_{(4)}$ is characterized by the following properties:
\item{$\bullet$} the point $(1,0,0)$ lies on $C_{(4)}$, i.e., $\Gamma_{300}^{(4)}=0$; 
\item{$\bullet$} the tangent to $C_{(4)}$ at $(1,0,0)$ is given by $Z_4=0$, i.e., 
   $\Gamma_{210}^{(4)}=0$. 
\item{$\bullet$} the point $(0,1,0)$ lies on $C_{(4)}$, i.e., $\Gamma_{030}^{(4)}=0$; 
\item{$\bullet$} the tangent to $C_{(4)}$ at $(0,1,0)$ is given by $X_4=0$, i.e., 
   $\Gamma_{021}^{(4)}=0$. 
\smallskip\noindent 
In our example, we have 
$$\eqalign{
  &\Gamma_{300}^{(4)}=0,\quad \Gamma_{210}^{(4)}=0,\quad \Gamma_{201}^{(4)}=-2,\quad \Gamma_{120}^{(4)}=6561, 
     \quad \Gamma_{111}^{(4)}=5184,\cr 
  &\Gamma_{102}^{(4)}=-4964,\quad \Gamma_{030}^{(4)}=0,\quad \Gamma_{021}^{(4)}=0,\quad \Gamma_{012}^{(4)}
     =-4\, 195\, 476,\cr 
  &\Gamma_{003}^{(4)}=-6\, 268\, 808.\cr}$$ 
\vfill\eject 

{\bf Step 5:} In this step we will transform $C_{(4)}$ to a Weierstra\ss{} cubic $C_{(5)}$ which is 
  characterized by the conditions that the point $(0,1,0)$ is an inflection point of $C_{(5)}$ such 
  that the tangent to $C_{(5)}$ at $(0,1,0)$ is given by the equation $Z=0$. In arithmetical terms, 
  we want to eliminate the monomial $XY^2$ from the equation of $C_{(4)}$. Once this is done, the 
  terms containing $Y$ are $Y^2Z$, $XYZ$ and $YZ^2$, so that we may split off the factor $Z$ and 
  in the remaining quadratic polynomial then complete the square to get a term $\widetilde{Y}^{2}$. 
  Elimination of $XY^2$ is realized by a quadratic transformation.
  \doppelbildbox{7 true cm}{C4gen_aff.eps}{C4gen_proj.eps}{{\bf Before:} 
   Affine view (left) and projective view (right) of the curve $C_{(4)}$.}
  \doppelbildbox{7 true cm}{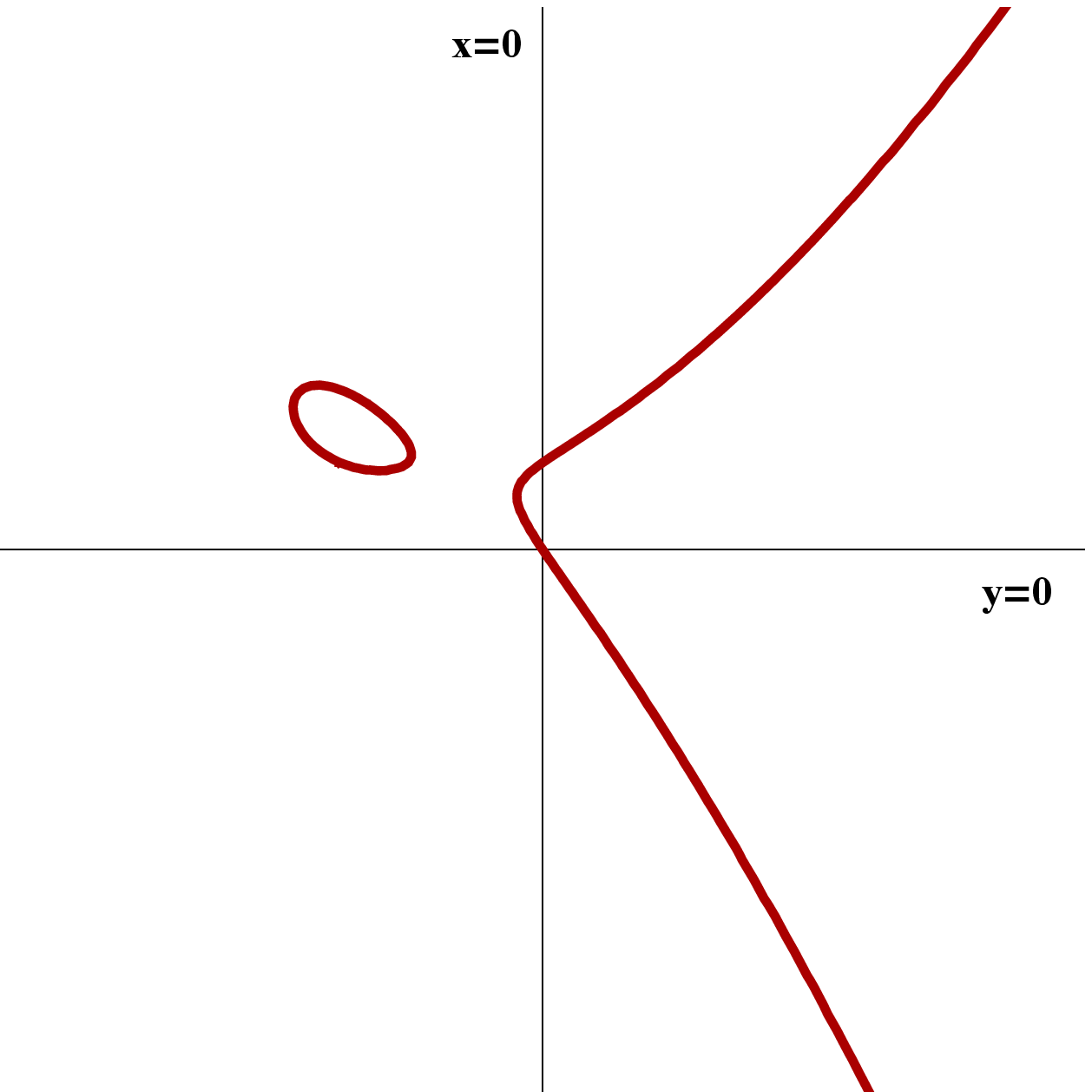}{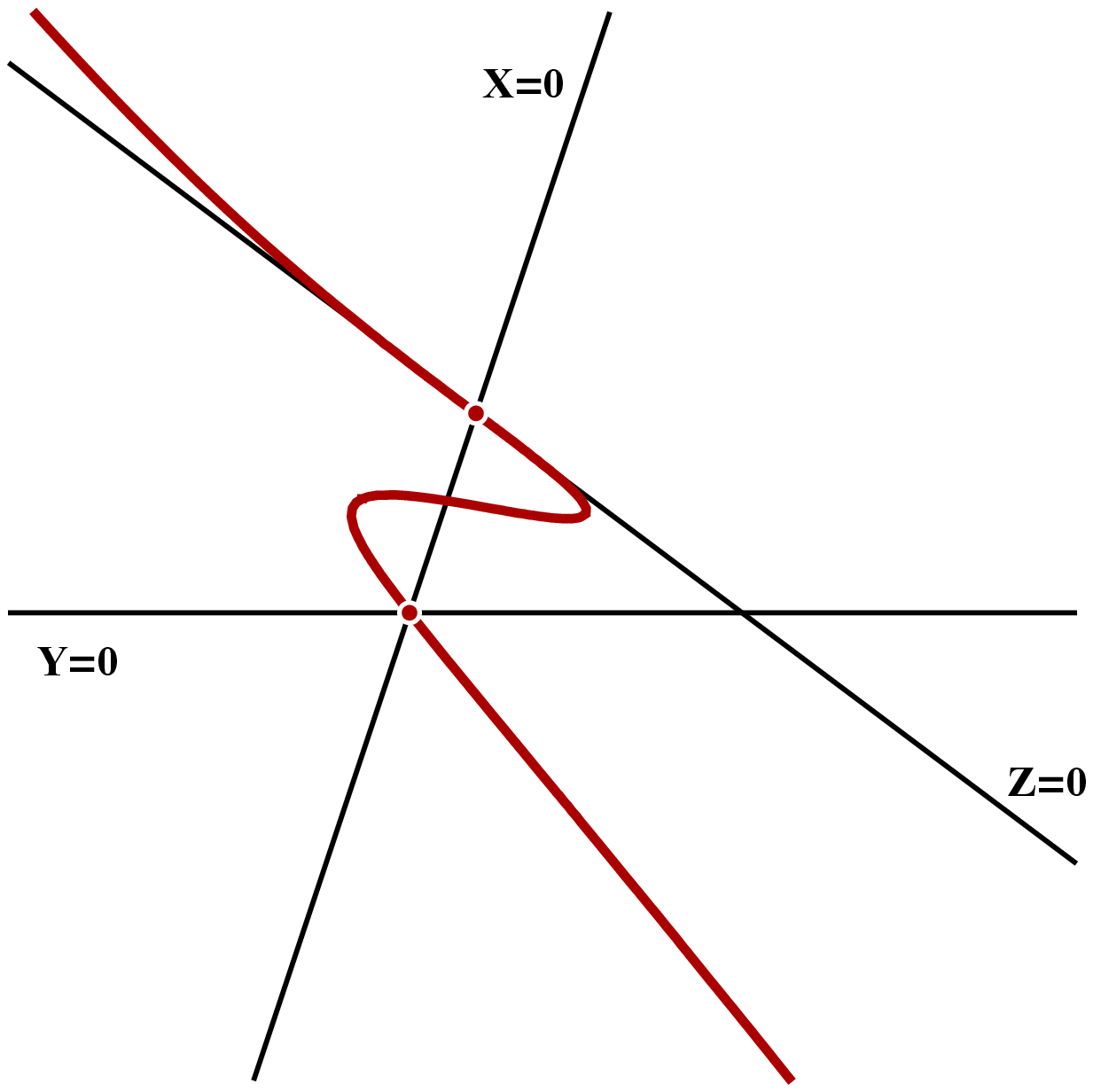}{}{{\bf After:} 
    Affine view (left) and projective view (right) of the curve $C_{(5)}$.}
  \eject 

{\bf Execution of step 5:} The cubic $C_{(4)}$ is given by an equation of the form 
$$0\ =\ \Gamma_{201}^{(4)}X_4^2Z_4 + \Gamma_{120}^{(4)}X_4Y_4^2 + \Gamma_{111}^{(4)}X_4Y_4Z_4 
  + \Gamma_{102}^{(4)}X_4Z_4^2 + \Gamma_{012}^{(4)}Y_4Z_4^2 + \Gamma_{003}^{(4)}Z_4^3. 
  \leqno{(38)}$$
Thus letting $C_{(5)}$ be the image of $C_{(4)}$ under the quadratic transformation 
$$\left[\matrix{X_5\cr Y_5\cr Z_5\cr}\right]\ :=\ \left[\matrix{X_4Z_4\cr X_4Y_4\cr Z_4^2\cr}\right]
  \ =:\ \rho(X_4,Y_4,Z_4)\leqno{(39)}$$ 
then $\rho:C_{(4)}\rightarrow C_{(5)}$ is an isomorphism whose inverse $\psi:C_{(5)}\rightarrow 
C_{(4)}$ is given by  
$$\left[\matrix{X_4\cr Y_4\cr Z_4\cr}\right]\ =\ \left[\matrix{X_5^2\cr Y_5Z_5\cr X_5Z_5\cr}\right]
  \ =:\ \psi (X_5,Y_5,Z_5).\leqno{(40)}$$ 
Note that $\rho$ is not defined at $(1,0,0)$ and $(0,1,0)$, whereas $\psi$ is not defined 
at $(0,1,0)$ and $(0,0,1)$. This is not a problem because of the regularity theorem formulated at 
the end of the introduction; however, we do not need to invoke this general result, but can explicitly 
write down the necessary redefinitions of $\rho$ and $\psi$; see below. Multiplying the equation (38) 
of the cubic $C_{(4)}$ by 
$X_4Z_4^2$ and substituting the transformation $\rho$ gives the cubic $C_{(5)}$ with the equation 
$$0\ =\ \Gamma_{300}^{(5)}X_5^3 + \Gamma_{201}^{(5)}X_5^2Z_5 + \Gamma_{111}^{(5)}X_5Y_5Z_5 
  + \Gamma_{102}^{(5)}X_5Z_5^2 + \Gamma_{021}^{(5)}Y_5^2Z_5 + \Gamma_{012}^{(5)}Y_5Z_5^2
  \leqno{(41)}$$ 
where
$$\Gamma_{300}^{(5)} = \Gamma_{201}^{(4)},  
  \ \Gamma_{201}^{(5)} = \Gamma_{102}^{(4)}, 
  \ \Gamma_{111}^{(5)} = \Gamma_{111}^{(4)},  
  \ \Gamma_{102}^{(5)} = \Gamma_{003}^{(4)},  
  \ \Gamma_{021}^{(5)} = \Gamma_{120}^{(4)},  
  \ \Gamma_{012}^{(5)} = \Gamma_{012}^{(4)}.
  \leqno{(42)}$$ 
Note that writing (41) in the form 
$$0\ =\ \Gamma_{021}^{(5)}Y_5^2Z_5 + \Gamma_{111}^{(5)}X_5Y_5Z_5 + \Gamma_{012}^{(5)}Y_5Z_5^2 
  + \Gamma_{300}^{(5)}X_5^3 + \Gamma_{201}^{(5)}X_5^2Z_5 + \Gamma_{102}^{(5)}X_5Z_5^2 \leqno{(43)}$$ 
shows that this is already a (general) Weierstra\ss{} equation with the additional property that 
$(0,0,1)$ lies on $C_{(5)}$. Furthermore, we have $\Gamma_{021}^{(5)}\neq 0$, since otherwise the 
curve would be singular. In our example, we have 
$$\eqalign{
  &\Gamma_{300}^{(5)}=-2,\quad \Gamma_{210}^{(5)}=0,\quad \Gamma_{201}^{(5)}=-4964,\quad \Gamma_{120}^{(5)}=0, 
     \quad \Gamma_{111}^{(5)}=5184,\cr 
  &\Gamma_{102}^{(5)}=-6\, 268\, 808,\quad \Gamma_{030}^{(5)}=0,\quad \Gamma_{021}^{(5)}=6561,
     \quad \Gamma_{012}^{(5)}=-4\, 195\, 476,\quad \Gamma_{003}^{(5)}=0.\cr}$$ 
\vfill\eject 

$\bullet$ {\bfall Redefinition of $\rho$ around $(1,0,0)$.} On $C_{(4)}$ we have 
$$(\Gamma_{120}^{(4)}Y_4 + \Gamma_{111}^{(4)}Z_4)X_4Y_4\ =
  \ -(\Gamma_{201}^{(4)}X_4^2 + \Gamma_{102}^{(4)}X_4Z_4 + 
  \Gamma_{012}^{(4)}Y_4Z_4 + \Gamma_{003}^{(4)}Z_4^2)Z_4.$$ 
Hence letting $\mu(X_4,Y_4,Z_4):=\Gamma_{120}^{(4)}Y_4 + \Gamma_{111}^{(4)}Z_4$, we have
$$\eqalign{
  &\rho_{5}(X_4,Y_4,Z_4)\ =\mu(X_4,Y_4,Z_4)\left[\matrix{X_4Z_4\cr X_4Y_4\cr Z_4^2\cr}\right] \cr 
  &=\ \left[\matrix{\mu(X_4,Y_4,Z_4)X_4Z_4\cr 
       -(\Gamma_{201}^{(4)}X_4^2 + \Gamma_{102}^{(4)}X_4Z_4 + \Gamma_{012}^{(4)}Y_4Z_4 
           +\Gamma_{003}^{(4)}Z_4^2)Z_4\cr 
       \mu(X_4,Y_4,Z_4)Z_4^2\cr}\right] \cr 
 &=\ Z_4\left[\matrix{\mu(X_4,Y_4,Z_4)X_4\cr 
             -(\Gamma_{201}^{(4)}X_4^2 + \Gamma_{102}^{(4)}X_4Z_4 
                   + \Gamma_{012}^{(4)}Y_4Z_4 + \Gamma_{003}^{(4)}Z_4^2)\cr 
             \mu(X_4,Y_4,Z_4)Z_4\cr}\right] \cr 
 &=\ \left[\matrix{\mu(X_4,Y_4,Z_4)X_4\cr 
          -(\Gamma_{201}^{(4)}X_4^2 + \Gamma_{102}^{(4)}X_4Z_4 + \Gamma_{012}^{(4)}Y_4Z_4 
          + \Gamma_{003}^{(4)}Z_4^2)\cr 
       \mu(X_4,Y_4,Z_4)Z_4\cr}\right] .\cr}$$ 
Since $\mu(1,0,0)=0$, evaluating this representation at $(1,0,0)$ yields  
$\rho(1,0,0)=(0,-\Gamma_{201}^{(4)},0) = (0,1,0)$. \hfill\break\smallskip 

$\bullet$ {\bfall Redefinition of $\rho$ at $(0,1,0)$.} On $C_{4}$ we also have 
$$(\Gamma_{102}^{(4)}X_4 + \Gamma_{012}^{(4)}Y_4 + \Gamma_{003}^{(4)}Z_4)Z_4^2\ =
  \ -(\Gamma_{201}^{(4)}X_4Z_4 + \Gamma_{120}^{(4)}Y_4^2 + \Gamma_{111}^{(4)}Y_4Z_4)X_4.$$
Hence letting $\lambda(X_4,Y_4,Z_4):=\Gamma_{102}^{(4)}X_4 + \Gamma_{012}^{(4)}Y_4 + 
\Gamma_{003}^{(4)}Z_4$, we have 
$$\eqalign{
  &\rho_{5}(X_4,Y_4,Z_4)\ =\ \lambda(X_4,Y_4,Z_4)\left[\matrix{
     X_4Z_4\cr X_4Y_4\cr Z_4^2\cr}\right]\cr 
  &=\ \left[\matrix{\lambda(X_4,Y_4,Z_4)X_4Z_4\cr 
      \lambda(X_4,Y_4,Z_4)X_4Y_4\cr 
      -(\Gamma_{201}^{(4)}X_4Z_4 + \Gamma_{120}^{(4)}Y_4^2 + \Gamma_{111}^{(4)}Y_4Z_4)X_4\cr}
     \right]\cr 
  &=\ X_4\left[\matrix{\lambda(X_4,Y_4,Z_4)Z_4\cr 
      \lambda(X_4,Y_4,Z_4)Y_4\cr 
      -(\Gamma_{201}^{(4)}X_4Z_4 + \Gamma_{120}^{(4)}Y_4^2 + \Gamma_{111}^{(4)}Y_4Z_4\cr}\right]\cr 
  &=\ \left[\matrix{\lambda(X_4,Y_4,Z_4)Z_4\cr 
             \lambda(X_4,Y_4,Z_4)Y_4\cr 
             -(\Gamma_{201}^{(4)}X_4Z_4 + \Gamma_{120}^{(4)}Y_4^2 + \Gamma_{111}^{(4)}Y_4Z_4\cr}
     \right] .\cr}$$
Since $\lambda(0,1,0)=\Gamma_{012}^{(4)}$, this representation yields $\rho(0,1,0) = 
(0,\Gamma_{012}^{(4)},-\Gamma_{120}^{(4)})$.
\vfill\eject 

$\bullet$ {\bfall Redefinition of $\psi$ at $(0,1,0)$.} On $C_{(5)}$ we have
$$(\Gamma_{300}^{(5)}X_5 + \Gamma_{201}^{(5)}Z_5)X_5^2\ =\ -(\Gamma_{111}^{(5)}X_5Y_5 + 
  \Gamma_{102}^{(5)}X_5Z_5 + \Gamma_{021}^{(5)}Y_5^2 + \Gamma_{012}^{(5)}Y_5Z_5)Z_5.$$ 
Hence letting $\sigma(X_5,Y_5,Z_5):=\Gamma_{300}^{(5)}X_5 + \Gamma_{201}^{(5)}Z_5$, we have
$$\eqalign{
  &\psi(X_5,Y_5,Z_5)\ =\ \sigma(X_5,Y_5,Z_5)\left[\matrix{X_5^2\cr 
        Y_5Z_5\cr X_5Z_5\cr}\right] \cr 
  &=\ \left[\matrix{-(\Gamma_{111}^{(5)}X_5Y_5 + \Gamma_{102}^{(5)}X_5Z_5 + 
               \Gamma_{021}^{(5)}Y_5^2 + \Gamma_{012}^{(5)}Y_5Z_5)Z_5\cr 
            \sigma(X_5,Y_5,Z_5)Y_5Z_5\cr 
            \sigma(X_5,Y_5,Z_5)X_5Z_5\cr}\right] \cr 
  &=\ Z_5\left[\matrix{-(\Gamma_{111}^{(5)}X_5Y_5 + \Gamma_{102}^{(5)}X_5Z_5 + 
           \Gamma_{021}^{(5)}Y_5^2 + \Gamma_{012}^{(5)}Y_5Z_5)\cr 
           \sigma(X_5,Y_5,Z_5)Y_5\cr 
           \sigma(X_5,Y_5,Z_5)X_5\cr}\right] \cr 
  &=\ \left[\matrix{-(\Gamma_{111}^{(5)}X_5Y_5 + \Gamma_{102}^{(5)}X_5Z_5 + \Gamma_{021}^{(5)}Y_5^2 
            + \Gamma_{012}^{(5)}Y_5Z_5)\cr \sigma(X_5,Y_5,Z_5)Y_5\cr 
            \sigma(X_5,Y_5,Z_5)X_5\cr}\right] .\cr}$$ 
Since $\sigma(0,1,0)=0$, this representation yields $\psi(0,1,0) = (-\Gamma_{021}^{(5)},0,0) 
= (1,0,0)$.\hfill\break\smallskip 

$\bullet$ {\bfall Redefinition of $\psi$ at $(0,0,1)$.} On $C_{5}$ we also have 
$$(\Gamma_{111}^{(5)}X_5 + \Gamma_{021}^{(5)}Y_5 + \Gamma_{012}^{(5)}Z_5)Y_5Z_5\ =
  \ -(\Gamma_{300}^{(5)}X_5^2 + \Gamma_{201}^{(5)}X_5Z_5 + \Gamma_{102}^{(5)}Z_5^2)X_5.$$ 
Hence letting $\tau(X_5,Y_5,Z_5):=\Gamma_{111}^{(5)}X_5 + \Gamma_{021}^{(5)}Y_5 + 
\Gamma_{012}^{(5)}Z_5$, we have
$$\eqalign{
  &\psi(X_5,Y_5,Z_5)\ =\ \tau(X_5,Y_5,Z_5)\left[\matrix{X_5^2\cr Y_5Z_5\cr X_5Z_5\cr}\right] \cr 
  &=\ \left[\matrix{\tau(X_5,Y_5,Z_5)X_5^2\cr 
        -(\Gamma_{300}^{(5)}X_5^2 + \Gamma_{201}^{(5)}X_5Z_5 + \Gamma_{102}^{(5)}Z_5^2)X_5\cr 
        \tau(X_5,Y_5,Z_5)X_5Z_5\cr}\right] \cr 
  &=\ X_5\left[\matrix{\tau(X_5,Y_5,Z_5)X_5\cr 
      -(\Gamma_{300}^{(5)}X_5^2 + \Gamma_{201}^{(5)}X_5Z_5 + \Gamma_{102}^{(5)}Z_5^2)\cr 
      \tau(X_5,Y_5,Z_5)Z_5\cr}\right] \cr 
  &=\ \left[\matrix{\tau(X_5,Y_5,Z_5)X_5\cr 
      -(\Gamma_{300}^{(5)}X_5^2 + \Gamma_{201}^{(5)}X_5Z_5 + \Gamma_{102}^{(5)}Z_5^2)\cr 
      \tau(X_5,Y_5,Z_5)Z_5\cr}\right] .\cr}$$ 
Since $\tau(0,0,1)=\Gamma_{012}^{(5)}$, this representation yields $\psi(0,0,1) = 
(0,-\Gamma_{102}^{(5)},\Gamma_{012}^{(5)})$.
\vfill\eject 

{\bf Step 6:} In this step we transform the generalized Weierstra\ss{} cubic $C_{(5)}$ into a Weierstra\ss{} 
  cubic $C_{(6)}$, i.e., a cubic for which the only nonzero term containing the variable $Y$ is 
  $Y^2Z$, so that the equation for $C_{(6)}$ has the form $0 = Y^2Z + P(X,Z)$ with a homogeneous 
  polynomial $P(X,Z)$ of degree 3. 
  \doppelbildbox{7 true cm}{C5gen_aff.eps}{C5gen_proj.eps}{{\bf Before:} 
   Affine view (left) and projective view (right) of the curve $C_{(5)}$.}
  \bigskip 
  \doppelbildbox{7 true cm}{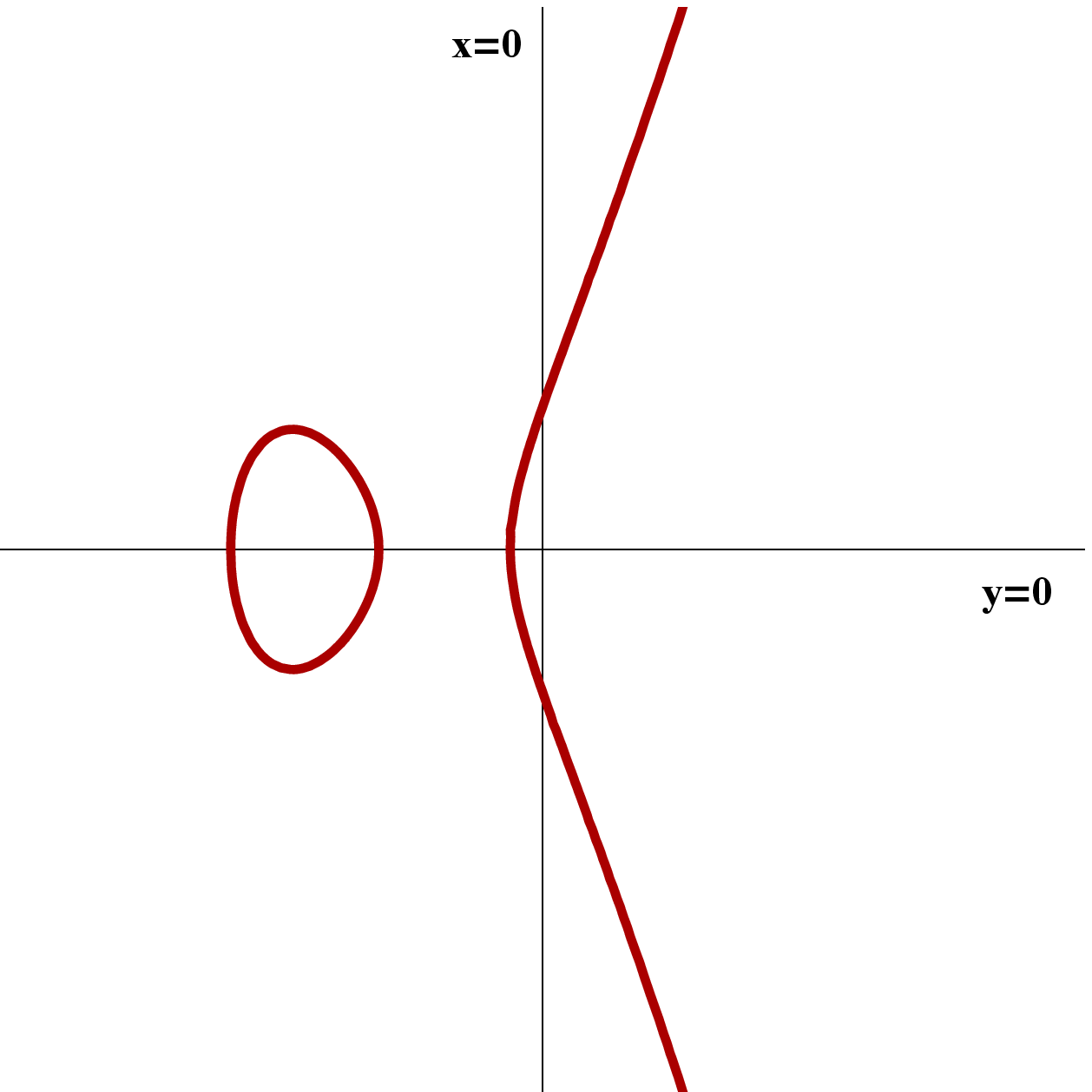}{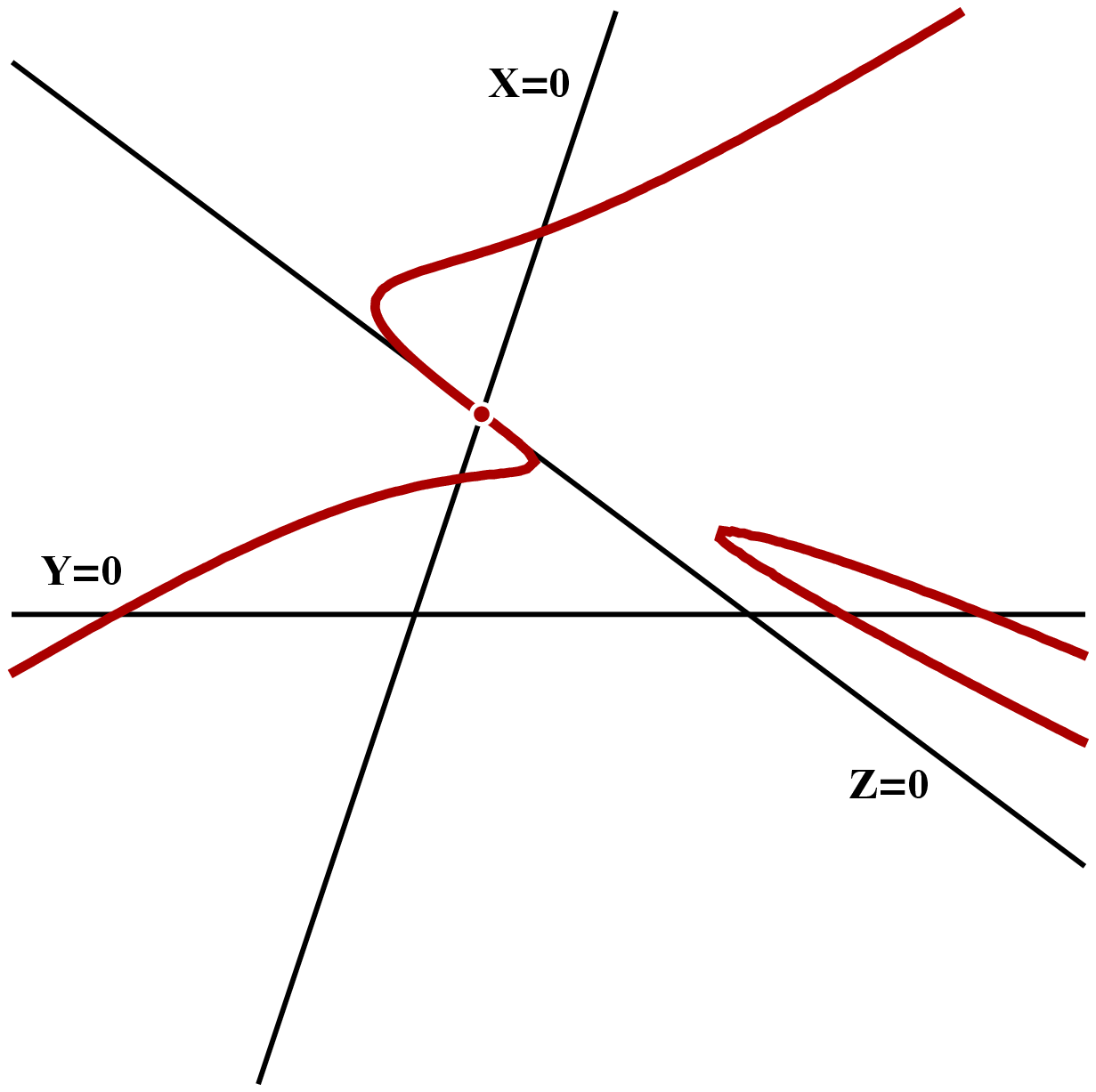}{}{{\bf After:} 
    Affine view (left) and projective view (right) of the curve $C_{(6)}$.}
  \vfill\eject 

{\bf Execution of step 6:} The purpose of this step is accomplished by completing the square in the terms 
containing $Y$. To do so, we multiply the above equation by $4\Gamma_{021}^{(5)}$ and get the 
complete square $(2\Gamma_{021}^{(5)}Y_5 + \Gamma_{111}^{(5)}X_5 + \Gamma_{012}^{(5)}Z_5)^{2}Z_5$ and a 
polynomial in $X_5$ and $Z_5$ of degree 3. More precisely, we use the linear transformation 
$$\left[\matrix{X_6\cr Y_6\cr Z_6\cr}\right] = \left[\matrix{1&0&0\cr 
  \Gamma_{111}^{(5)}&2\Gamma_{021}^{(5)}&\Gamma_{012}^{(5)}\cr 0&0&1\cr}\right]
  \left[\matrix{X_5\cr Y_5\cr Z_5\cr}\right]\leqno{(44)}$$ 
with inverse
$$\left[\matrix{X_5\cr Y_5\cr Z_5\cr}\right] = \left[\matrix{2\Gamma_{021}^{(5)}&0&\phantom{-}0\cr 
  -\Gamma_{111}^{(5)}&1&-\Gamma_{012}^{(5)}\cr \phantom{-}0&\phantom{-}0&
  \phantom{-}2\Gamma_{021}^{(5)}\cr}\right]\left[\matrix{X_6\cr Y_6\cr Z_6\cr}\right].
  \leqno{(45)}$$ 
Then the transformed cubic is given by $0 = \Gamma_{300}^{(6)}X_6^3 + \Gamma_{201}^{(6)}X_6^2Z_6 
+ \Gamma_{102}^{(6)}X_6Z_6^2 + \Gamma_{021}^{(6)}Y_6^2Z_6 + \Gamma_{003}^{(6)}Z_5^3$ with
$$\eqalign{
  \Gamma_{300}^{(6)}\ &=\ 4\Gamma_{021}^{(5)}\Gamma_{300}^{(5)},\cr 
  \Gamma_{201}^{(6)}\ &=\ 4\Gamma_{021}^{(5)}\Gamma_{201}^{(5)}-(\Gamma_{111}^{(5)})^2,\cr 
  \Gamma_{102}^{(6)}\ &=\ 4\Gamma_{021}^{(5)}\Gamma_{102}^{(5)}-2\Gamma_{111}^{(5)}\Gamma_{012}^{(5)},\cr 
  \Gamma_{021}^{(6)}\ &=\ 1,\cr 
  \Gamma_{003}^{(6)}\ &=\ 4\Gamma_{021}^{(5)}\Gamma_{003}^{(5)}-(\Gamma_{012}^{(5)})^2.\cr}
\leqno{(46)}$$
The resulting Weierstra\ss{} cubic $C_{(6)}$ is characterized by the following properties:
\item{$\bullet$} the point $(0,1,0)$ lies on $C_{(6)}$, i.e., $\Gamma_{030}^{(6)}=0$; 
\item{$\bullet$} the tangent to $C_{(6)}$ at $(0,1,0)$ is given by $Z_6=0$, i.e., $\Gamma_{210}^{(6)}=0$; 
\item{$\bullet$} the point $(0,1,0)$ is an inflection point, i.e., $\Gamma_{120}^{(6)}=0$; 
\item{$\bullet$} the only monomial containing $Y$ is $Y^2Z$, i.e., $\Gamma_{111}^{(6)}=0$ and 
   $\Gamma_{012}^{(6)}=0$; 
\item{$\bullet$} the coefficient of $Y^2Z$ is $1$, i.e., $\Gamma_{021}^{(6)}=1$.
\smallskip\noindent 
In our example, we have 
$$\eqalign{
  &\Gamma_{300}^{(6)}=-52\, 488,\quad \Gamma_{210}^{(6)}=0,\quad \Gamma_{201}^{(6)}=-157\, 149\, 072,
     \quad \Gamma_{120}^{(6)}=0,\cr  
  &\Gamma_{111}^{(6)}=0,\quad \Gamma_{102}^{(6)}=-121\, 019\, 901\, 984,\quad \Gamma_{030}^{(6)}=0,\cr 
  &\Gamma_{021}^{(6)}=1,\quad \Gamma_{012}^{(6)}=0,\quad \Gamma_{003}^{(6)}=-17\, 602\, 018\, 866\, 576.\cr}$$ 
\vfill\eject  

{\bf Step 7:} We transform the Weierstra\ss{} cubic $C_{(6)}$ into a Weierstra\ss{} cubic $C_{(7)}$ in 
  normal form with the best possible reduction of the coefficients. The goal is to obtain an equation of 
  the form $Y^2Z + P(X,Z)=0$ such that additionally 
  \item{$\bullet$} the coefficient $\Gamma_{300}$ of $X^3$ is $-1$; 
  \item{$\bullet$} the coefficients $\Gamma_{201}$, $\Gamma_{102}$, $\Gamma_{003}$ do not contain a 
  factor which can be cancelled by a scaling factor of $X$. 
  \doppelbildbox{7 true cm}{C6gen_aff.eps}{C6gen_proj.eps}{{\bf Before:} 
   Affine view (left) and projective view (right) of the curve $C_{(6)}$.}
  \bigskip 
  \doppelbildbox{7 true cm}{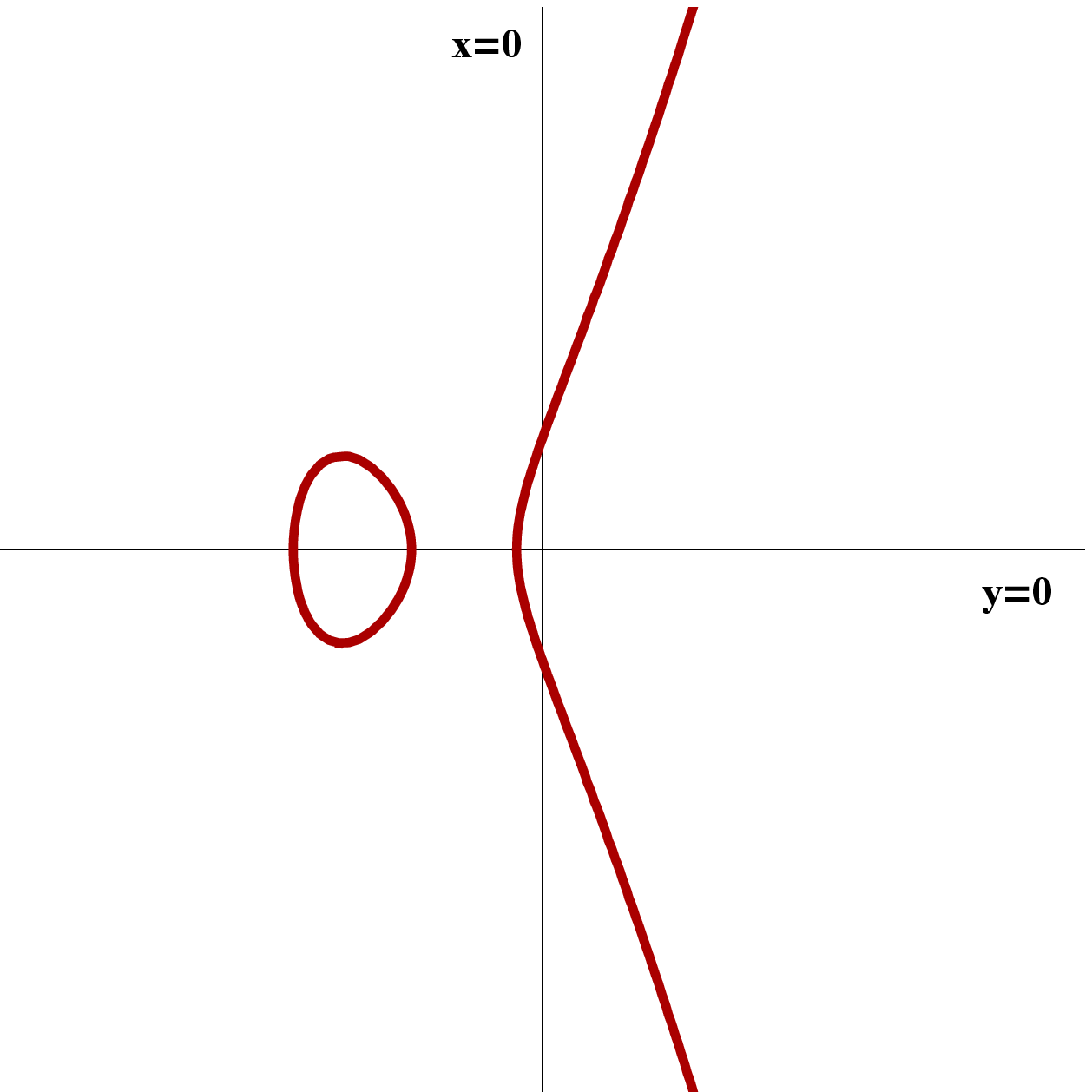}{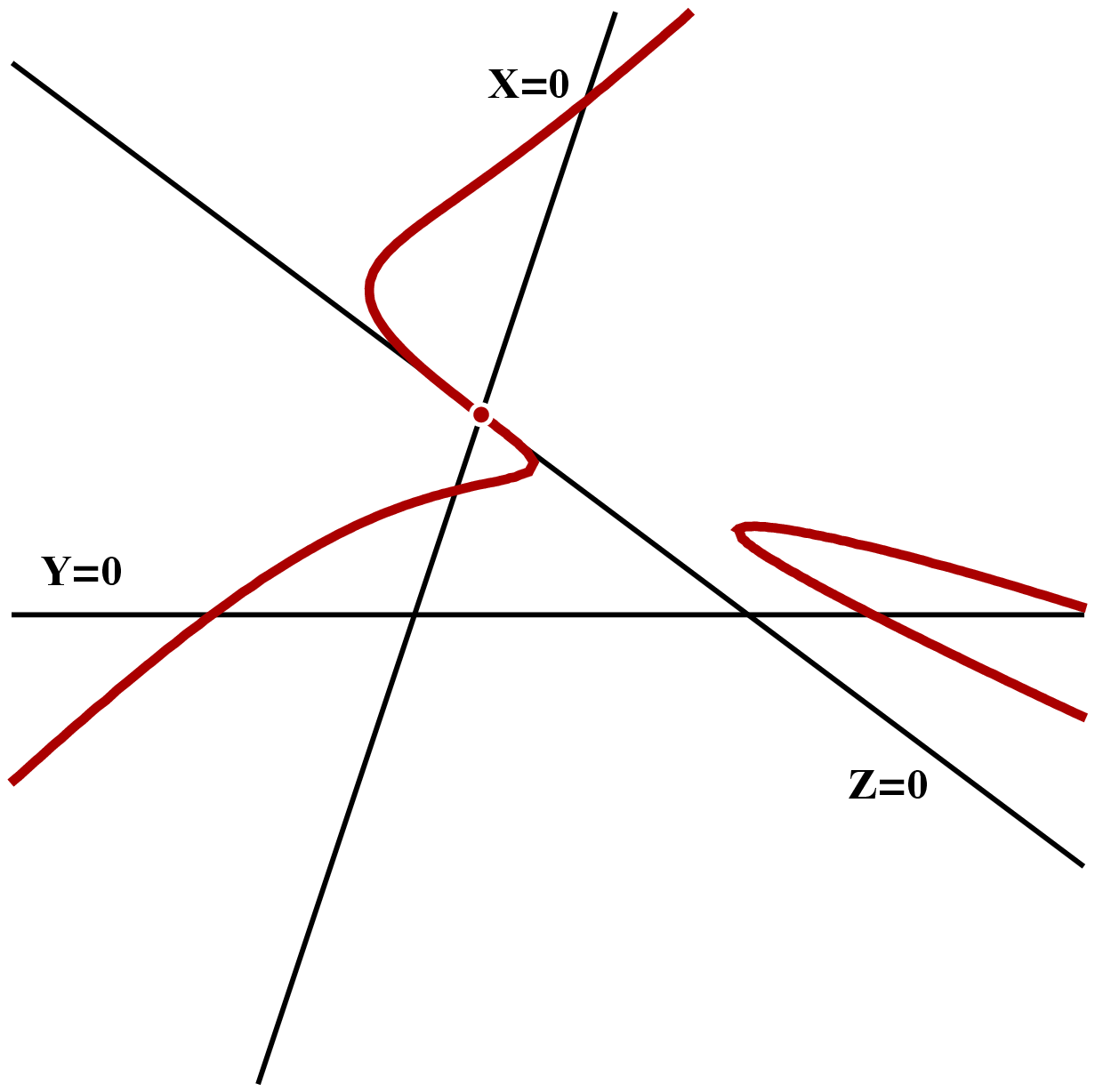}{}{{\bf After:} 
    Affine view (left) and projective view (right) of the curve $C_{(7)}$.}
  \vfill\eject 
´
{\bf Execution of step 7:} Let $\delta :=-\Gamma_{300}^{(6)}$ and let $\varphi$ be a (maximal) factor 
such that 
$$\varphi^2|\Gamma_{201}^{(6)}, \quad 
  \varphi^4|\Gamma_{300}^{(6)}\Gamma_{102}^{(6)},\quad 
  \varphi^6|(\Gamma_{300}^{(6)})^2\Gamma_{003}^{(6)}.\leqno{(47)}$$ 
We then consider the linear transformation
$$\left[\matrix{X_7\cr Y_7\cr Z_7\cr}\right] = \left[\matrix{\delta\phi&0&0\cr 0&\delta&0\cr 
  0&0&\varphi^3\cr}\right]\left[\matrix{X_6\cr Y_6\cr Z_6\cr}\right]\quad\hbox{with inverse}\quad 
  \left[\matrix{X_6\cr Y_6\cr Z_6\cr}\right] = \left[\matrix{\varphi^2&0&0\cr 0&\varphi^3&0\cr 
  0&0&\delta\cr}\right]\left[\matrix{X_7\cr Y_7\cr Z_7\cr}\right] \leqno{(48)}$$ 
Then the transformed cubic $C_{(7)}$ has the only nonzero coefficients 
$$\Gamma_{201}^{(7)} = \Gamma_{201}^{(6)}/\varphi^2, \quad  
  \Gamma_{102}^{(7)} = \delta\Gamma_{102}^{(6)}/\varphi^4, \quad 
  \Gamma_{003}^{(7)} = \delta^2\Gamma_{003}^{(6)}/\varphi^6.\leqno{(49)}$$ 
Dehomogenization with respect to the variable $Z=Z_7$, i.e., letting $x=X/Z$ and $y=Y/Z$, 
yields a classical Weierstra\ss{} equation of the form $y^2 = x^3 + a_2x^2 + a_4x + a_6$ 
where $a_2=-\Gamma_{201}^{(7)}$, $a_4=-\Gamma_{102}^{(7)}$ and $a_6=-\Gamma_{003}^{(7)}$. 
In the general situation no further transformations are possible which lead to essential 
simplifications. Only if the polynomial $x^3 + a_2x^2 + a_4x + a_6$ has three rational (in 
fact integral) roots the equation can be further simplified, namely to the form $y^2 = 
x(x+A)(x+B)$. In our example, we have 
$$\eqalign{
  &\Gamma_{300}^{(7)}=-1,\ \Gamma_{210}^{(7)}=0,\ \Gamma_{201}^{(7)}=-5988,
     \ \Gamma_{120}^{(7)}=0,\ \Gamma_{111}^{(7)}=0,\cr 
  &\Gamma_{102}^{(7)}=-9\, 222\, 672,\ \Gamma_{030}^{(7)}=0,\ \Gamma_{021}^{(7)}=1, 
     \ \Gamma_{012}^{(7)}=0,\ \Gamma_{003}^{(7)}=-2\, 682\, 825\, 616.\cr}$$ 
\hfill\break\smallskip 

This completes the transformation. Let us review what happens to the distinguished 
rational point on the intersection of quadrics we started with. This point is first 
mapped to the point $z\in C_{(0)}$ given by formula $(16)$ above. Next, the point $z$ 
is mapped to $(1,0,0)\in C_{(1)}$, and this point is mapped to itself in the transformation 
from $C_{(1)}$ to $C_{(2)}$. If $(1,0,0)$ happens to be an inflection point of $C_{(2)}$ 
we apply a coordinate exchange which maps $(1,0,0)$ to $(0,1,0)$; in the generic 
case, the point $(1,0,0)$ remains fixed both during the transition from $C_{(2)}$ to 
$C_{(3)}$ and during the transition from $C_{(3)}$ to $C_{(4)}$ and is then mapped 
to $(0,1,0)\in C_{(5)}$ by the subsequent quadratic transformation. Thus in all cases, 
the original distinguished point is mapped to $(0,1,0)\in C_{(5)}$, and this point 
stays fixed under the remaining transformations. 
\vfill\eject 

\leftline{\titlefont 4. Example: Euler's concordant forms}
\bigskip 

As a first example, let us apply the above theory to the intersection $Q_{M,N}$ of the two quadrics 
$X^2+MY^2 =Z^2$ and $X^2+NY^2=W^2$ where $M\not= N$ are nonzero integers; according to Euler, the 
numbers $M$ and $N$ are called concordant if $Q_{M,N}$ possesses a rational point which is nontrivial 
in the sense that $Y\not= 0$. As a consequence of the previous discussion, the quadric intersection 
$Q_{M,N}$ is isomorphic to a plane elliptic curve $E_{M,N}$ given by a Weierstra\ss{} equation. Note 
that in [3] and [7] the authors use a mapping of degree $4$ instead of a biregular morphism; this is 
somewhat surprising, since in this way some information on torsion points is lost and since, as we 
shall see, the biregular morphism is given by a rather simple linear mapping. Also note that in [8] a 
linear isomorphism from $Q_{M,N}$ to a smooth plane cubic is given which is rather similar to the 
first one of our curves, but no transformation of this curve to Weierstra\ss{} form is carried out. 
To apply the theory developed in the previous paragraph we write $(X_0,X_1,X_2,X_3):=(Y,X,Z,W)$ and 
thus consider the quadrics $Q_1$ given by $MX_0^2+X_1^2-X_2^2=0$ and $Q_2$ given by $NX_0^2+X_1^2-X_3^2=0$. 
This corresponds to the equations $(5)$ with $A=\hbox{\rm diag}(M,1,-1,0)$ and $B=\hbox{\rm diag}(N,1,
0,-1)$. The intersection $Q_1\cap Q_2$ contains the four trivial rational points $(0,1,\pm 1,\pm 1)$; 
we choose $x:=(0,1,1,1)$ as the distinguished rational point. We carry out the procedure described in 
the previous paragraph and visualize the various steps by images generated for the special case 
$(M,N)=(3,2)$. The isomorphism used in Section 2 maps $Q_1\cap Q_2$ to the plane cubic whose nonzero 
coefficients are 
$$\eqalign{\Gamma_{210}^{(0)}\ &=\ N-M,\cr 
           \Gamma_{201}^{(0)}\ &=\ -N,\cr 
           \Gamma_{021}^{(0)}\ &=\ -1,\cr 
           \Gamma_{012}^{(0)}\ &=\ 1,\cr}$$ 
and $\varphi$ maps $x=(0,1,1,1)$ to $y=(1,0,0)$. In  the special case $(M,N)=(3,2)$ this means that 
$$\Gamma_{210}^{(0)}=-1,\qquad \Gamma_{201}^{(0)}=-2, \qquad \Gamma_{021}^{(0)}=-1,\qquad 
  \Gamma_{012}^{(0)}=1.$$ 
\vfill\eject 
\doppelbildbox{7 true cm}{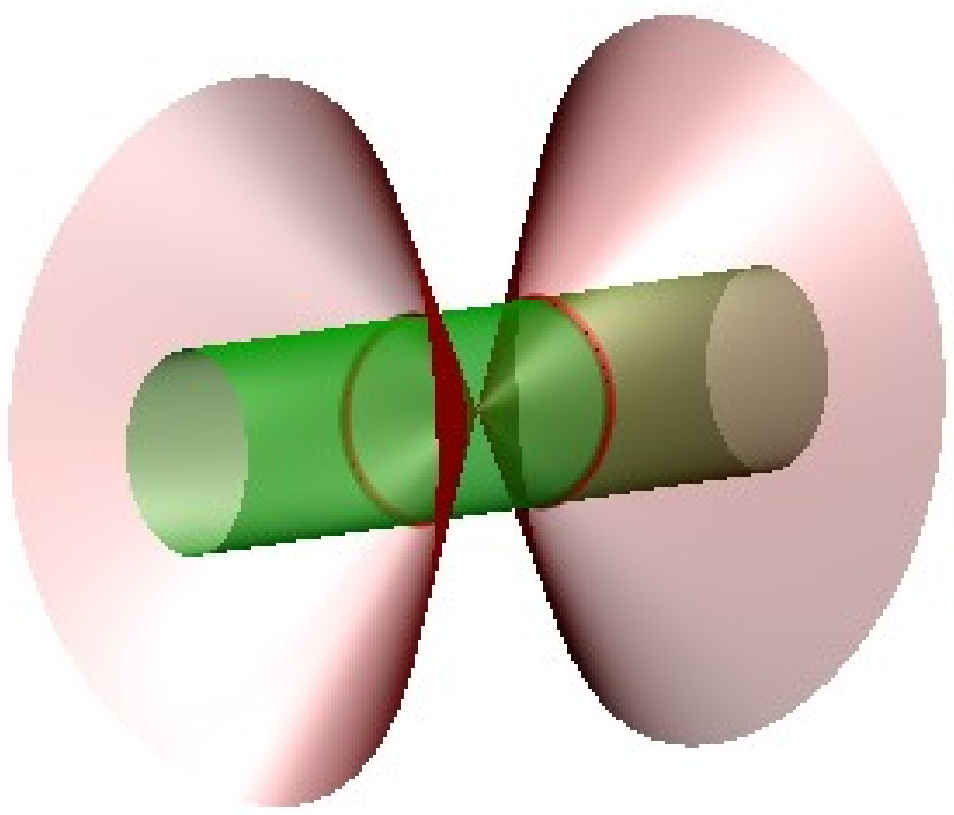}{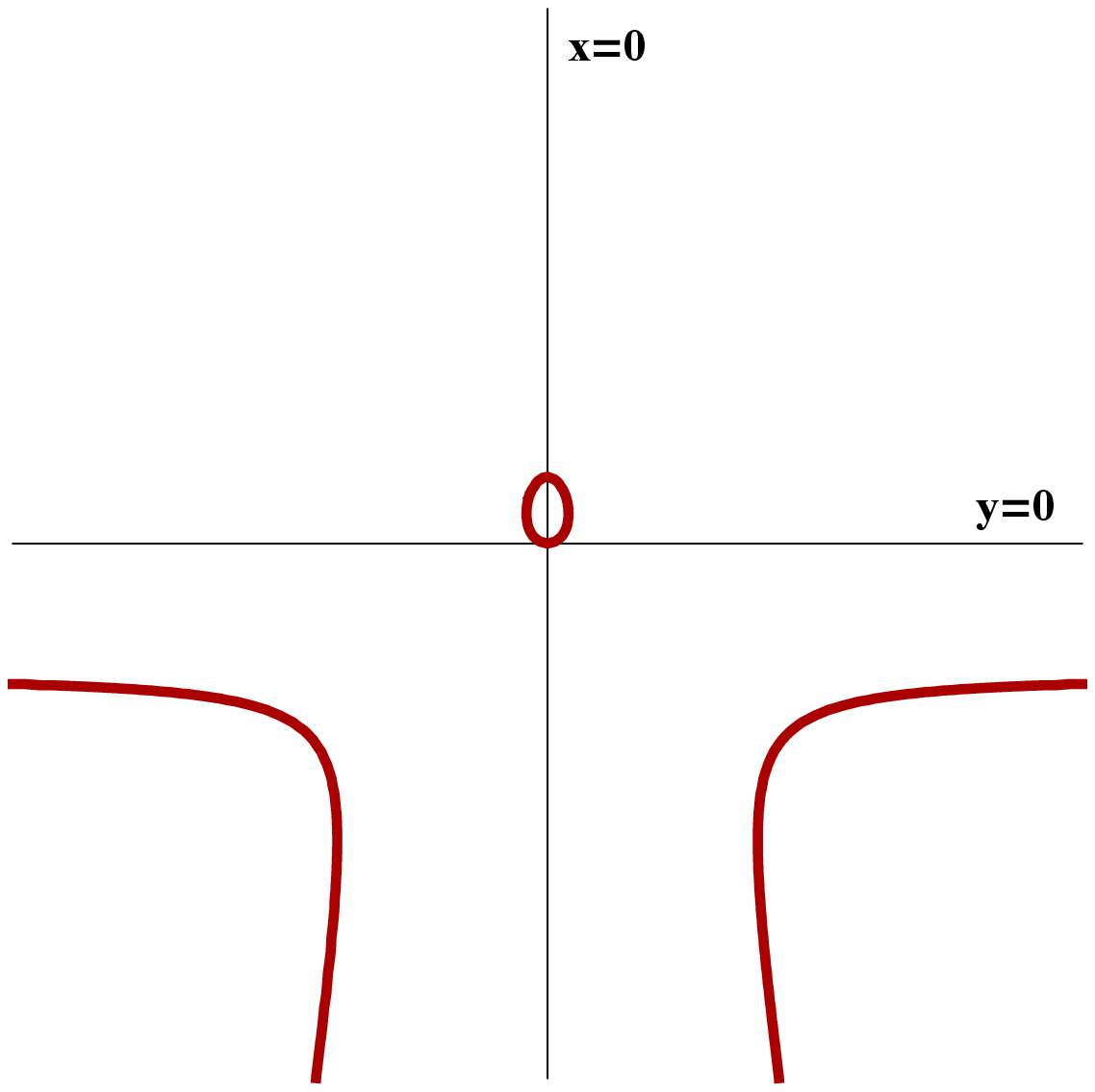}{{\bf Left:} Affine view $X_3=1$ of 
the intersection of the quadrics $MX_0^2+X_1^2=X_2^2$ and $NX_0^2+X_1^2=X_3^2$. {\bf Right:} Affine view 
of the cubic $C_{(0)}$ to which this quadric intersection is initially transformed.} 
\vfill\eject 

{\bf Step 1:} This step is superfluous, because the distinguished point is already $(1,0,0)$.
   Hence $C_{(1)}=C_{(0)}$. \smallskip 
{\bf Step 2:} The tangent of the curve $C_{(0)}=C_{(1)}$ at the point $y=(1,0,0)$
   is given by $0=(N-M)Y_1 - NZ_1$. The curve $C_{(2)}$ has the nonzero coefficients 
   $\Gamma_{201}^{(2)}=N^2$, $\Gamma_{030}^{(2)}=M(M-N)$, $\Gamma_{021}^{(2)}=2M-N$ and 
   $\Gamma_{012}^{(2)}=1$. In the special case $(M,N)=(3,2)$ this means that $\Gamma_{201}^{(2)}
  =-4$, $\Gamma_{030}^{(2)}=-3$, $\Gamma_{021}^{(2)}=-4$ and $\Gamma_{012}^{(2)}=-1$. 
  \doppelbildbox{7 true cm}{C0cf_aff.eps}{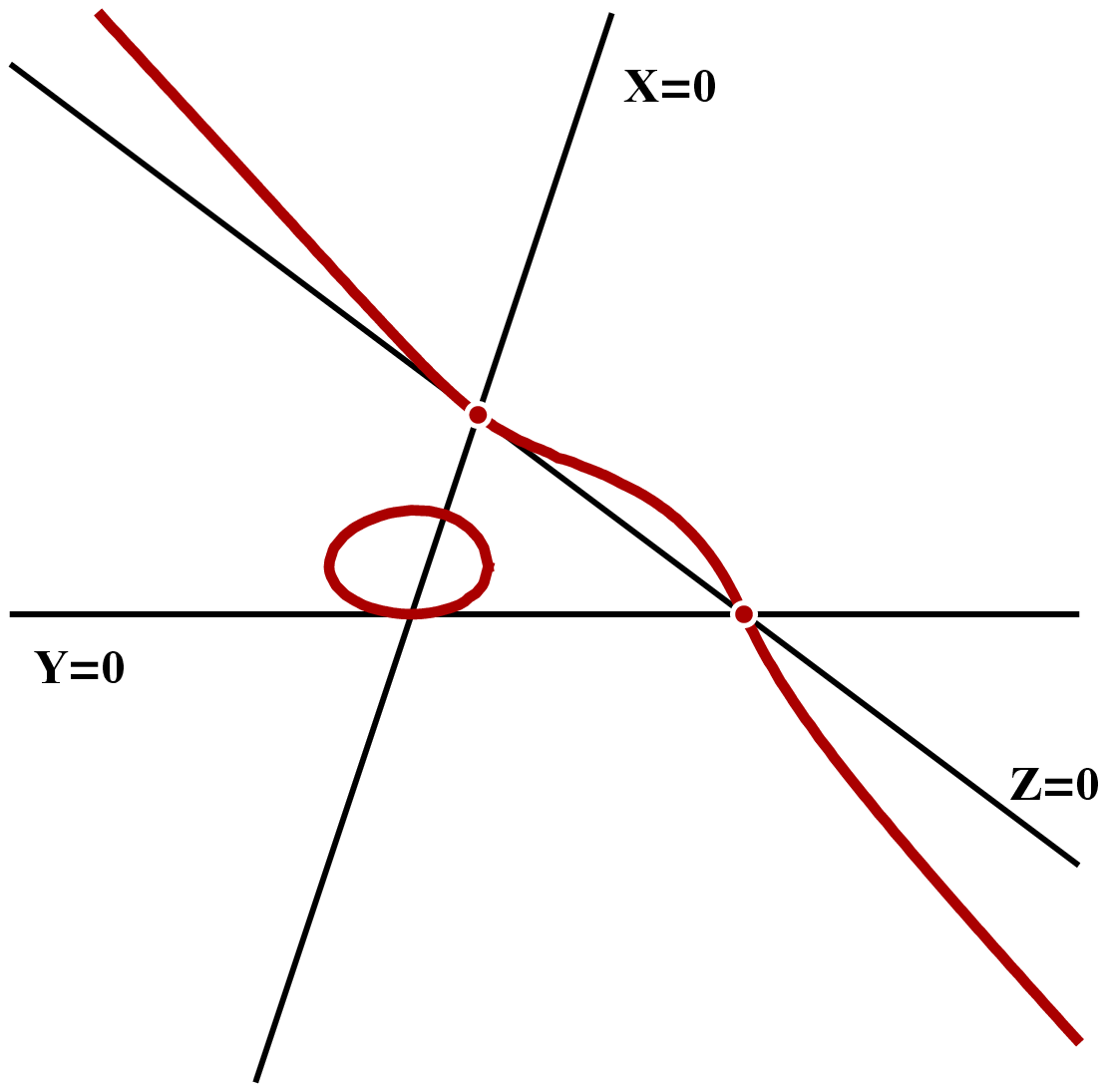}{{\bf Before:} 
  Affine view (left) and projective view (right) of the curve $C_{(0)}$.}
  \doppelbildbox{7 true cm}{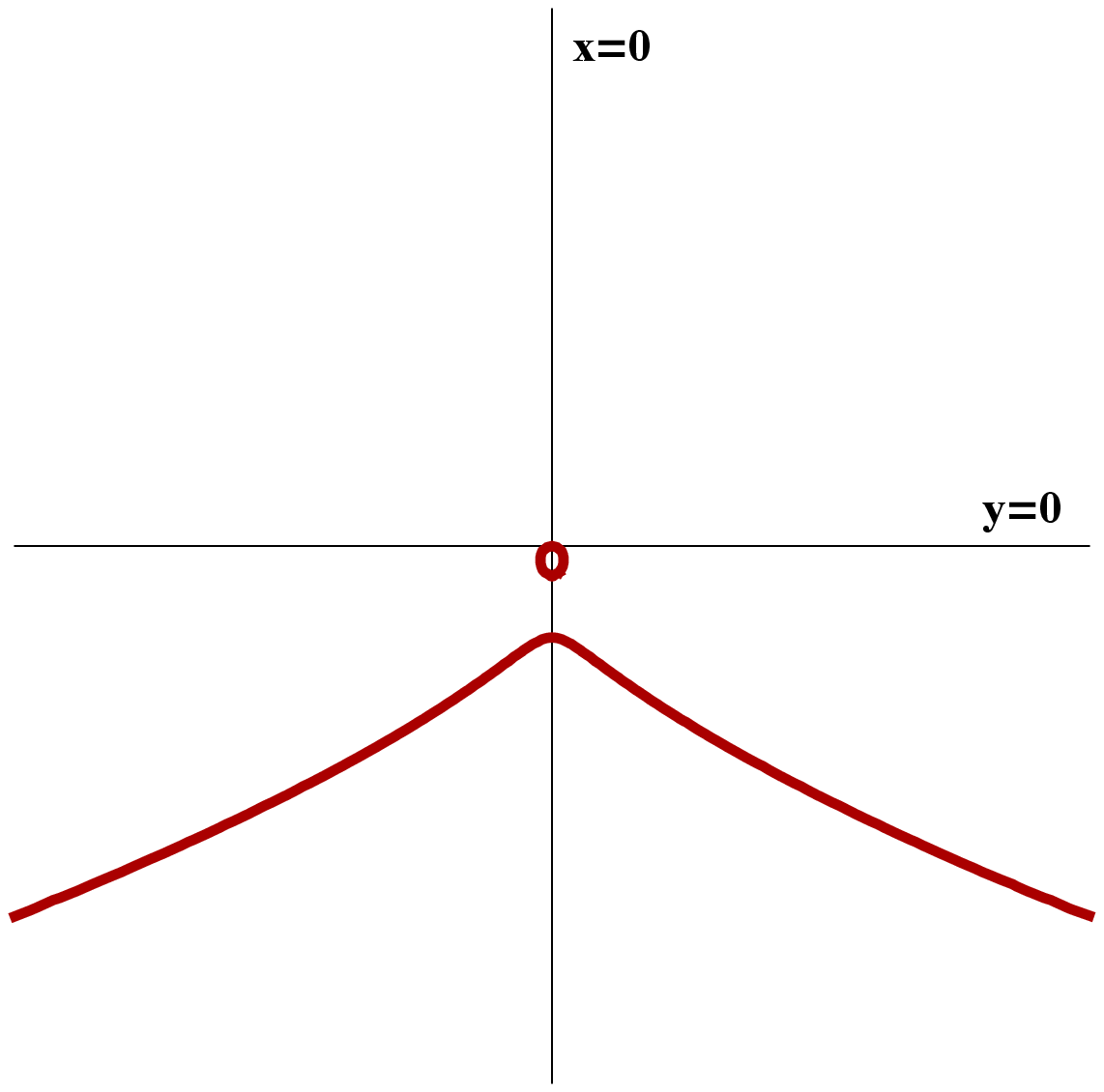}{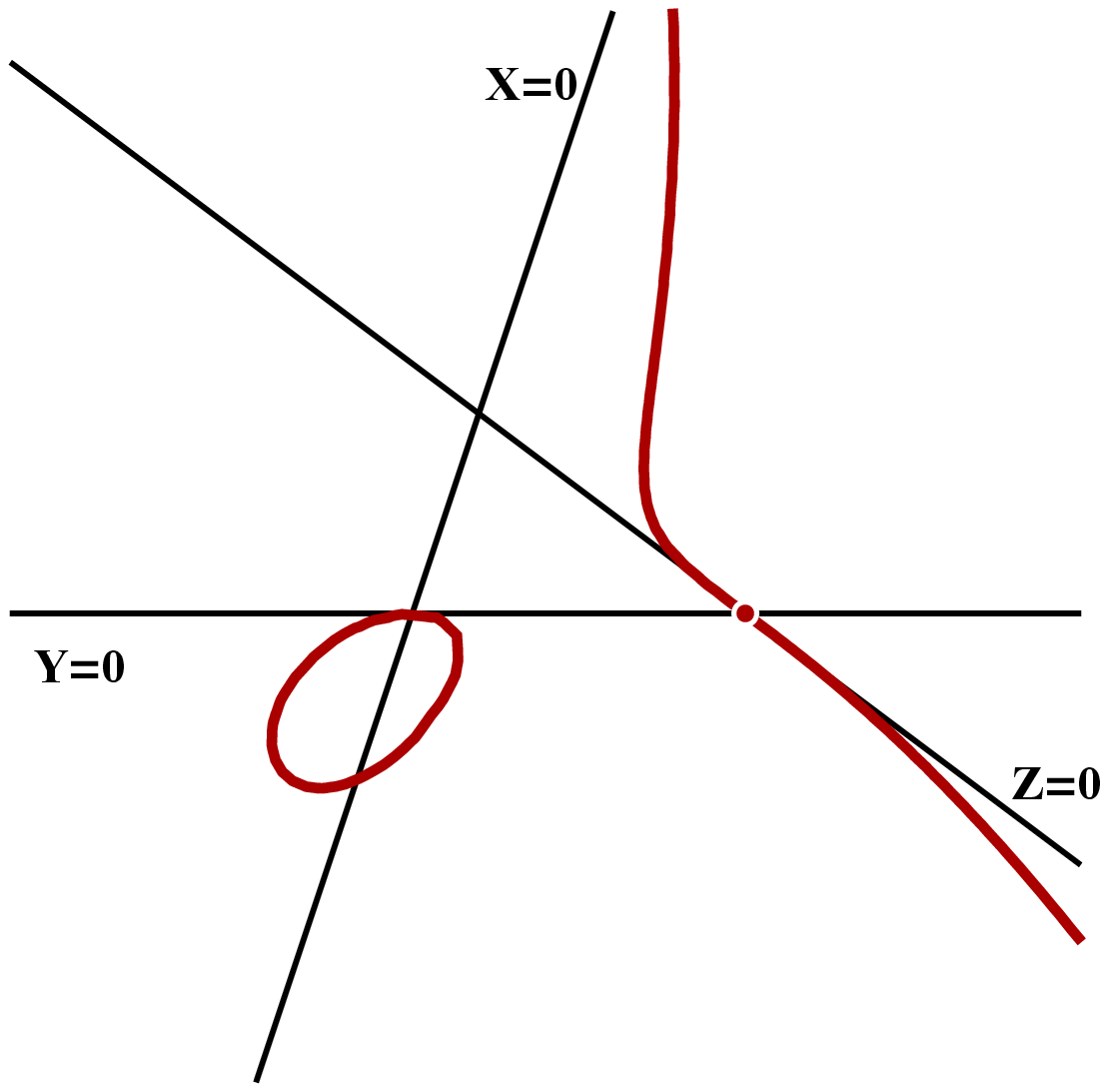}{}{{\bf After:} 
  Affine view (left) and projective view (right) of the curve $C_{(2)}$.}
  \eject   

{\bf Steps 3, 4 and 5:} These steps are superfluous, because the distinguished point  
   $p=(1,0,0)$ is already an inflection point; simply exchanging $X_2$ and $Y_2$ already 
   yields a general Weierstra\ss{} equation with the nonzero coefficients 
   $\Gamma_{300}^{(5)}=\Gamma_{030}^{(2)}=M(M-N)$, $\Gamma_{201}^{(5)}=\Gamma_{021}^{(2)}=2M-N$, 
   $\Gamma_{102}^{(5)}=\Gamma_{012}^{(2)}=1$ and $\Gamma_{021}^{(5)}=\Gamma_{201}^{(2)}=N^2$. 
   In the special case $(M,N)=(3,2)$ this means that $\Gamma_{300}^{(5)}=-3$, $\Gamma_{201}^{(5)}
   =-4$, $\Gamma_{102}^{(5)}=-1$ and $\Gamma_{021}^{(5)}=-4$. 
  \doppelbildbox{7 true cm}{C2cf_aff.eps}{C2cf_proj.eps}{{\bf Before:} 
  Affine view (left) and projective view (right) of the curve $C_{(2)}$.}
  \bigskip 
  \doppelbildbox{7 true cm}{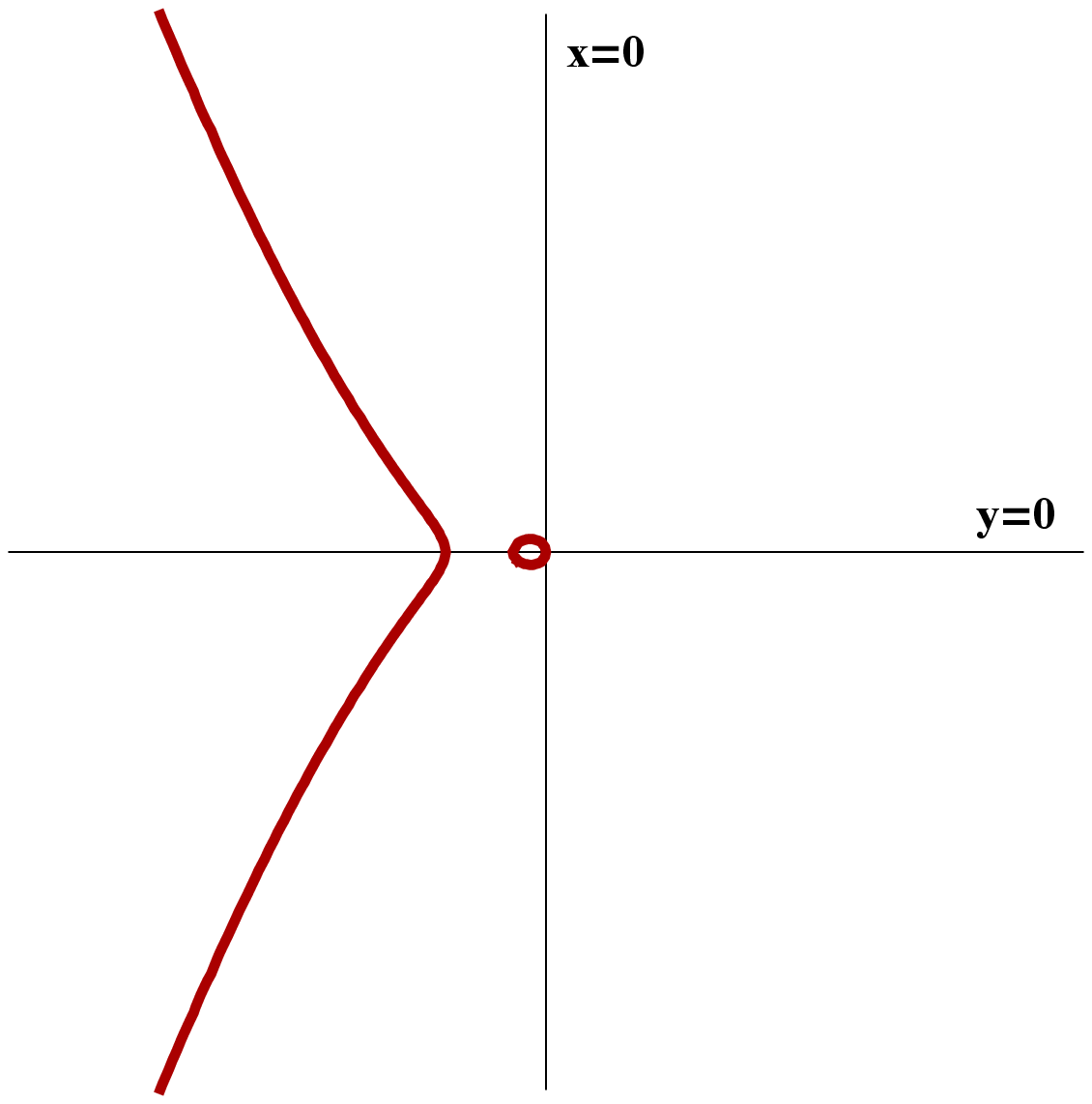}{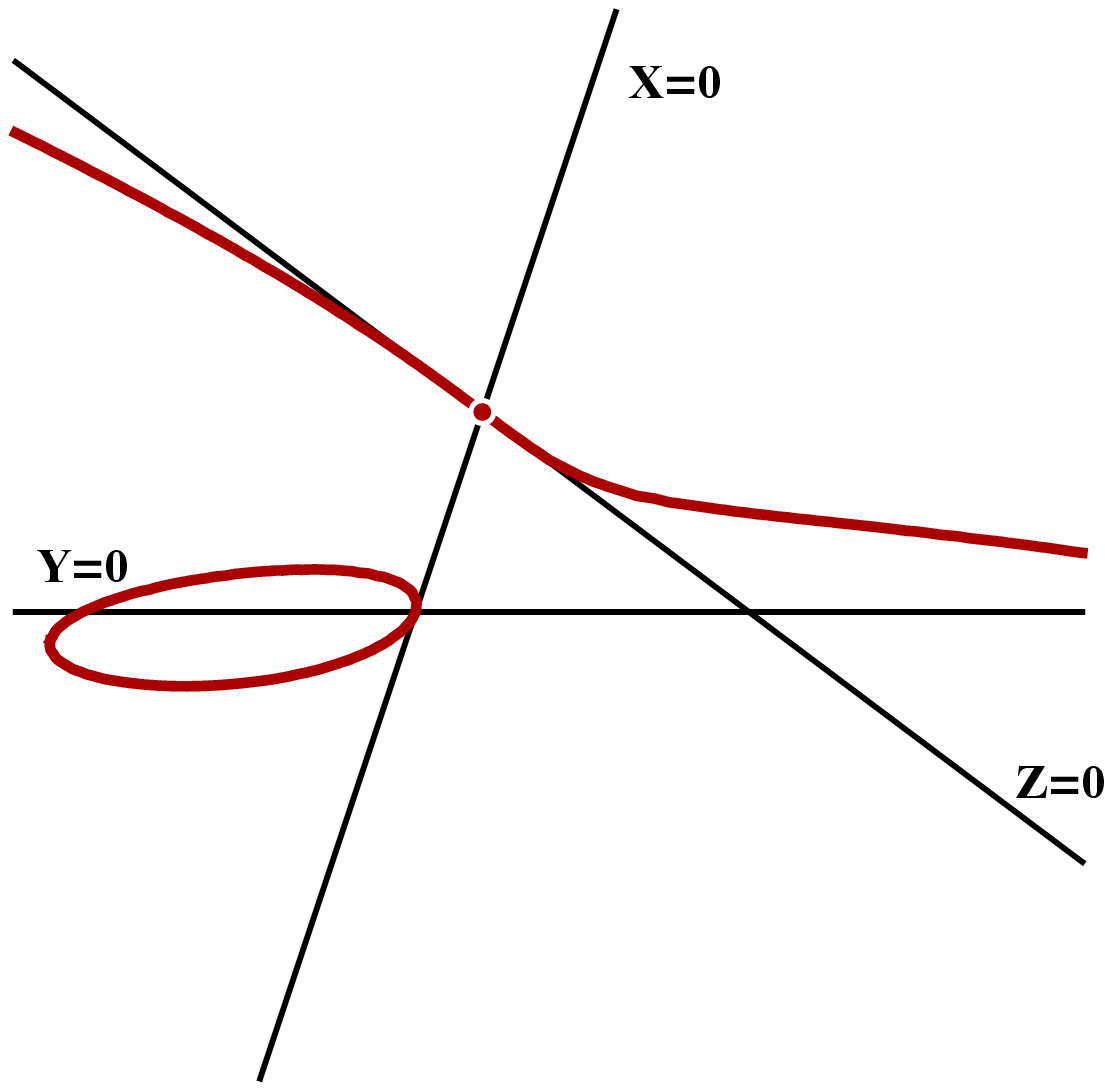}{}{{\bf After:} 
  Affine view (left) and projective view (right) of the curve $C_{(5)}$.}
  \vfill\eject 

{\bf Step 6:} The curve $C_{(6)}$ has the nonzero coefficients $\Gamma_{300}^{(6)}=4N^2M(M-N)$, 
   $\Gamma_{201}^{(6)}=4N^2(2M-N)$, $\Gamma_{021}^{(6)}=1$ and $\Gamma_{102}^{(6)}=4N^2$. 
   In the special case $(M,N)=(3,2)$ this means that $\Gamma_{300}^{(6)}=48$, 
   $\Gamma_{201}^{(6)}=64$, $\Gamma_{102}^{(6)}=16$ and $\Gamma_{021}^{(6)}=1$. 
   \doppelbildbox{7 true cm}{C5cf_aff.eps}{C5cf_proj.eps}{{\bf Before:} 
  Affine view (left) and projective view (right) of the curve $C_{(5)}$.}
  \bigskip 
  \doppelbildbox{7 true cm}{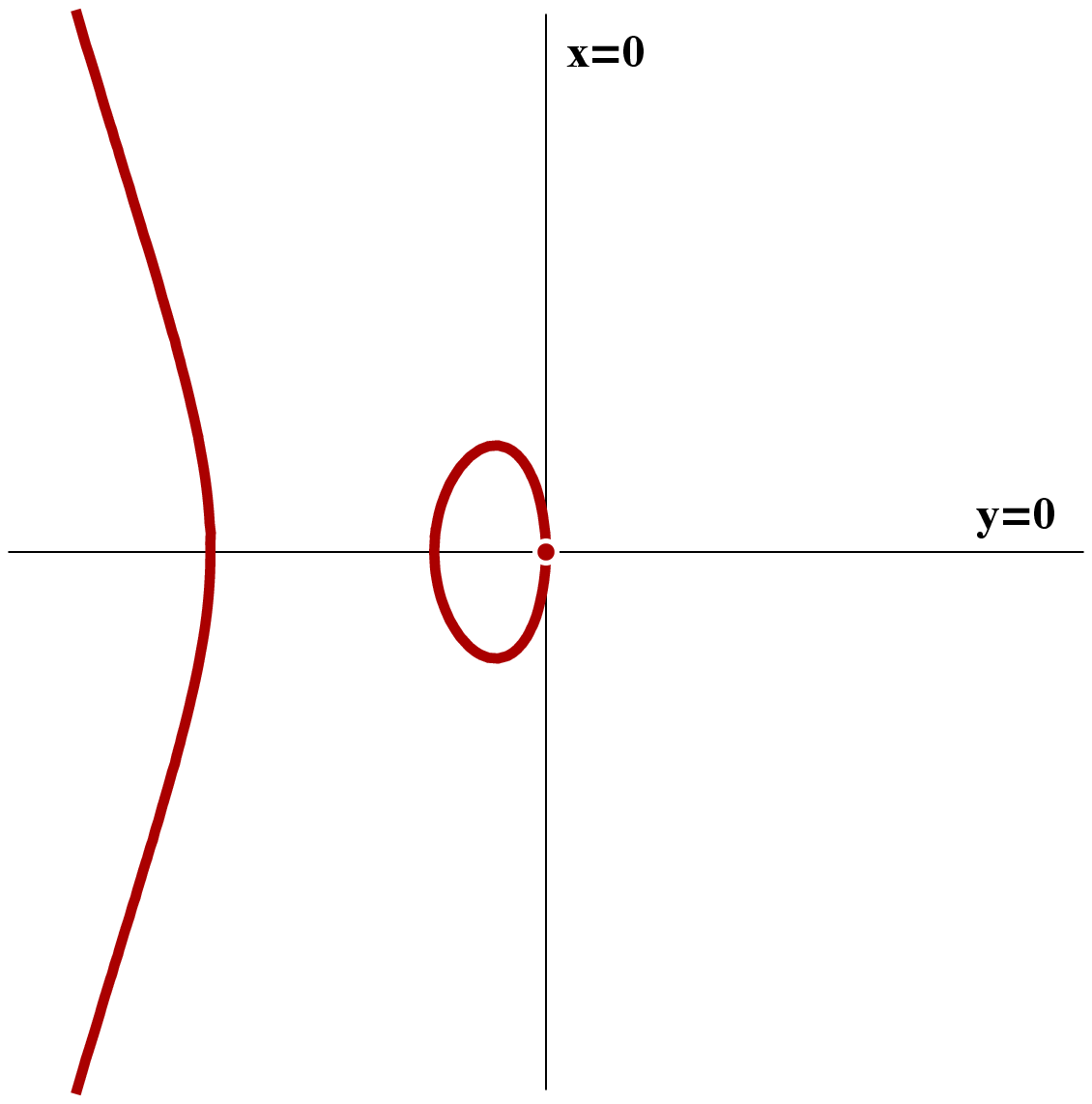}{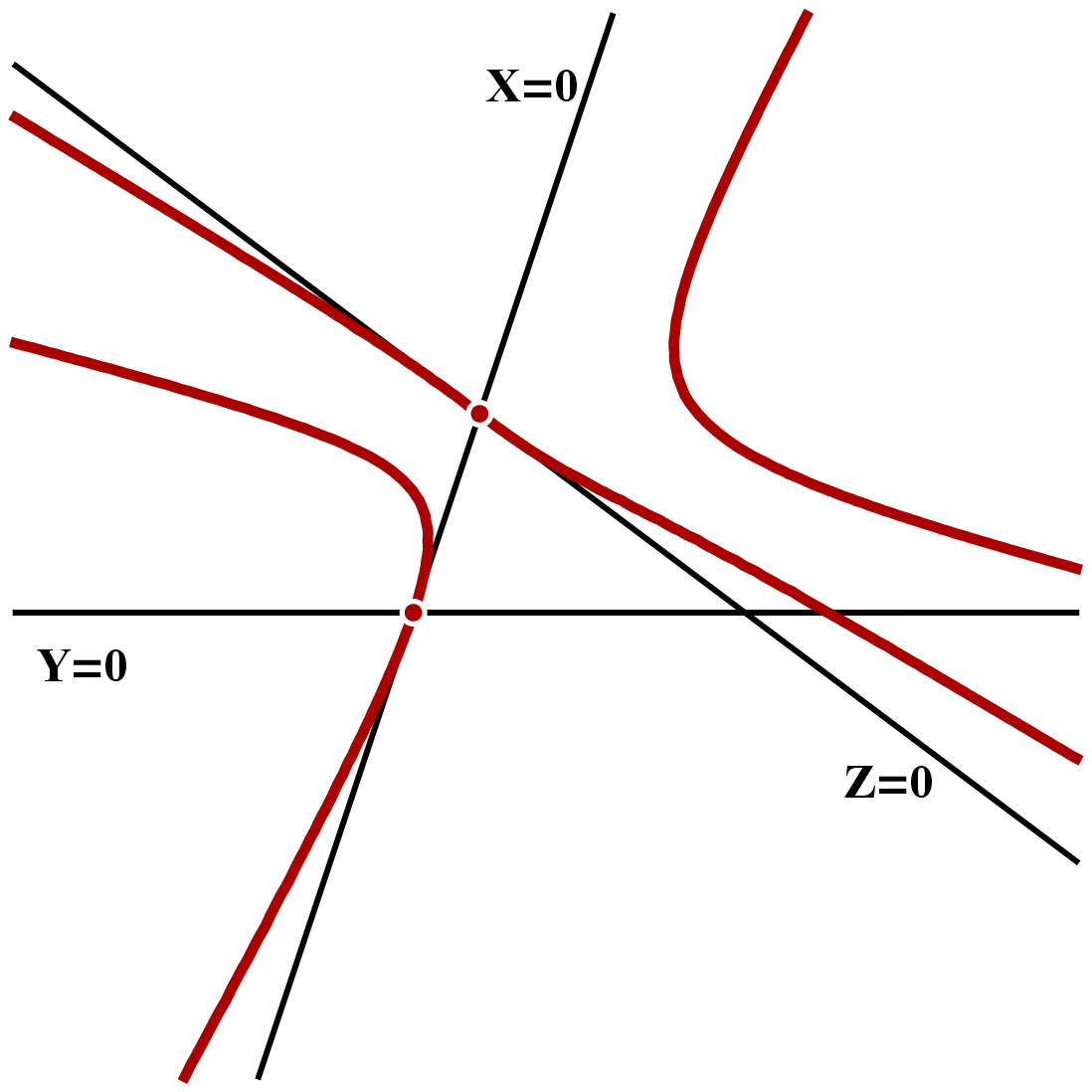}{}{{\bf After:} 
  Affine view (left) and projective view (right) of the curve $C_{(6)}$.}
  \vfill\eject 

{\bf Step 7:} We have $\delta=-4N^2M(M-N)$ and $\varphi=2N$. Hence the curve $C_{(7)}$ has the 
   nonzero coefficients $\Gamma_{300}^{(7)}=-1$, $\Gamma_{201}^{(7)}=2M-N$, $\Gamma_{102}^{(7)}
   =-M(M-N)$ and $\Gamma_{021}^{(7)}=1$ and hence, letting $(X,Y,Z):=(X_7,Y_7,Z_7)$ is given by 
   the equation $Y^2Z = X^3 - (2M-N)X^2Z + M(M-N)XZ^2 = X(X-M)\bigl( X-(M\! -\! N)Z\bigr)$. In 
   affine form, this reads $y^2 = x(x-M)\bigl(x- (M\! -\! N)\bigr)$. In the special case 
   $(M,N)=(3,2)$ this means that $\Gamma_{300}^{(7)}=-1$, $\Gamma_{201}^{(7)}=4$, $\Gamma_{102}^{(7)}
   =-3$ and $\Gamma_{021}^{(7)}=1$, and the curve $C_{(7)}$ is in affine form given by 
   $y^2=x(x-3)(x-1)$. 
  \doppelbildbox{6.8 true cm}{C6cf_aff.eps}{C6cf_proj.eps}{{\bf Before:} 
  Affine view (left) and projective view (right) of the curve $C_{(6)}$.} \vskip -0.1 true cm
  \doppelbildbox{6.8 true cm}{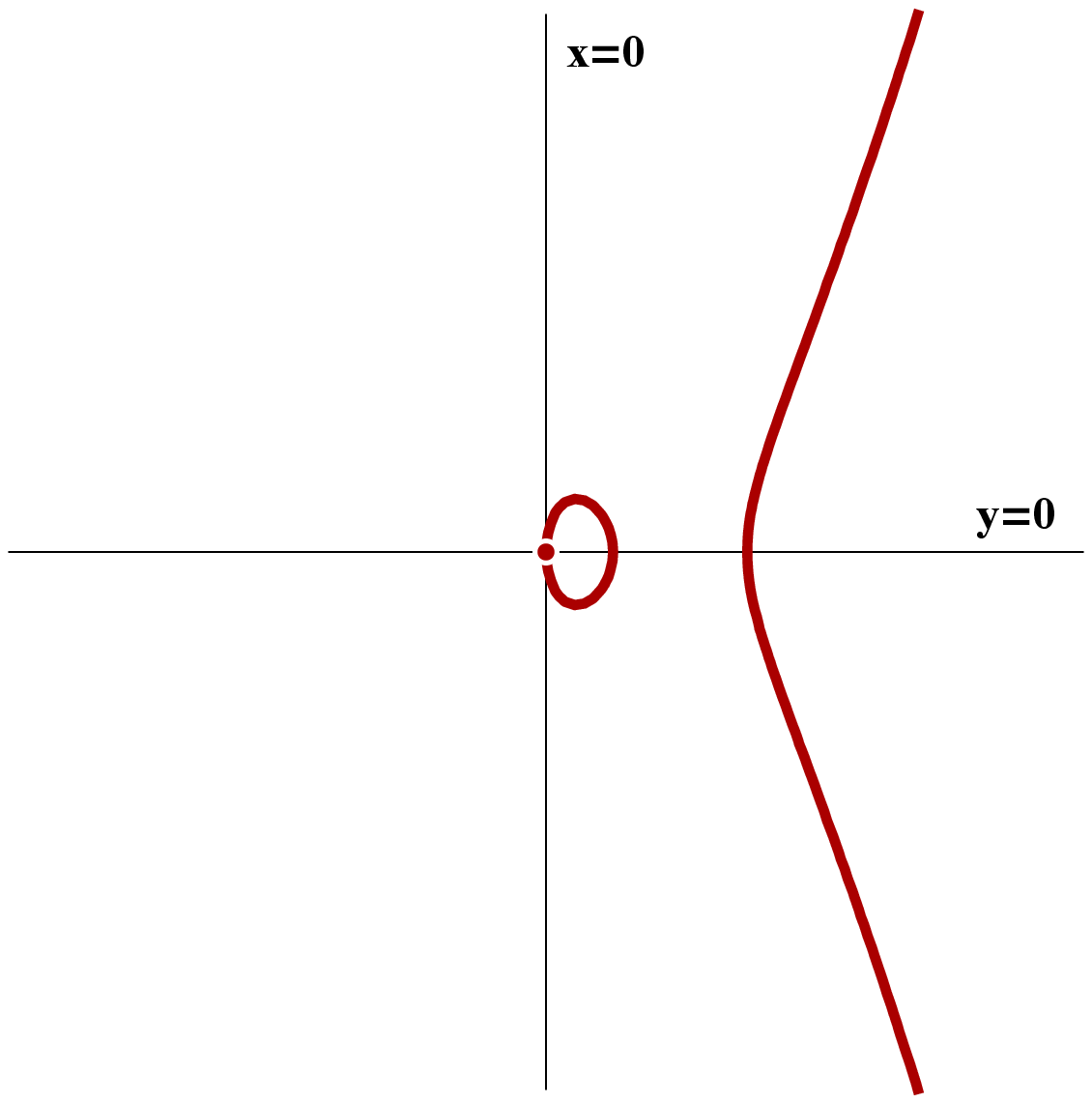}{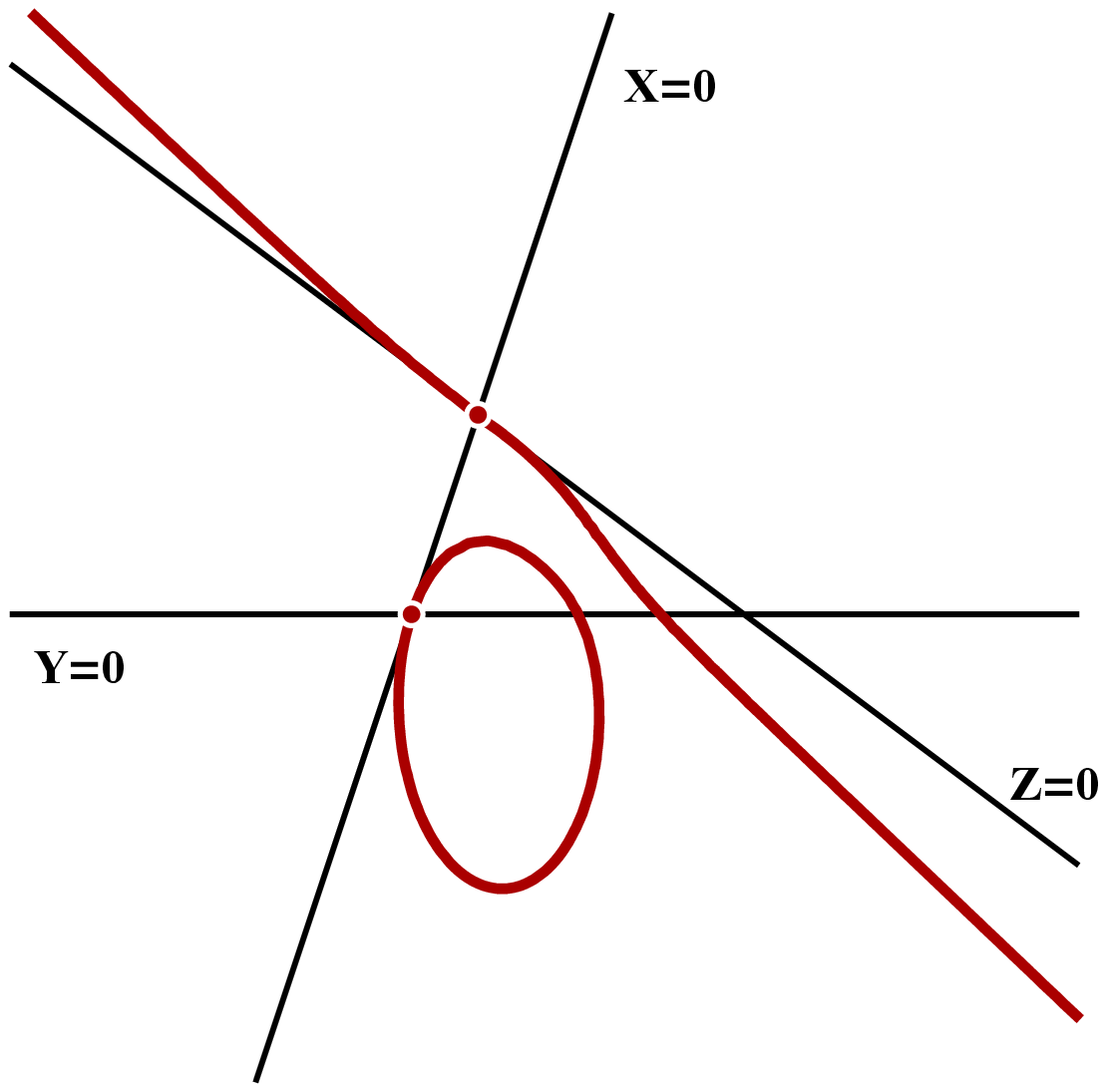}{}{{\bf After:} 
  Affine view (left) and projective view (right) of the curve $C_{(7)}$.}
  \eject 

We observe that the concatenation of the above coordinate transformations is given by 
$$\left[\matrix{X\cr Y\cr Z\cr}\right]\ =\ \left[\matrix{0&-M(M\! -\! N)&\phantom{-}0&M(M\! -\! N)\cr 
   MN(M\! -\! N)&0&\phantom{-}0&0\cr 0&-(M\! -\! N)&-N&M\cr}\right]\left[\matrix{X_0\cr X_1\cr X_2\cr 
   X_3\cr}\right].$$
Thus the isomorphism $\Phi$ from the quadric intersection $Q_{M,N}$ to the elliptic curve $E_{M,N}$ 
with the affine equation $y^2 = x(x-M)\bigl( x-(M\! -\! N)\bigr)$ is accomplished by a simple 
linear map. The original quadric intersection $Q_{M,N}$ possesses four trivial points with the 
projective coordinates $(0,1,\pm 1,\pm 1)$. These are mapped by $\Phi$ to four rational points on 
the curve $E_{M,N}$, namely 
$$\eqalign{\Phi(0,1,1,1)\ &=\ (0,1,0)\quad\hbox{(point at infinity)},\cr 
           \Phi(0,1,1,-1)\ &=\ (M-N,0,1), \cr  
           \Phi(0,1,-1,1)\ &=\ (0,0,1),\cr 
           \Phi(0,1,-1,-1)\ &=\ (M,0,1).\cr}$$ 
Not surprisingly, these are exactly the $2$-torsion points on $E_{M,N}$. \vfill\eject 

\leftline{\titlefont 5. Example: Rational squares in arithmetic progression} \smallskip 
In the special situation of the $(k,\ell,m)$-problem which was one of the motivations of this paper 
mentioned in the introduction, we have $A=\hbox{\rm diag}(k+\ell, -k,-\ell, 0)$, $B=\hbox{\rm 
diag}(-m,m+\ell, 0,-\ell)$ and $x=(1,1,1,1)$. The cubic curve $C_{(0)}$ to which the quadric 
intersection $Q_1\cap Q_2$ is transformed is given by $\sum_{i+j+k=3} \Gamma_{ijk} X^iY^jZ^k=0$ 
where the only nonzero coefficients $\Gamma_{ijk} = \Gamma^{(0)}_{ijk}$ are as follows: 
$$\eqalign{\Gamma_{210}^{(0)}\ &=\ -(k+\ell +m),\cr 
           \Gamma_{201}^{(0)}\ &=\ m,\cr 
           \Gamma_{120}^{(0)}\ &=\ k+\ell +m,\cr 
           \Gamma_{102}^{(0)}\ &=\ -m,\cr 
           \Gamma_{021}^{(0)}\ &=\ -(\ell +m),\cr 
           \Gamma_{012}^{(0)}\ &=\ \ell +m.\cr}$$ 
The distinguished rational point is $p^{(0)}=(\ell\! +\! m, m, k\! +\! \ell\! +\! m)$. We now follow 
the general procedure described in the previous paragraph, exemplifying all results for the case 
$(k,\ell, m) = (2,3,5)$ in which we have the coefficients 
$$\eqalign{
  &C_{300}=0,\quad C_{210}=-10,\quad C_{201}=5,\quad C_{120}=10,\quad C_{111}=0,\cr 
  &C_{102}=-5,\quad C_{030}=0,\quad C_{021}=-8,\quad C_{012}=8,\quad C_{003}=0\cr}$$ 
and the distinguished rational point $(8,5,10)$.\smallskip 

\doppelbildbox{7 true cm}{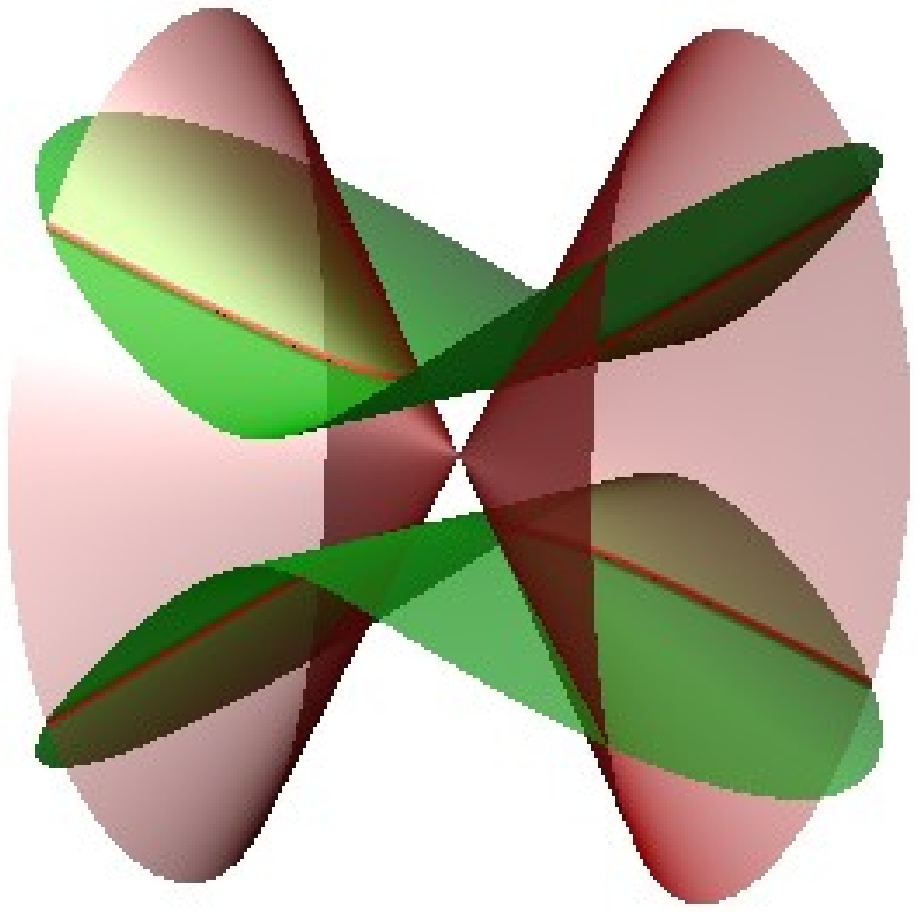}{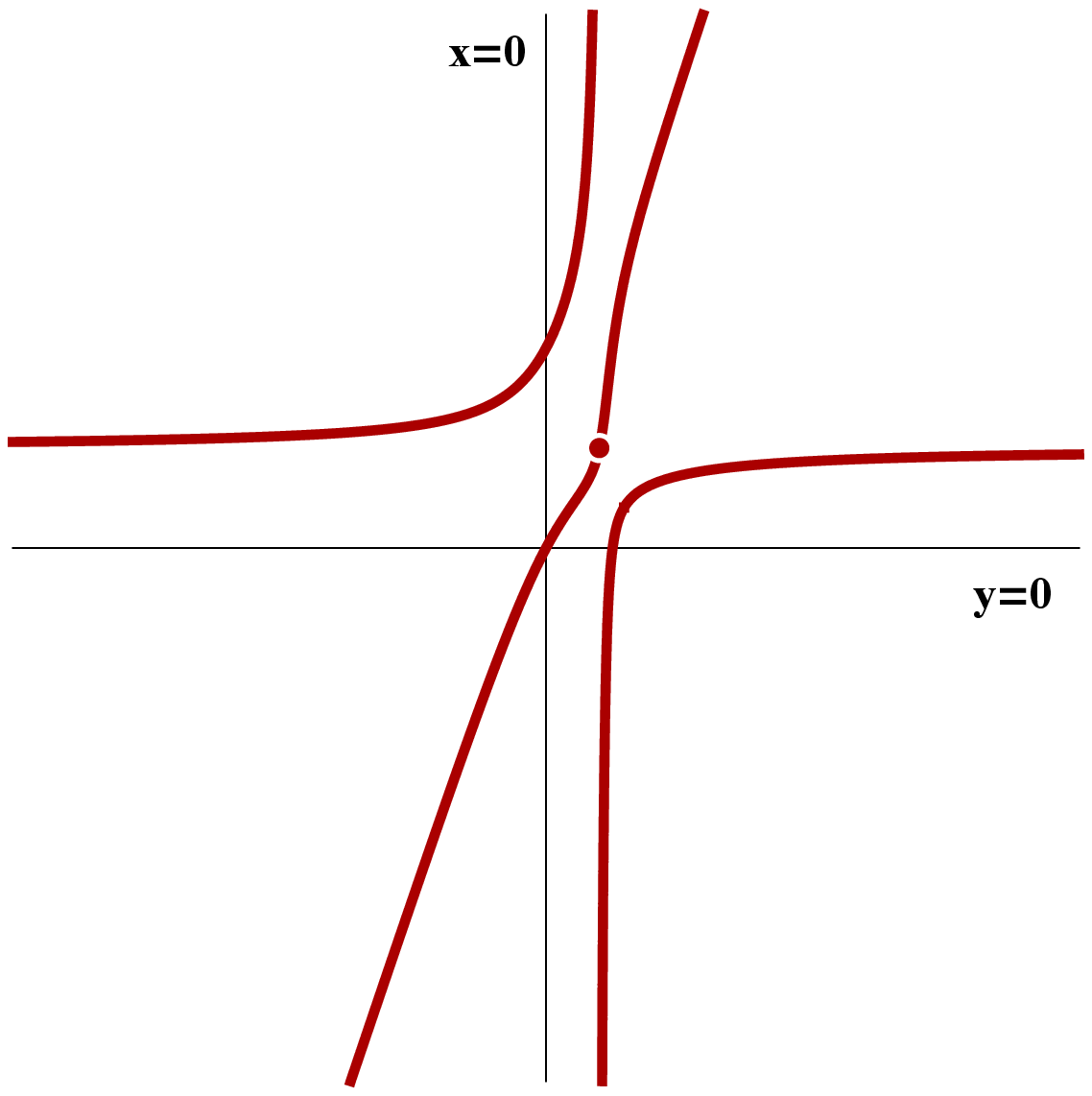}{{\bf Left:} Affine view 
$X_3=1$ of the intersection of the quadrics $(k+\ell)X_0^2 - kX_1^2 - \ell X_2^2 = 0$ and 
$-mX_0^2 + (m+\ell) X_1^2 - \ell X_3^2 = 0$, shown for $(k,\ell, m) = (2,3,5)$. {\bf Right:} 
Affine view of the cubic $C_{(0)}$ to which this quadric intersection is initially transformed.} 
\vfill\eject 

{\bf Step 1:} Translating the distinguished rational point to $(1,0,0)$. 
  \doppelbildbox{7 true cm}{C0klm_aff.eps}{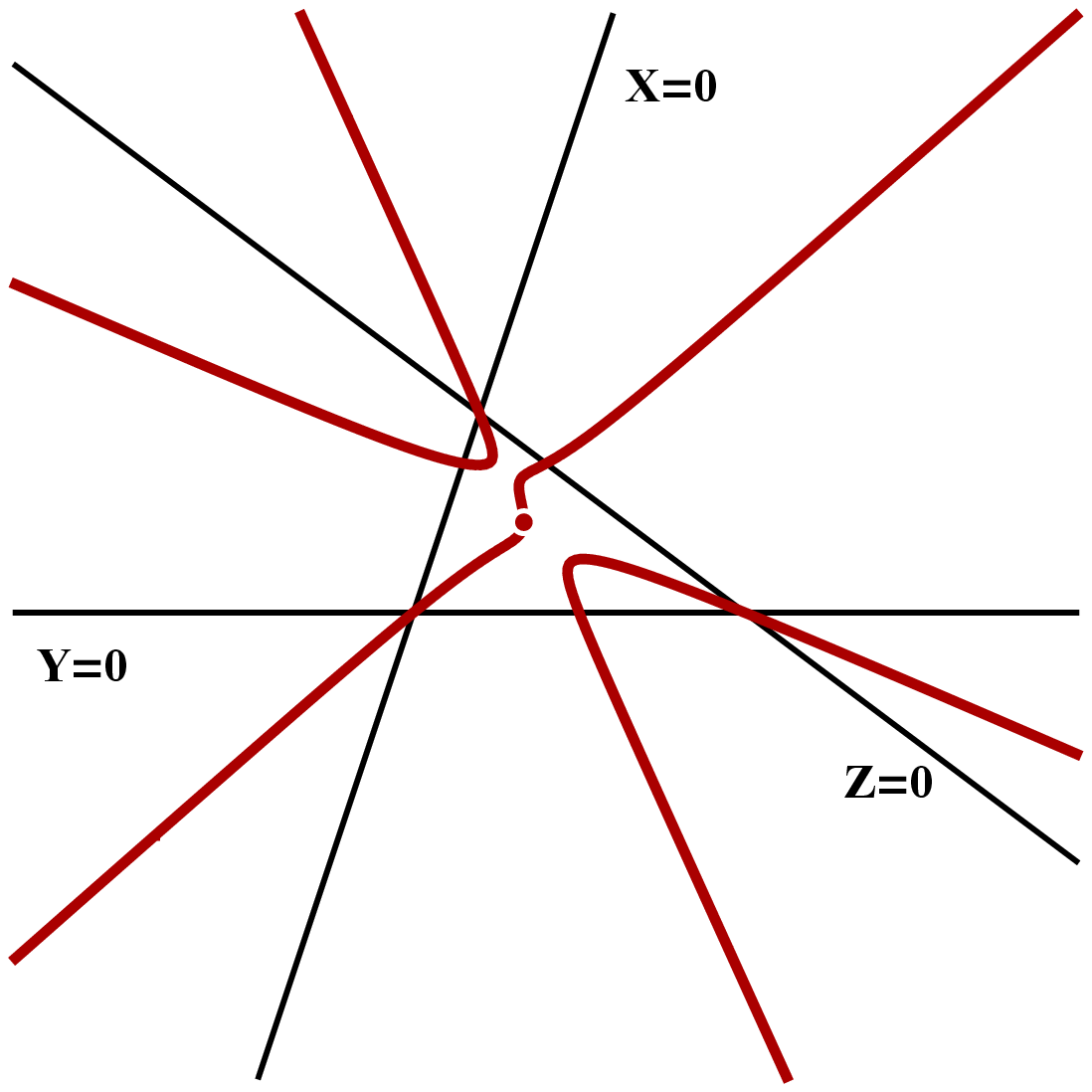}{{\bf Before:} 
  Affine view (left) and projective view (right) of the curve $C_{(0)}$.}
  \bigskip 
  \doppelbildbox{7 true cm}{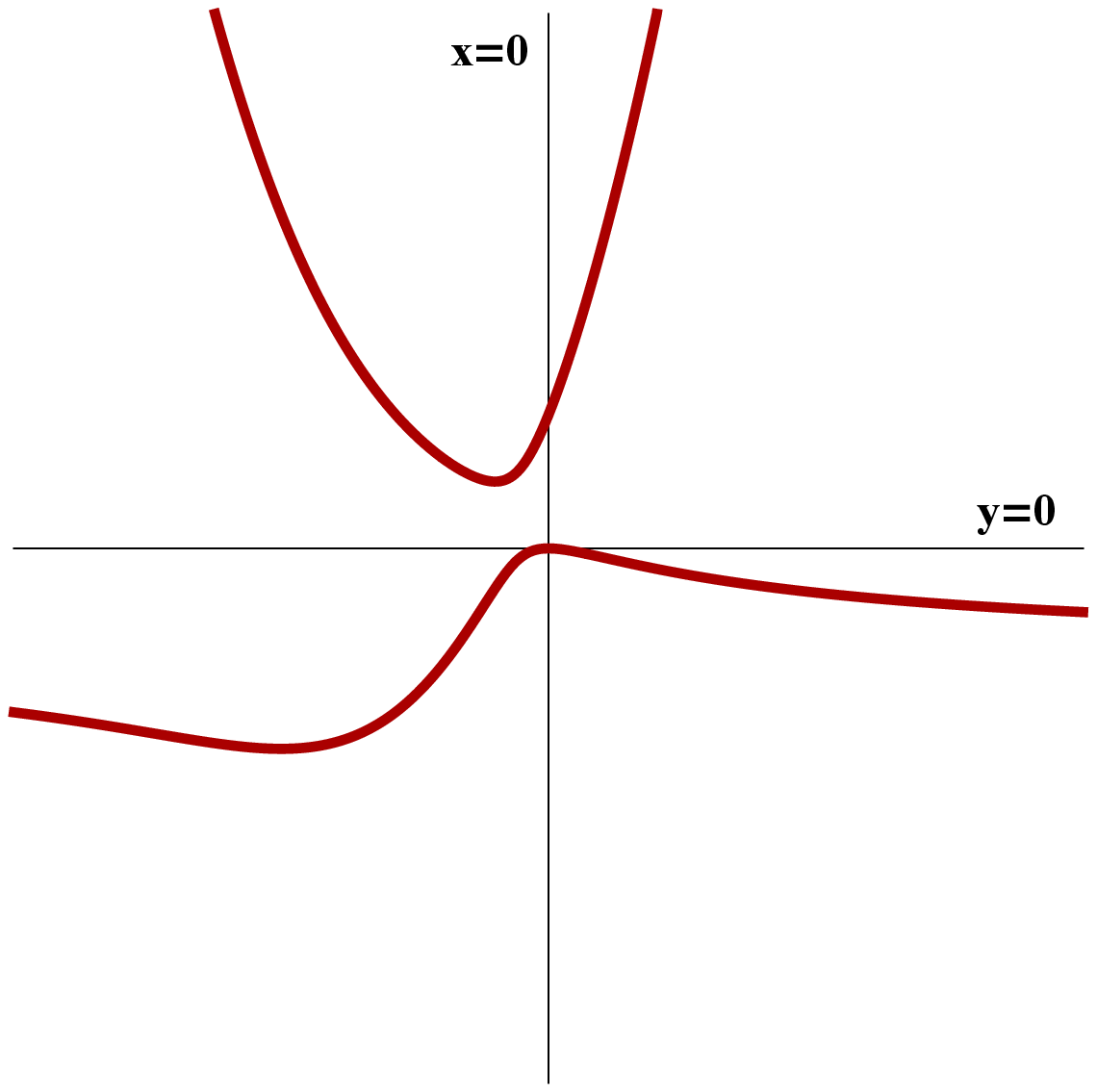}{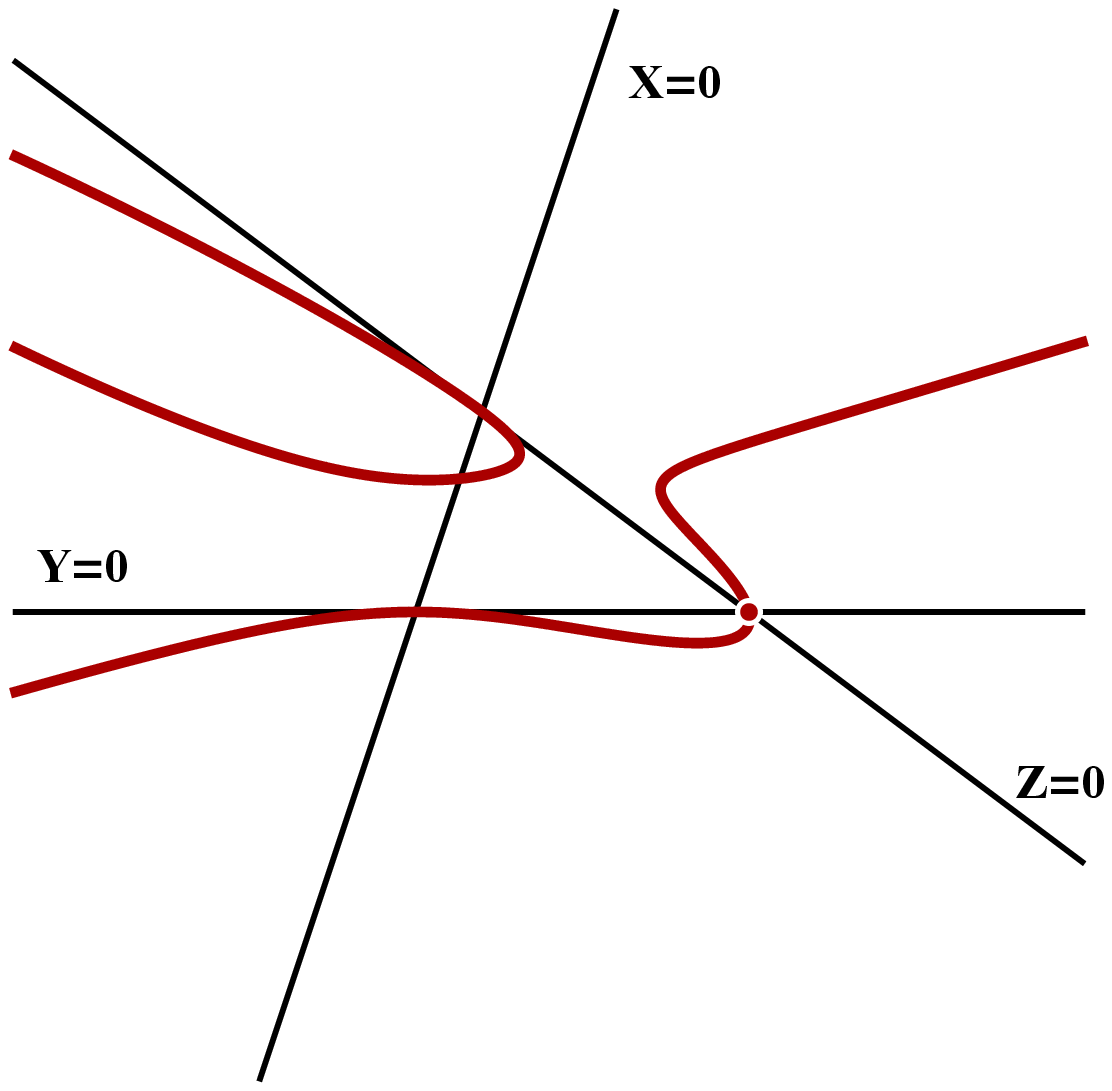}{}{{\bf After:} 
  Affine view (left) and projective view (right) of the curve $C_{(1)}$.}
  \vfill\eject 
  The curve $C_{(1)}$ has the nonzero coefficients 
  $$\eqalign{
    \Gamma_{210}^{(1)}\ &=\ k(k+\ell +m),\cr 
    \Gamma_{201}^{(1)}\ &=\ \ell m,\cr 
    \Gamma_{111}^{(1)}\ &=\ 2(k+\ell)(\ell +m),\cr 
    \Gamma_{021}^{(1)}\ &=\ -(\ell +m)^2,\cr 
    \Gamma_{012}^{(1)}\ &=\ (\ell +m)^2.\cr}$$ 
  In the special case $(k,\ell, m) = (2,3,5)$ this yields 
  $$\eqalign{
    &\Gamma_{300}^{(1)}=0,\quad \Gamma_{210}^{(1)}=20,\quad \Gamma_{201}^{(1)}=15,
        \quad \Gamma_{120}^{(1)}=0,\quad \Gamma_{111}^{(1)}=80,\cr 
    &\Gamma_{102}^{(1)}=0,\quad \Gamma_{030}^{(1)}=0,\quad \Gamma_{021}^{(1)}=-64, 
        \quad \Gamma_{012}^{(1)}=64,\ \Gamma_{003}^{(1)}=0.\cr}$$ 
  The distinguished rational point after this step is $(1,0,0)$.
  \vfill\eject 

{\bf Step 2:} Adapting the tangent of the distinguished rational point $(1,0,0)$. 
  \doppelbildbox{7 true cm}{C1klm_aff.eps}{C1klm_proj.eps}{{\bf Before:} 
  Affine view (left) and projective view (right) of the curve $C_{(1)}$.}
  \bigskip 
  \doppelbildbox{7 true cm}{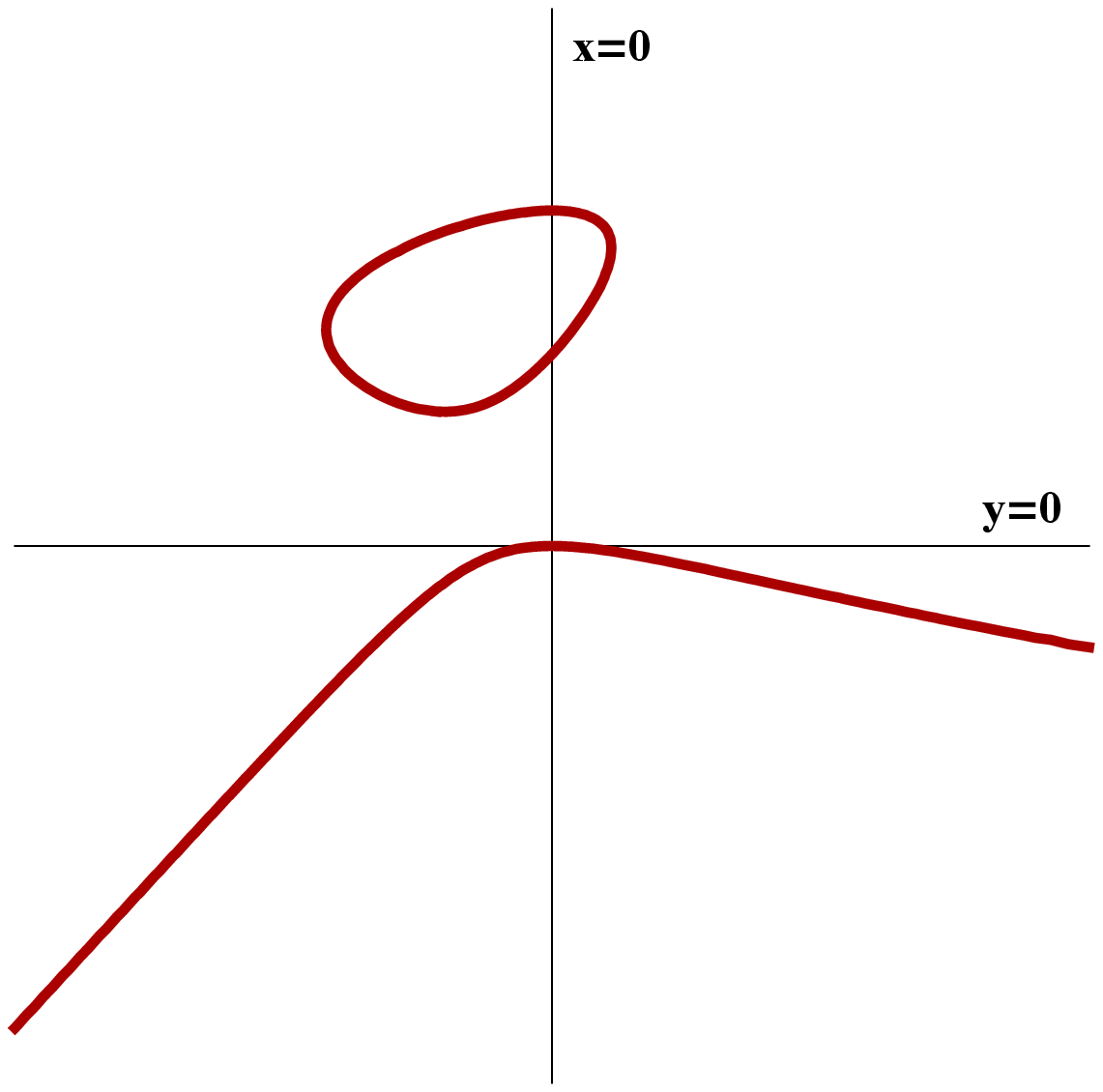}{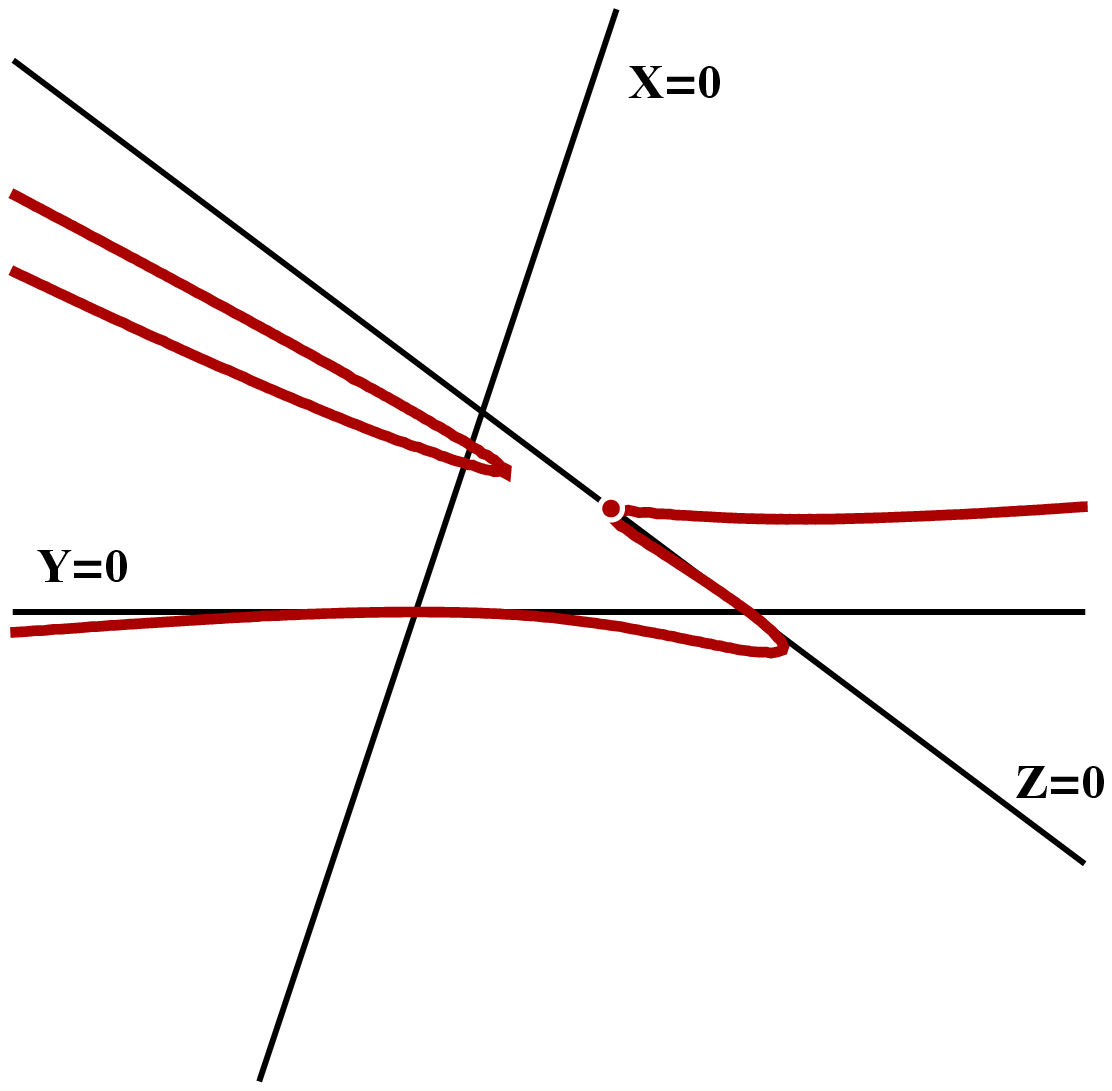}{}{{\bf After:} 
  Affine view (left) and projective view (right) of the curve $C_{(2)}$.}
  \vfill\eject 
  The curve $C_{(2)}$ has the nonzero coefficients
  $$\eqalign{
    \Gamma_{201}^{(2)}\ &=\ \ell^2m^2,\cr 
    \Gamma_{120}^{(2)}\ &=\ -2k\ell m(k+\ell)(\ell +m)(k+\ell +m),\cr 
    \Gamma_{111}^{(2)}\ &=\ 2\ell m(k+\ell )(\ell +m),\cr 
    \Gamma_{030}^{(2)}\ &=\ k(k+\ell)(k+m)(k+\ell +m)(\ell +m)^2,\cr 
    \Gamma_{021}^{(2)}\ &=\ -(\ell +m)^2\bigl( 2k(k\! +\! \ell\! +\! m) + \ell m\bigr),\cr 
    \Gamma_{012}^{(2)}\ &=\ (\ell +m)^2.\cr}$$ 
  The second point of intersection of $C_{(2)}$ with the tangent $Z_2=0$ is 
  $$p^{(2)}\ =\ \bigl( k\ell + (k\! +\!\ell\! +\!  m)m,\, 2\ell m,\, 0\bigr).$$ 
  In our example we have 
  $$\eqalign{
    &\Gamma_{300}^{(2)}=0,\quad \Gamma_{210}^{(2)}=0,\quad \Gamma_{201}^{(2)}=45,
         \quad \Gamma_{120}^{(2)}=-960,\quad \Gamma_{111}^{(2)}=240,\cr 
    &\Gamma_{102}^{(2)}=0,\quad \Gamma_{030}^{(2)}=1792,\quad \Gamma_{021}^{(2)}=-704, 
         \quad \Gamma_{012}^{(2)}=64,\quad \Gamma_{003}^{(2)}=0\cr}$$ 
  and $p^{(2)}=(28,15,0)$. 
  \vfill\eject 

{\bf Step 3:} Arranging $(0,1,0)$ to be the second intersection point of the tangent at $(1,0,0)$ 
  with the cubic. 
  \doppelbildbox{7 true cm}{C2klm_aff.eps}{C2klm_proj.eps}{{\bf Before:} 
  Affine view (left) and projective view (right) of the curve $C_{(2)}$.}
  \bigskip 
  \doppelbildbox{7 true cm}{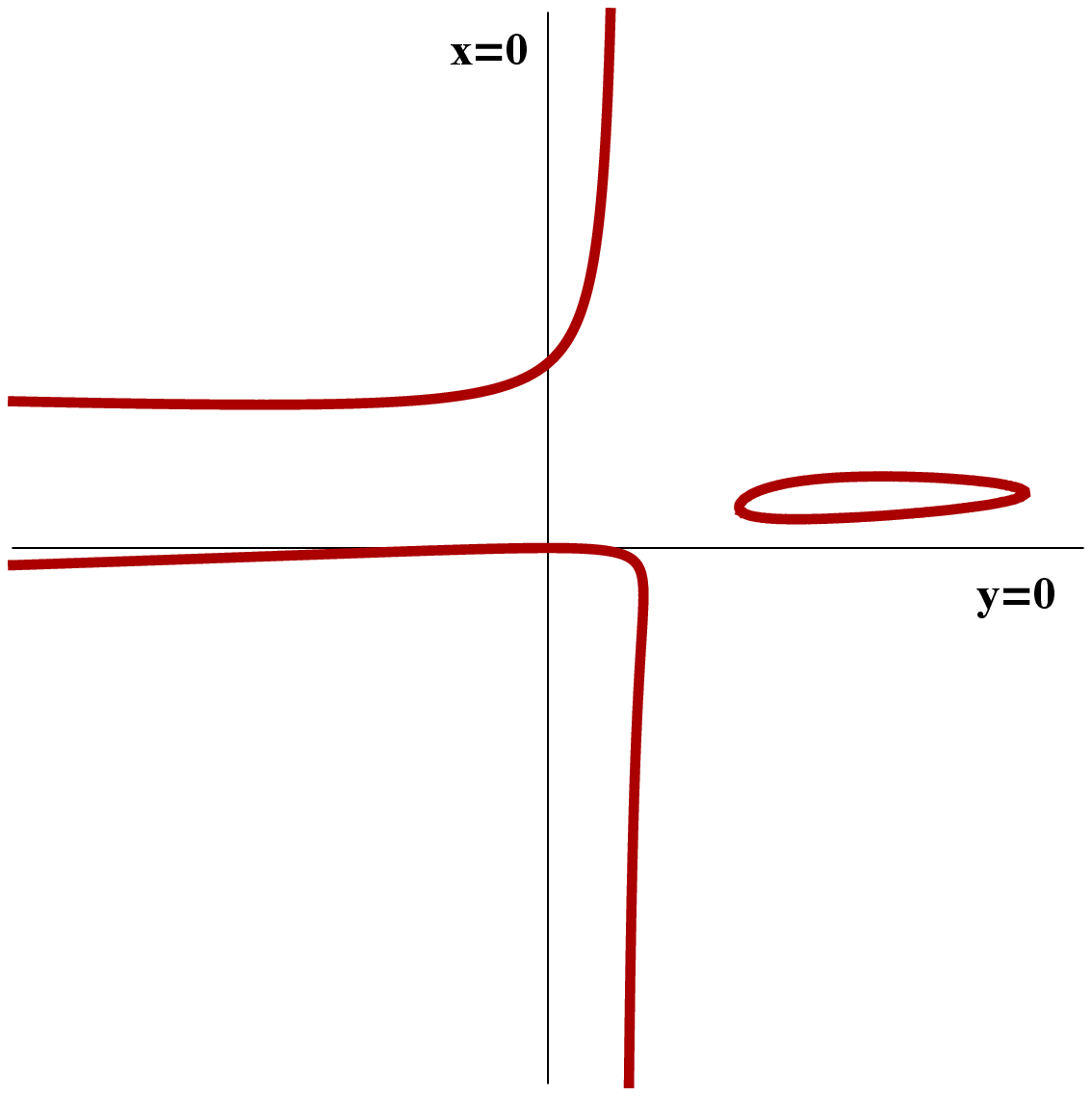}{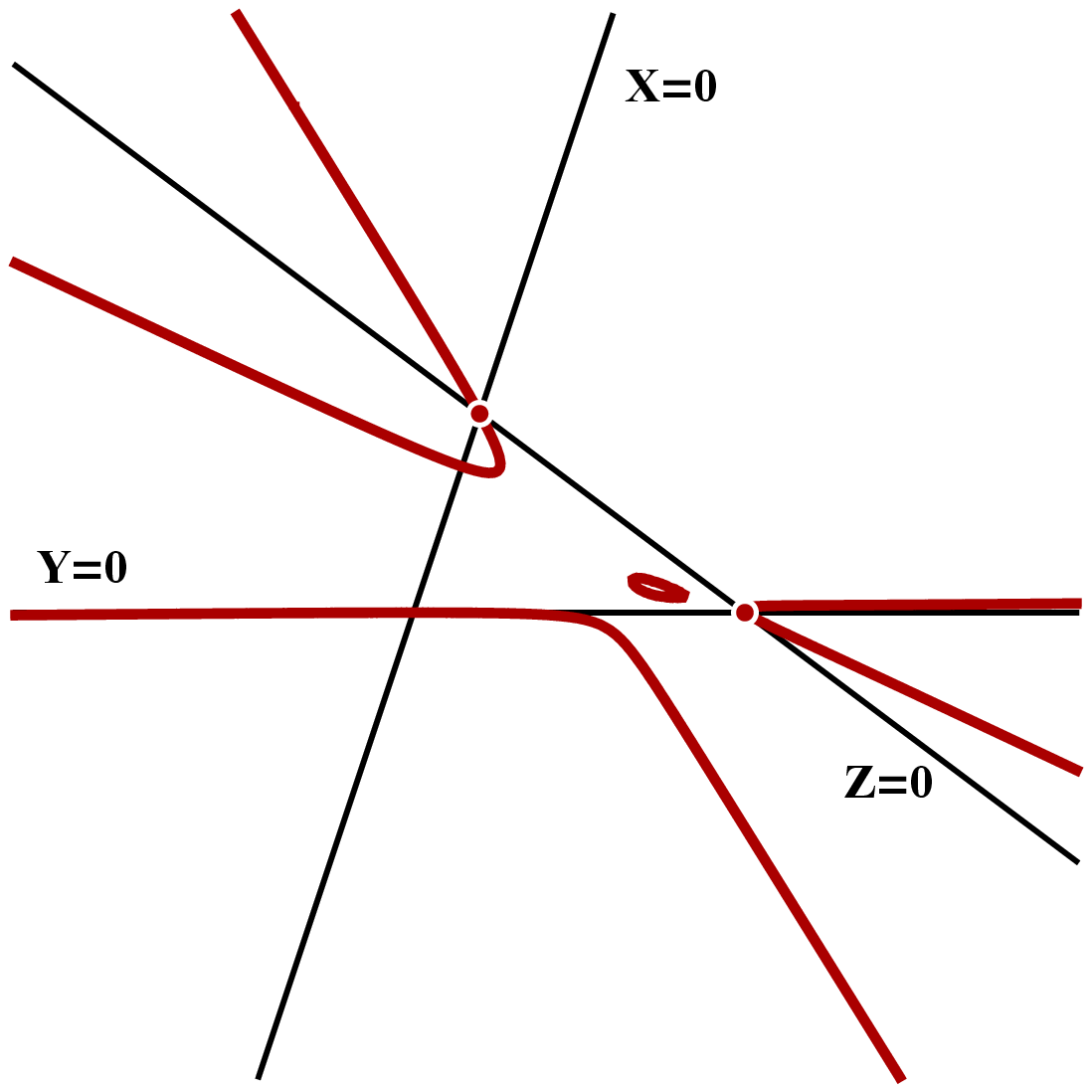}{}{{\bf After:} 
  Affine view (left) and projective view (right) of the curve $C_{(3)}$.}
  \vfill\eject 
  The curve $C_{(3)}$ has the nonzero coefficients
  $$\eqalign{
    \Gamma_{201}^{(3)}\ &=\ 1,\cr 
    \Gamma_{120}^{(3)}\ &=\ 4k(k+\ell )(\ell +m)(k+\ell +m),\cr 
    \Gamma_{111}^{(3)}\ &=\ -2(\ell +m)(3k+2\ell +m),\cr 
    \Gamma_{021}^{(3)}\ &=\ -(\ell +m)^2(3k^2+4k\ell +2km-m^2), \cr 
    \Gamma_{012}^{(3)}\ &=\ 4(\ell +m)^2. \cr}$$  
  In our example this yields
  $$\eqalign{
    &\Gamma_{300}^{(3)}=0,\quad \Gamma_{210}^{(3)}=0,\quad \Gamma_{201}^{(3)}=1,
        \quad \Gamma_{120}^{(3)}=320,\quad \Gamma_{111}^{(3)}=-136,\cr 
    &\Gamma_{102}^{(3)}=0,\quad \Gamma_{030}^{(3)}=0,\quad \Gamma_{021}^{(3)}=-496, 
        \quad \Gamma_{012}^{(3)}=320,\quad \Gamma_{003}^{(3)}=0.\cr}$$ 
  \vfill\eject 

{\bf Step 4:} Adapting the tangent at the point $(0,1,0)$. 
  \doppelbildbox{7 true cm}{C3klm_aff.eps}{C3klm_proj.eps}{{\bf Before:} 
  Affine view (left) and projective view (right) of the curve $C_{(3)}$.}
  \bigskip 
  \doppelbildbox{7 true cm}{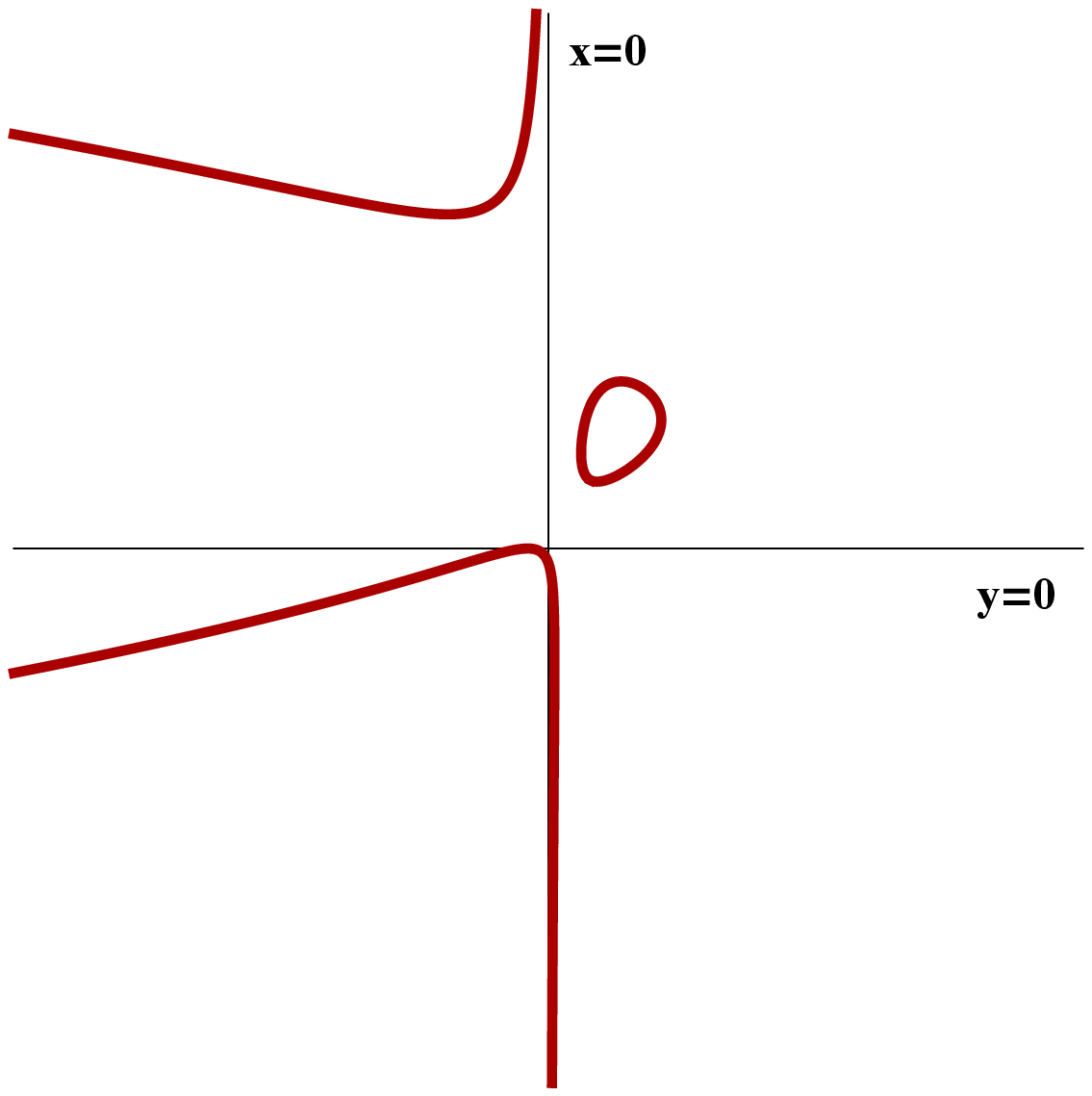}{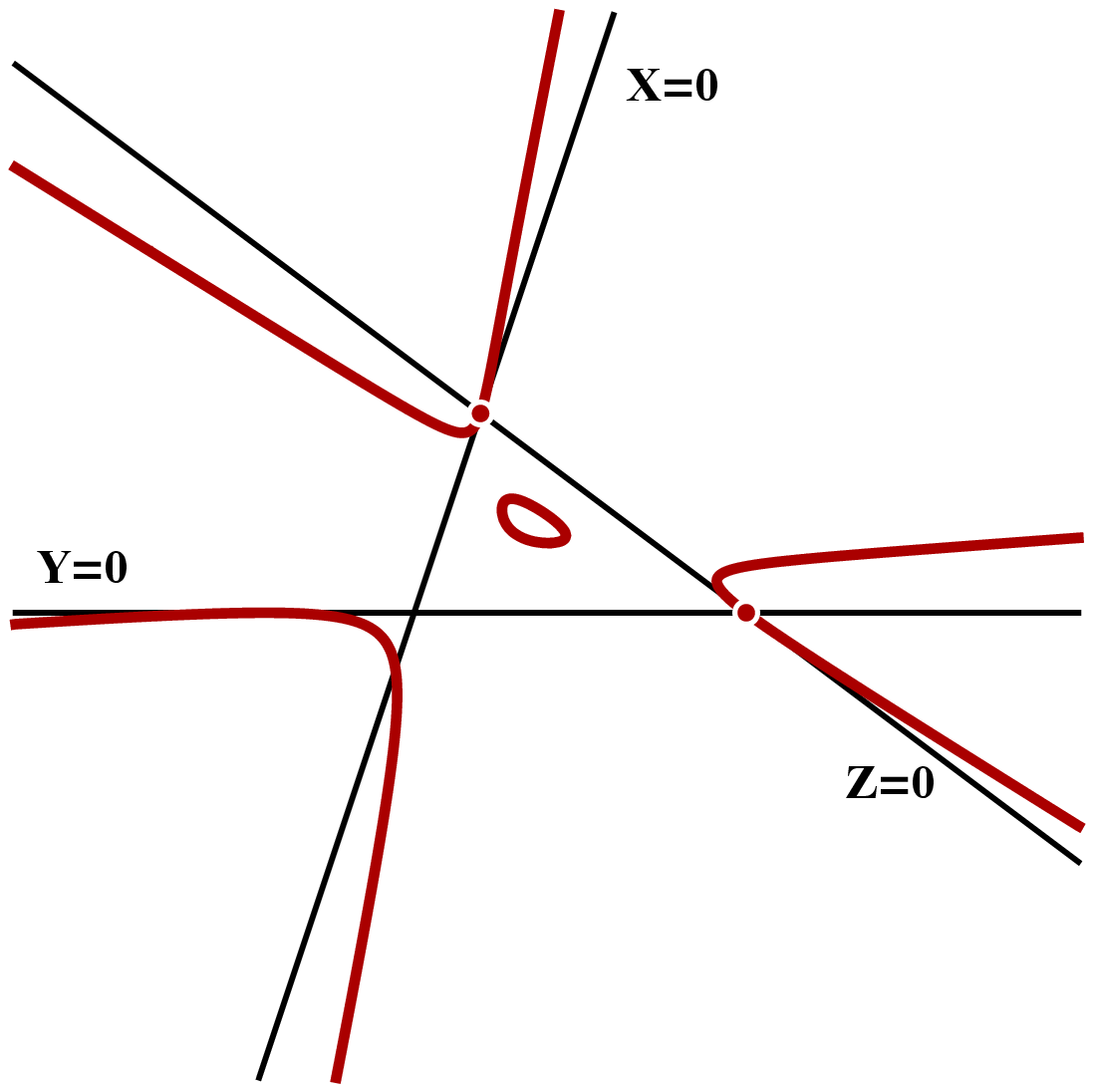}{}{{\bf After:} 
  Affine view (left) and projective view (right) of the curve $C_{(4)}$.}
  \vfill\eject 
  The curve $C_{(4)}$ has the nonzero coefficients 
  $$\eqalign{
    \Gamma_{201}^{(4)}\ &=\ 1,\cr 
    \Gamma_{120}^{(4)}\ &=\ 16(\ell +m)k^2(k+\ell )^2(k+\ell +m)^2,\cr 
    \Gamma_{111}^{(4)}\ &=\ -8k(k+\ell )(\ell +m)(k+\ell +m)(3k+2\ell +m),\cr 
    \Gamma_{102}^{(4)}\ &=\ 2(\ell +m)(3k^2 + 4k\ell + 2km - m^2), \cr 
    \Gamma_{012}^{(4)}\ &=\ -8k(k-m)(k+m)(k+\ell )(k+\ell +m)(k+2\ell +m)(\ell +m)^2,\cr 
    \Gamma_{003}^{(4)}\ &=\ (\ell +m)^2(3k^2+4k\ell +2km-m^2)^2. \cr}$$ 
  In our example this means 
  $$\eqalign{
    &\Gamma_{300}^{(4)}=0,\quad \Gamma_{210}^{(4)}=0,\quad \Gamma_{201}^{(4)}=1,
        \quad \Gamma_{120}^{(4)}=6400,\quad \Gamma_{111}^{(4)}=-2720,\cr 
    &\Gamma_{102}^{(4)}=62,\quad \Gamma_{030}^{(4)}=0,\quad \Gamma_{021}^{(4)}=0, 
        \quad \Gamma_{012}^{(4)}=43\, 680,\quad \Gamma_{003}^{(4)}=961.\cr}$$ 
  \vfill\eject 

{\bf Step 5:} Transformation to a general Weierstra\ss{} cubic. 
  \doppelbildbox{7 true cm}{C4klm_aff.eps}{C4klm_proj.eps}{{\bf Before:} 
  Affine view (left) and projective view (right) of the curve $C_{(4)}$.}
  \bigskip 
  \doppelbildbox{7 true cm}{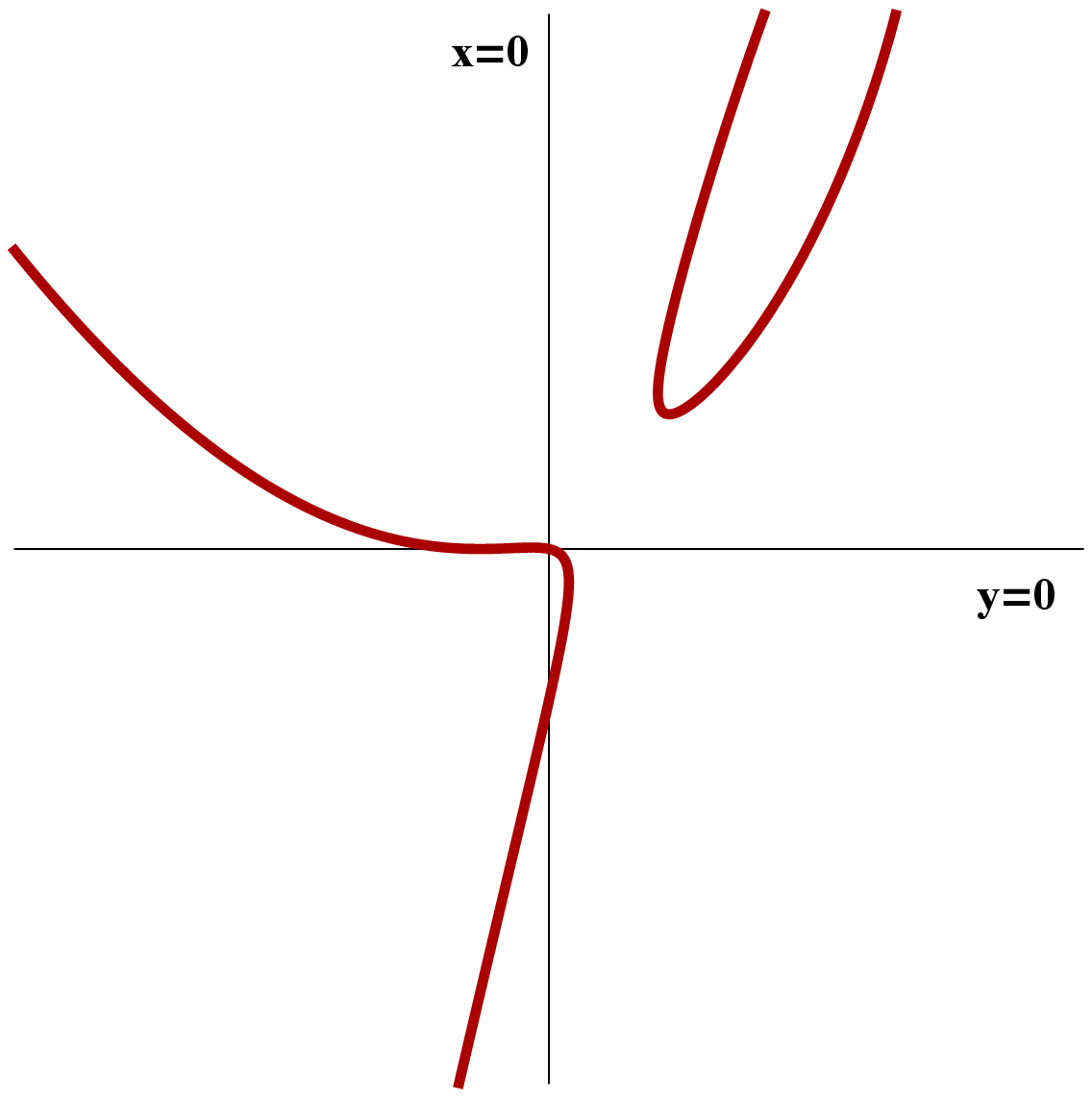}{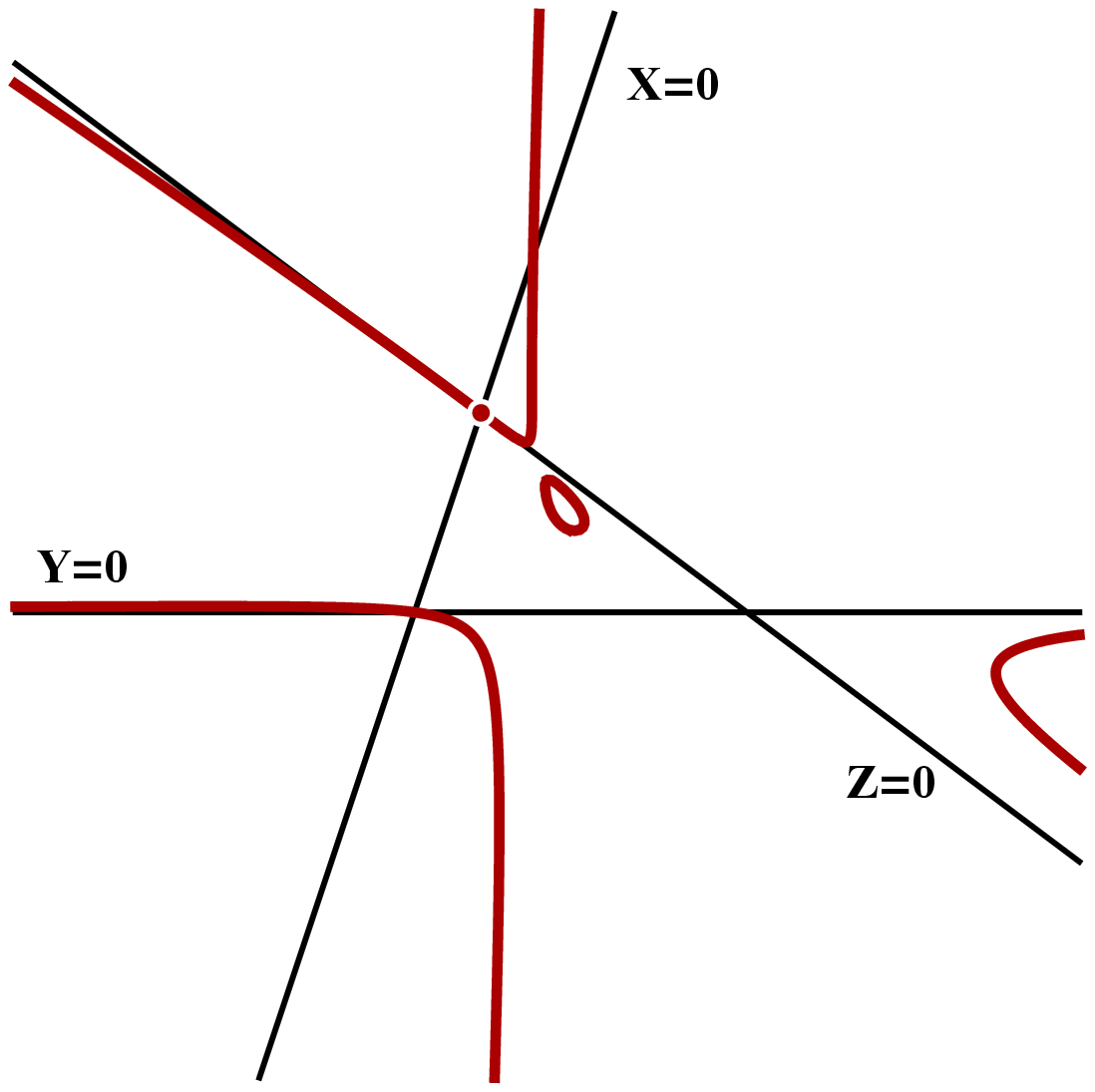}{}{{\bf After:} 
  Affine view (left) and projective view (right) of the curve $C_{(5)}$.}
  \vfill\eject 
  The curve $C_{(5)}$ has the nonzero coefficients 
  $$\eqalign{
    \Gamma_{300}^{(5)}\ &=\ \Gamma_{201}^{(4)}\ =\ 1,\cr 
    \Gamma_{201}^{(5)}\ &=\ \Gamma_{102}^{(4)}\ =\ 2(\ell+m)(3k^2+4k\ell +2km-m^2), \cr 
    \Gamma_{111}^{(5)}\ &=\ \Gamma_{111}^{(4)}\ =\ -8k(k+\ell )(\ell +m)(k+\ell +m)(3k+2\ell +m), \cr 
    \Gamma_{102}^{(5)}\ &=\ \Gamma_{003}^{(4)}\ =\ (\ell +m)^2(3k^2+4k\ell +2km-m^2)^2, \cr 
    \Gamma_{021}^{(5)}\ &=\ \Gamma_{120}^{(4)}\ =\ 16k^2(k+\ell )^2(k+\ell +m)^2(\ell +m), \cr 
    \Gamma_{012}^{(5)}\ &=\ \Gamma_{012}^{(4)}\ =\ -8k(k-m)(k+m)(k+\ell)(k+\ell+m)(k+2\ell+m)(\ell+m)^2.\cr}$$
  In our example this means 
  $$\eqalign{
    &\Gamma_{300}^{(5)}=1,\quad \Gamma_{210}^{(5)}=0,\quad \Gamma_{201}^{(5)}=62,
        \quad \Gamma_{120}^{(5)}=0,\quad \Gamma_{111}^{(5)}=-2720,\quad \Gamma_{102}^{(5)}=961,\cr 
    &\Gamma_{030}^{(5)}=0,\quad \Gamma_{021}^{(5)}=6400,\quad \Gamma_{012}^{(5)}=43680,
        \quad \Gamma_{003}^{(5)}=0.\cr}$$ 
  \vfill\eject 

{\bf Step 6:} Transformation to a special Weierstra\ss{} cubic. 
  \doppelbildbox{7 true cm}{C5klm_aff.eps}{C5klm_proj.eps}{{\bf Before:} 
  Affine view (left) and projective view (right) of the curve $C_{(5)}$.}
  \bigskip 
  \doppelbildbox{7 true cm}{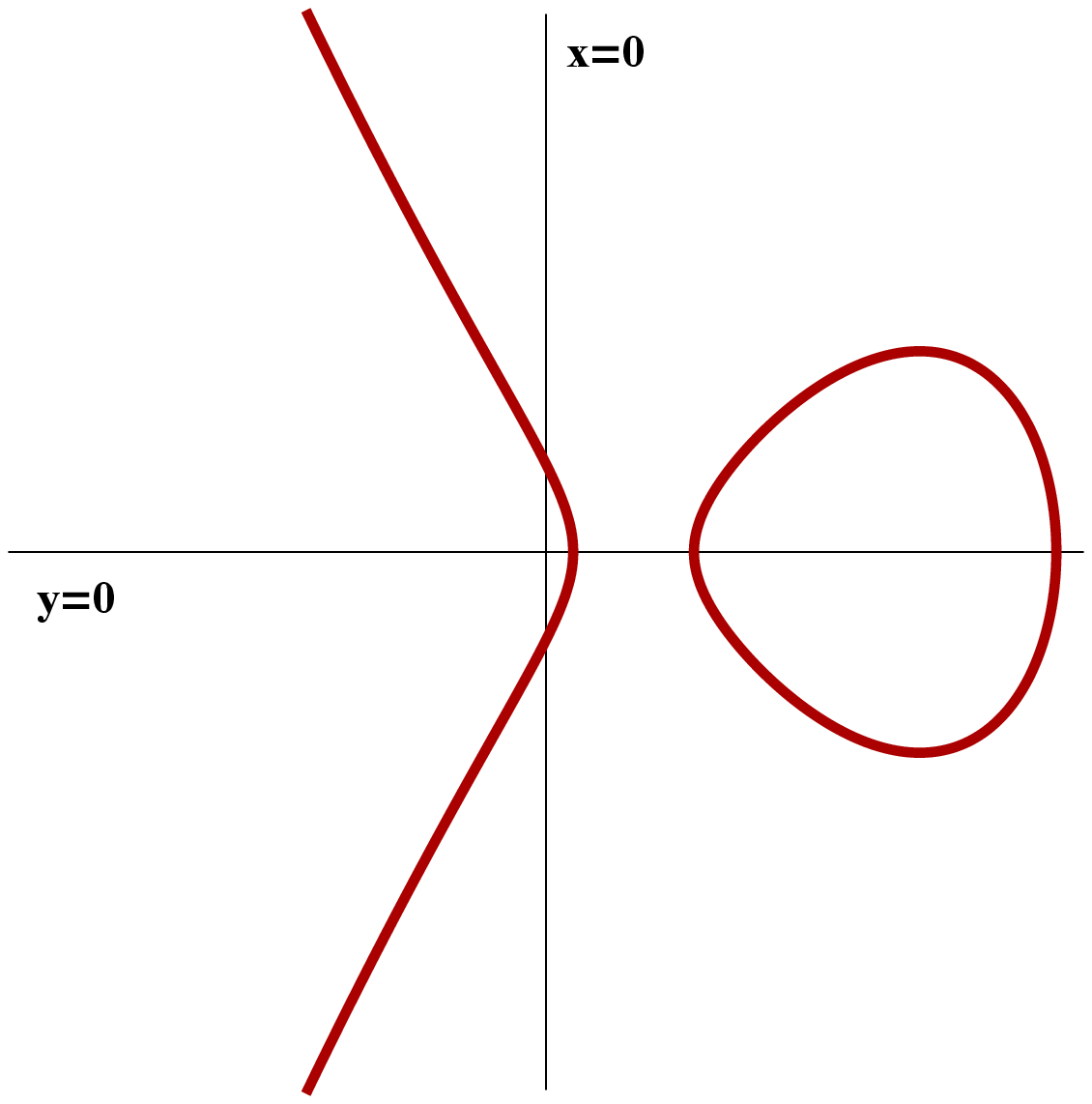}{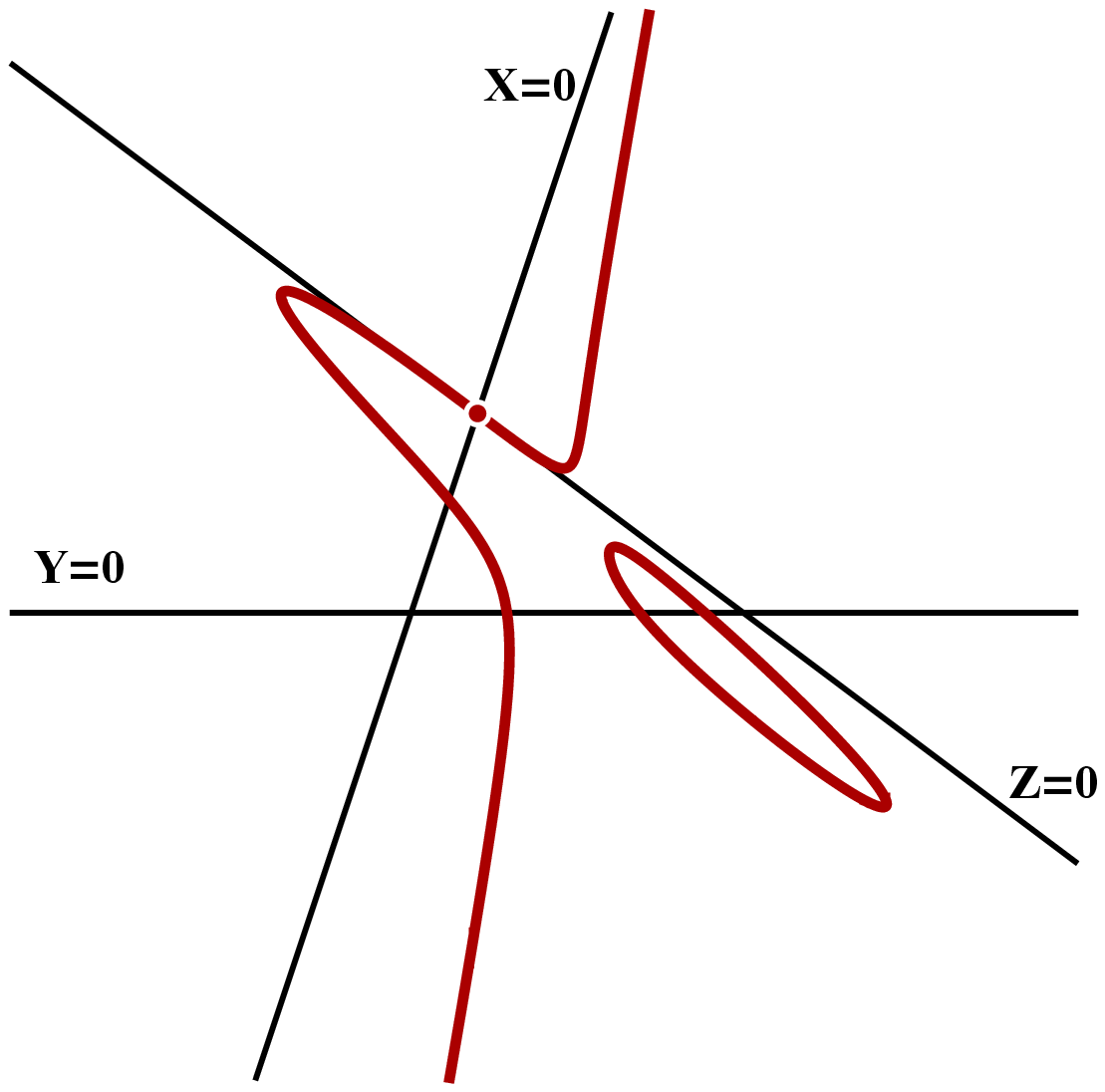}{}{{\bf After:} 
  Affine view (left) and projective view (right) of the curve $C_{(6)}$.}
  \vfill\eject 
  The curve $C_{(6)}$ has the nonzero coefficients 
  $$\eqalign{
    \Gamma_{300}^{(6)}\ &=\ 64k^2(k+\ell)^2(k+\ell +m)^2(\ell +m),\cr 
    \Gamma_{201}^{(6)}\ &=\ 64k^2(k+\ell)^2(k+\ell +m)^2(\ell +m)^2(3k^2+4k\ell +4\ell^2
        +2km+4\ell m+3m^2),\cr 
    \Gamma_{102}^{(6)}\ &=\ 64k^2(k+\ell )^2(k+\ell +m)^2(\ell +m)^3\cdot\bigl( 
         3k^4 + 8k^3\ell + 8k^2\ell^2\ \cdots \cr 
    &\phantom{Platz} \cdots\ + 4k^3m + 8k^2\ell m + 2k^2m^2 + 8k\ell m^2 + 8\ell^2m^2 + 4km^3 
         + 8\ell m^3 + 3m^4\bigr),\cr 
    \Gamma_{021}^{(6)}\ &=\ 1,\cr 
    \Gamma_{003}^{(6)}\ &=\ -64k^2(k-m)^2(k+m)^2(k+\ell )^2(k+\ell +m)^2(k+2\ell +m)^2(\ell +m)^4.\cr}$$ 
  In our example we have 
  $$\eqalign{
    &\Gamma_{300}^{(6)}=25600,\quad \Gamma_{210}^{(6)}=0,\quad \Gamma_{201}^{(6)}=-5\, 811\, 200,\cr 
    &\Gamma_{120}^{(6)}=0,\quad \Gamma_{111}^{(6)}=0,\quad \Gamma_{102}^{(6)}=262\, 220\, 800,\cr 
    &\Gamma_{030}^{(6)}=0,\quad \Gamma_{021}^{(6)}=1,\quad \Gamma_{012}^{(6)}=0,
       \quad \Gamma_{003}^{(6)}=-1\, 907\, 942\, 400.\cr}$$ 
  \vfill\eject 

{\bf Step 7:} Transformation to a Weierstra\ss{} cubic in normal form.  
  \doppelbildbox{7 true cm}{C6klm_aff.eps}{C6klm_proj.eps}{{\bf Before:} 
  Affine view (left) and projective view (right) of the curve $C_{(6)}$.}
  \bigskip 
  \doppelbildbox{7 true cm}{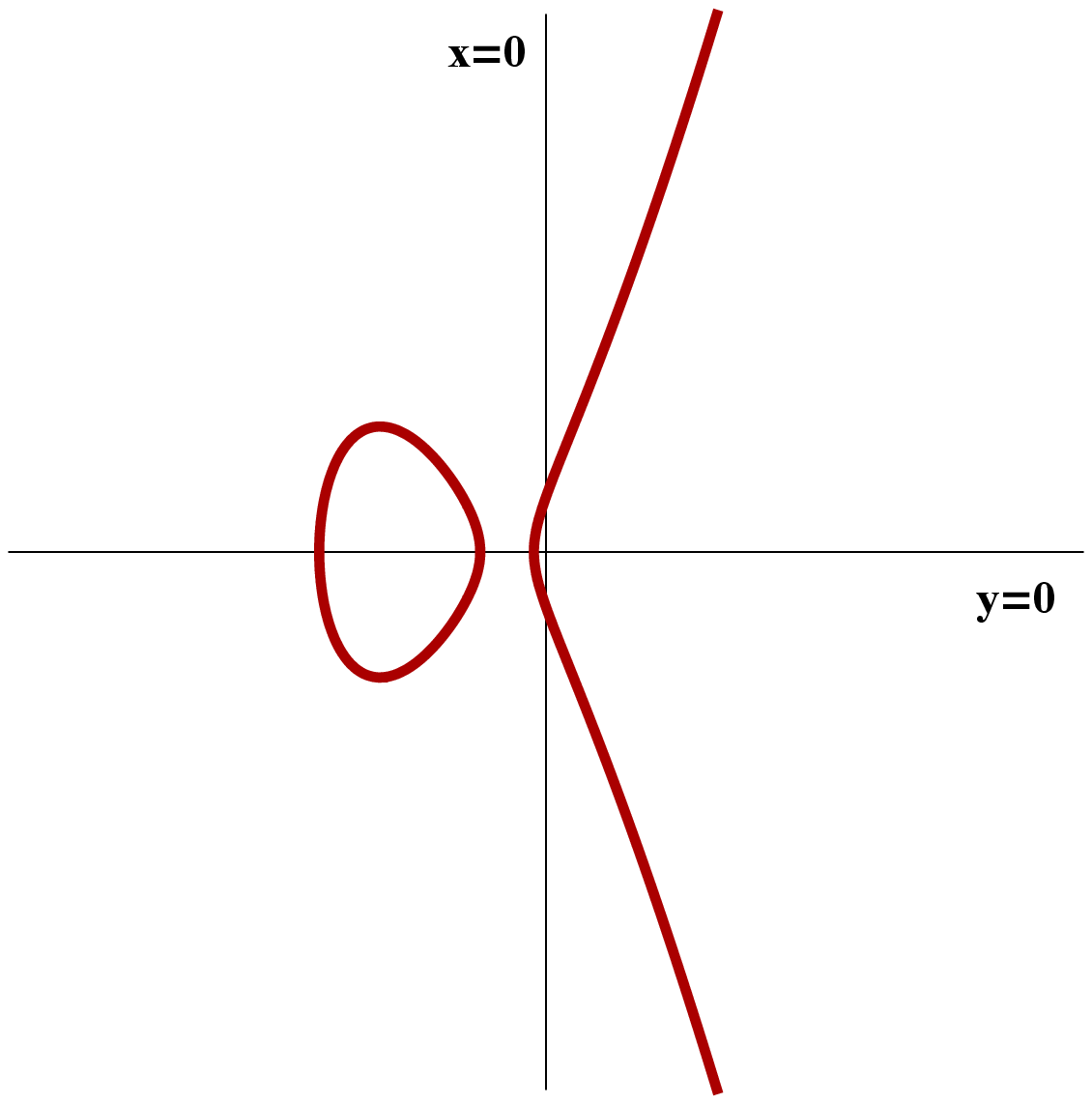}{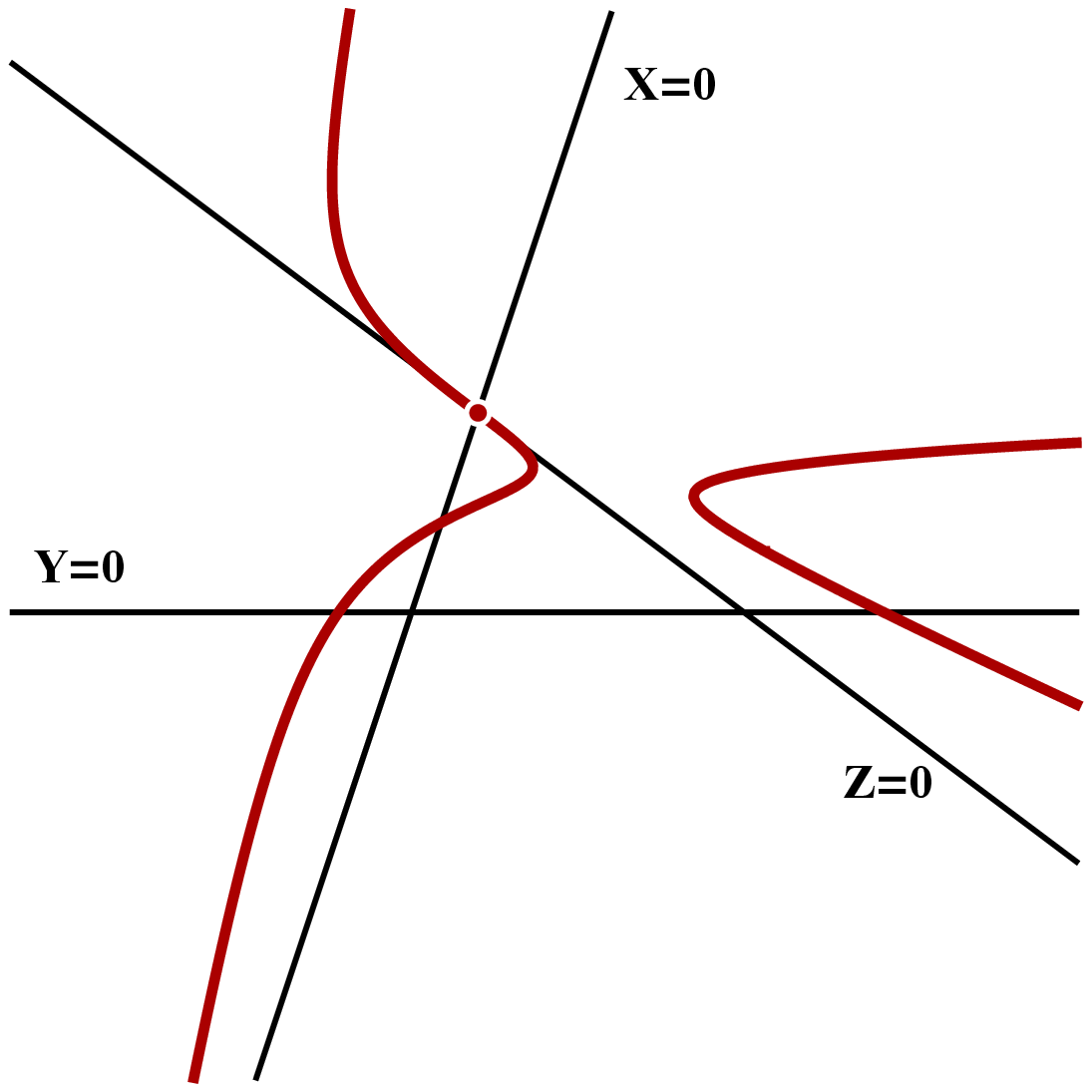}{}{{\bf After:} 
  Affine view (left) and projective view (right) of the curve $C_{(7)}$.}
  \vfill\eject 
  In our situation we have $\delta=-64(\ell +m)k^2(k+\ell)^2(k+\ell+m)^2$ and $\varphi = 
  8(\ell +m)k(k+\ell)(k+\ell +m)$, and the curve $C_{(7)}$ has the nonzero coefficients 
  $$\eqalign{
    \Gamma_{300}^{(7)}\ &=\ -1,\cr 
    \Gamma_{201}^{(7)}\ &=\ -(3k^2 + 4k\ell + 4\ell^2 + 2km + 4\ell m + 3m^2), \cr 
    \Gamma_{102}^{(7)}\ &=\ -\bigl( 3k^4 + 8k^3\ell + 8k^2\ell^2 + 4k^3m + 8k^2\ell m + 2k^2m^2 
            + 8k\ell m^2\ \cdots\cr 
    &\phantom{Platz} \cdots\ + 8\ell^2m^2 + 4km^3 + 8\ell m^3 + 3m^4\bigr), \cr 
    \Gamma_{021}^{(7)}\ &=\ 1,\cr 
    \Gamma_{003}^{(7)}\ &=\ -(k-m)^2(k+m)^2(k+2\ell +m)^2.\cr}$$
  In our example this means 
  $$\eqalign{
    &\Gamma_{300}^{(7)}=-1,\quad \Gamma_{210}^{(7)}=0,\quad \Gamma_{201}^{(7)}=-227, 
      \quad \Gamma_{120}^{(7)}=0,\quad \Gamma_{111}^{(7)}=0,\quad \Gamma_{102}^{(7)}=-10243,\cr 
    &\Gamma_{030}^{(7)}=0,\quad \Gamma_{021}^{(7)}=1,\quad \Gamma_{012}^{(7)}=0, 
      \quad \Gamma_{003}^{(7)}=-74529.\cr}$$ 
  \vfill\eject 

{\bf Step 8.} In our special situation we can do more than in the general case, since the
  cubic polynomial in $X_7$ and $Z_7$ splits over the rationals. One easily computes the roots 
  of this polynomial to be $-(k-m)^2$, $-(k+m)^2$ and $-(k+2\ell +m)^2$. The simple substitution 
  $x=\widehat{x}-(k-m)^2$ and a further simple substitution to eliminate a common factor $4$ in 
  the roots yields the equation 
  $$\widehat{y}^2\ =\ \widehat{x}(\widehat{x}+km)\bigl(\widehat{x}+(k+\ell)(\ell+m)\bigr).$$ 
    \doppelbildbox{6.9 true cm}{C7klm_aff.eps}{C7klm_proj.eps}{{\bf Before:} 
    Affine view (left) and projective view (right) of the curve $C_{(7)}$.}
    \bigskip 
    \doppelbildbox{6.9 true cm}{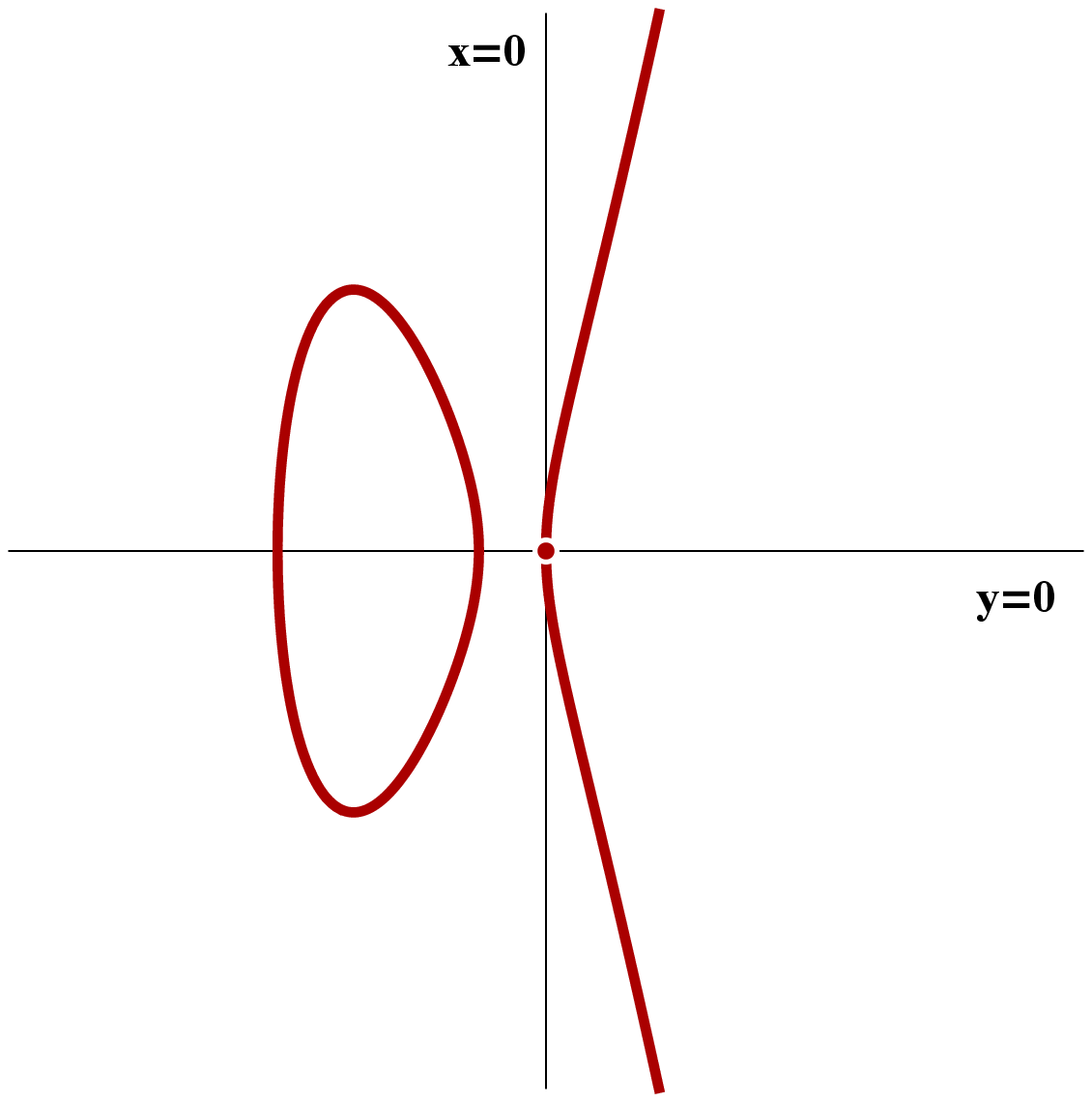}{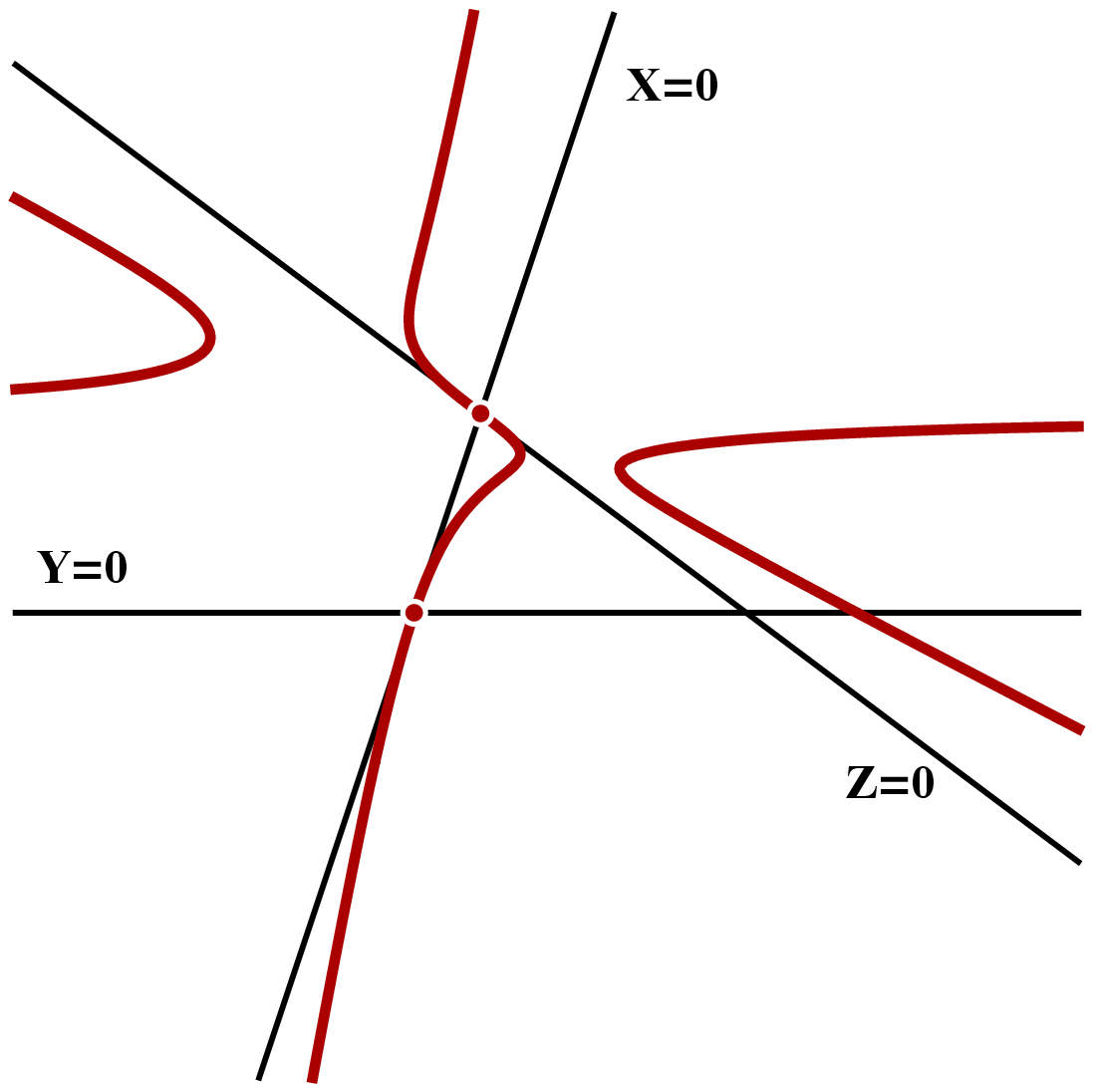}{}{{\bf After:} 
    Affine view (left) and projective view (right) of the curve $C_{(8)}$.}
  \eject 
  The equation obtained in the previous step has the affine form $y^2 = p(x)$ where 
  $p$ has the roots $-(k-m)^2$, $-(k+m)^2$ and $-(k+2\ell +m)^2$. The simple substitution $x = 
  \widetilde{x}-(k\! -\! m)^2$ and a further simple substitution to eliminate the common factor 
  $4$ in the roots yields the equation $\widetilde{y}^{2} = \widetilde{x}(\widetilde{x}+km)
  \bigl(\widetilde{x}+(k\! +\!\ell )(\ell +m))$. Thus we see that the problem of finding 
  rational squares in an arithmetic progression of type $(k,\ell,m)$ is equivalent to finding 
  rational points on the curve $E_{k,\ell ,m}$ given by the Weierstra\ss{} equation 
  $$y^2\ =\ x(x+km)\bigl( x+(k+\ell )(\ell +m)\bigr)\, .\leqno{(\star )}$$ 
  In projective form, this curve is given by the homogeneous equation 
  $$Y^2Z\ =\ X(X+kmZ)\bigl( X+(k+\ell )(\ell +m)Z\bigr)$$ 
  where $x=X/Z$ and $y=Y/Z$. Thus in our example we end up with the coefficients 
  $$\eqalign{
    &\Gamma_{300}^{(8)}=-1,\quad \Gamma_{210}^{(8)}=0,\quad \Gamma_{201}^{(8)}=-50,\quad 
        \Gamma_{120}^{(8)}=0, \quad \Gamma_{111}^{(8)}=0,\cr 
    &\Gamma_{102}^{(8)}=-400,\quad \Gamma_{030}^{(8)}=0,\quad \Gamma_{021}^{(8)}=1,
        \quad \Gamma_{012}^{(8)}=0,\quad \Gamma_{003}^{(8)}=0.\cr}$$ 
From $(\star )$ one recognizes the remarkable fact that the torsion group of the resulting curve 
contains ${\Z}_2\times{\Z}_2$ as a subgroup. Moreover, the final coefficients are surprisingly 
small, due to the cancellation of the terms which was possible in Step 7 (with rather large 
factors $\delta$ and $\varphi$). 
 \vfill\eject 
The composition of the above transformations results in an overall transformation 
$\Phi:{\P}^3\rightarrow{\P}^2$, say $\Phi(x_0,x_1,x_2,x_3)=(X,Y,Z)$, which induces 
a birational (thus {\it a fortiori} biregular) mapping of the quadric intersection 
$Q_{k,\ell,m}$ to the elliptic curve $E_{k,\ell,m}$. This transformation $\Phi$ is 
represented by a quadratic transformation, defined by 
$$X = \sum_{i,j} X_{ij}x_ix_j, \quad 
  Y = \sum_{i,j} Y_{ij}x_ix_j, \quad   
  Z = \sum_{i,j} Z_{ij}x_ix_j$$ 
where the sums are formed over all indices $0\leq i,j\leq 3$ such that $i\leq j$ and 
where the coefficients are given by  
$$\eqalign{
  X_{00}\ &=\ -km^2(k+\ell)^2(\ell +m)(k+\ell +m)^2,\cr 
  X_{01}\ &=\ km(k+\ell )(\ell +m)(k+\ell +m)^2(k\ell +2km+\ell m),\cr 
  X_{02}\ &=\ -k\ell m^2(k+\ell )(k-m)(\ell +m)(k+\ell +m),\cr 
  X_{03}\ &=\ -k\ell m(k-m)(k+\ell)^2(\ell +m)(k+\ell +m),\cr 
  X_{11}\ &=\ -k^2m(k+\ell)(\ell +m)^2(k+\ell +m)^2,\cr 
  X_{12}\ &=\ k\ell m(k+\ell )(k-m)(\ell +m)^2(k+\ell +m),\cr 
  X_{13}\ &=\ k^2m\ell (k+\ell )(k-m)(\ell +m)(k+\ell +m),\cr 
  X_{22}\ &=\ k\ell^2m^2(k+\ell )(\ell +m)^2,\cr 
  X_{23}\ &=\ -k\ell^2m(k+\ell )(\ell +m)(k^2+k\ell +\ell m+m^2),\cr 
  X_{33}\ &=\ k^2\ell^2m(k+\ell )^2(\ell +m),\cr}$$ 
by 
$$\eqalign{
  Y_{00}\ &=\ k\ell m^2(k+\ell )^2(\ell +m)(k+\ell +m)^2, \cr 
  Y_{01}\ &=\ k\ell^2m(k+\ell )(k-m)(\ell +m)(k+\ell +m)^2, \cr 
  Y_{02}\ &=\ -k^2\ell m^2(k+\ell )(\ell +m)(k+\ell +m)(k+2\ell +m), \cr 
  Y_{03}\ &=\ -k\ell m(k+\ell )^2(k+m)(\ell +m)^2(k+\ell +m), \cr 
  Y_{11}\ &=\ -k^2\ell m(k+\ell )(\ell +m)^2(k+\ell +m)^2, \cr 
  Y_{12}\ &=\ k\ell m(k+\ell )^2(k+m)(\ell +m)^2(k+\ell +m), \cr 
  Y_{13}\ &=\ k^2\ell m^2(k+\ell )(\ell +m)(k+\ell +m)(k+2\ell +m), \cr 
  Y_{22}\ &=\ -k\ell^2m^2(k+\ell)(\ell +m)^2(k+\ell +m), \cr 
  Y_{23}\ &=\ -k\ell^2m(k+\ell )(k-m)(\ell +m)(k+\ell +m)^2, \cr 
  Y_{33}\ &=\ k^2\ell^2m(k+\ell )^2(\ell +m)(k+\ell +m) \cr}$$ 
and by 
$$\eqalign{
  Z_{00}\ &=\ m^2(k+\ell)^2(k+\ell +m)^2, \cr 
  Z_{01}\ &=\ -2km(k+\ell)(\ell +m)(k+\ell +m)^2, \cr  
  Z_{02}\ &=\ -2\ell m^2(k+\ell)(\ell +m)(k+\ell +m), \cr 
  Z_{03}\ &=\ 2k\ell m(k+\ell )^2(k+\ell +m), \cr 
  Z_{11}\ &=\ k^2(\ell +m)^2(k+\ell +m)^2, \cr 
  Z_{12}\ &=\ 2k\ell m(\ell +m)^2(k+\ell +m), \cr 
  Z_{13}\ &=\ -2k^2\ell (k+\ell )(\ell +m)(k+\ell +m), \cr 
  Z_{22}\ &=\ \ell^2m^2(\ell +m)^2, \cr 
  Z_{23}\ &=\ -2k\ell^2m(k+\ell )(\ell +m), \cr 
  Z_{33}\ &=\ k^2\ell^2(k+\ell)^2.\cr}$$ 
The original quadric $Q_{k,\ell,m}$ possesses eight trivial points with
coordinates $(1,\pm 1,\pm 1,\pm 1)$. These are mapped to eight rational
points on the elliptic curve $E_{k,\ell,m}$ via the transformation $\Phi$.
These points are given in projective form $(X,Y,T)$ with integer 
coordinates as follows: 
\def\phm{\phantom{-}}
$$\eqalign{
  \Phi(1,{\phm}1,{\phm}1,{\phm}1)\ &=\ (0,1,0)\quad \hbox{(point at infinity)}, \cr 
  \Phi(1,{\phm}1,{\phm}1,-1)\ &=\ \bigl( m(\ell+m),\, m(\ell +m)(k+\ell +m),\, 1\bigr), \cr 
  \Phi(1,{\phm}1,-1,{\phm}1)\ &=\ \bigl( k(k+\ell),\, -k(k+\ell)(k+\ell +m),\, 1\bigr), \cr 
  \Phi(1,{\phm}1,-1,-1)\ &=\ (0,0,1)\quad\hbox{(2-torsion point)},\cr 
  \Phi(1,-1,{\phm}1,{\phm}1)\ &=\ \bigl( -m(k+\ell),\, -m\ell (k+\ell),\, 1\bigr), \cr 
  \Phi(1,-1,{\phm}1,-1)\ &=\ \bigl( -(k+\ell)(\ell +m),\, 0,\, 1)\quad\hbox{(2-torsion point)}, \cr 
  \Phi(1,-1,-1,{\phm}1)\ &=\ (-km,0,1)\quad\hbox{(2-torsion point)},\cr 
  \Phi(1,-1,-1,-1)\ &=\ \bigl( -k(\ell +m),\, k\ell (\ell +m),\, 1\bigr). \cr}$$ 
\vfill\eject 

\leftline{\titlefont 6. References:} \bigskip 
\item{[1]} John William Scott Cassels, {\it Lectures on Elliptic Curves}; 
  Cambridge University Press, Cambridge 1991 \smallskip\noindent  
\item{[2]} Ian Connell, {\it Elliptic Curve Handbook}; 
  http://webs.ucm.es/BUCM/mat/doc8354.pdf\smallskip\noindent 
\item{[3]} Leonhard Euler, {\it De binis formulis speciei xx+myy et xx+nyy inter se 
  concordibus et discordibus}, Mem. Acad. Sci. St.-Petersbourg 1780 (Opera Omnia: 
  Ser. 1, Vol. 5, pp. 406-413) \smallskip\noindent 
\item{[4]} Robin C. Hartshorne, {\it Algebraic Geometry}, Springer, 1977
  \smallskip\noindent 
\item{[5]} Neal Koblitz, {\it Introduction to Elliptic Curves and Modular Forms}, 
  Springer, New York/Berlin/Heidelberg 1984  \smallskip\noindent 
\item{[6]} Trygve Nagell, {\it Sur les propri\'e{}t\'e{}s arithm\'e{}tiques 
  des cubiques planes du premier genre}, Acta Mathematica {\bf 52}, 1929, 
  pp- 93-126 \smallskip\noindent 
\item{[7]} Ken Ono, {\it Euler's Concordant Forms}, Acta arithmetica LXXVIII (2), 
  1996, pp. 101-123 \smallskip\noindent 
\item{[8]} Takashi Ono, {\it Variations on a Theme of Euler}, Plenum Press, 
  New York/London 1994 \smallskip\noindent 
\item{[9]} John G. Semple, Geoffrey T. Kneebone, {\it Algebraic Projective Geometry}, 
  Clarendon Press, Oxford 1952\smallskip\noindent  
\item{[10]} Igor R. Shafarevich, {\it Basic Algebraic Geometry 1}, Springer, 
  Berlin/Heidelberg 1977 \smallskip\noindent 
\item{[11]} Joseph H. Silverman, {\it The Arithmetic of Elliptic Curves}, Springer, 
  Dordrecht$\,$/$\,$Heidel\-berg$\,$/$\,$London$\,$/$\,$New York 2009  \smallskip\noindent 
\item{[12]} Joseph H. Silverman, John Tate, {\it Rational Points on Elliptic Curves}, 
  Springer, New York, 2nd edition, 2015
\par\hfill\break\smallskip 

\noindent 
Hagen Knaf, Karlheinz Spindler \hfill\break
Hochschule RheinMain, Germany \hfill\break
Applied Mathematics \hfill\break
{\tt hagen.knaf@hs-rm.de}, {\tt karlheinz.spindler@hs-rm.de} 
\par\bigskip\noindent 
Erich Selder \hfill\break
Frankfurt University of Applied Sciences, Germany \hfill\break
Computer Science and Engineering \hfill\break 
{\tt e\_selder@fb2.fra-uas.de} 
\bye

%% file: colorbox.tex

\newdimen\colorboxtextdistance
\colorboxtextdistance=10pt

%% file: bfall.tex

\font\tenbi=cmmib10      \font\tenbsy=cmbsy10
\font\sevenbi=cmmib7     \font\sevenbsy=cmbsy7
\font\fivebi=cmmib5      \font\fivebsy=cmbsy5
\font\tenbit=cmbxti10    \font\tenbsl=cmbxsl10

\def\boldmath{
  \textfont0=\tenbf   \scriptfont0=\sevenbf
     \scriptscriptfont0=\fivebf
  \textfont1=\tenbi   \scriptfont1=\sevenbi
     \scriptscriptfont1=\fivebi
  \textfont2=\tenbsy  \scriptfont2=\sevenbsy
     \scriptscriptfont2=\fivebsy
  \textfont\itfam=\tenbit  \textfont\slfam=\tenbsl}

\def\bfall{\boldmath \let\tenrm\tenbf \let\tenit\tenbit
  \let\tensl\tenbsl \rm}